\pdfoutput=1
\RequirePackage{ifpdf}
\ifpdf 
\documentclass[pdftex]{sigma}
\else
\documentclass{sigma}
\fi

\usepackage{mathrsfs,amsmath,fancyhdr,amsthm,tikz-cd}
\usepackage{cancel,rotating}
\usepackage{amsbsy}

\usepackage{mathrsfs,fancyhdr,amsthm,tikz-cd}
\usepackage[all]{xy}

\numberwithin{equation}{section}

\newtheorem{Theorem}{Theorem}[section]
\newtheorem*{Theorem*}{Theorem}
\newtheorem{Corollary}[Theorem]{Corollary}
\newtheorem{Lemma}[Theorem]{Lemma}
\newtheorem{Proposition}[Theorem]{Proposition}
 { \theoremstyle{definition}
\newtheorem{Definition}[Theorem]{Definition}

\newtheorem{Remark}[Theorem]{Remark} }

\newcommand*\Bell{\ensuremath{\boldsymbol\ell}}

 \def\leq{\leqslant} \def\geq{\geqslant}

\def\Hom{\operatorname{Hom}} 
 \def\im{\operatorname{Im}}

\def\Ker{\operatorname{Ker}}
\def\dim{\operatorname{dim}} \def\End{\operatorname{End}}

\def\mod{\operatorname{mod}}  
 
\def\e{{\bf e}}
\def\f{{\bf f}}
\def\g{{\bf d}}
\def\l{{\bf l}}

\def\m{{\bf m}}
\def\n{{\bf n}}

\def\v{{\bf v}}
\def\j{{\bf j}}

\def\w{{\bf w}}
\def\k{{\bf k}}

\newcommand \MF{\operatorname{MF}}
\newcommand \MG{\operatorname{MG}}
\newcommand \EExt{\operatorname{Ext}}
\newcommand \Ext{\operatorname{Ext}}
\newcommand \KKer{\operatorname{Ker}}
\newcommand \EF{\operatorname{EF}}

\begin{document}
\allowdisplaybreaks

\newcommand{\arXivNumber}{2110.12429}

\renewcommand{\PaperNumber}{097}

\FirstPageHeading

\ShortArticleName{The Multiplication Formulas of Weighted Quantum Cluster Functions}

\ArticleName{The Multiplication Formulas of Weighted Quantum\\ Cluster Functions}

\Author{Zhimin CHEN~$^{\rm a}$, Jie XIAO~$^{\rm b}$ and Fan XU~$^{\rm c}$}

\AuthorNameForHeading{Z.~Chen, J.~Xiao and F.~Xu}

\Address{$^{\rm a)}$~Department of Mathematics, Tsinghua University, Beijing 100084, P.~R.~China}
\EmailD{\href{mailto:chen-zm15@mails.tsinghua.edu.cn}{chen-zm15@mails.tsinghua.edu.cn}}

\Address{$^{\rm b)}$~School of Mathematical Sciences, Beijing Normal University, Beijing 100875, P.~R.~China}
\EmailD{\href{mailto:jxiao@bnu.edu.cn}{jxiao@bnu.edu.cn}}

\Address{$^{\rm c)}$~Department of Mathematical Sciences, Tsinghua University, Beijing 100084, P.~R.~China}
\EmailD{\href{mailto:fanxu@mail.tsinghua.edu.cn}{fanxu@mail.tsinghua.edu.cn}}

\ArticleDates{Received June 02, 2022, in final form November 26, 2023; Published online December 13, 2023}

\Abstract{By applying the property of Ext-symmetry and the affine space structure of certain fibers, we introduce the notion of weighted quantum cluster functions and prove their multiplication formulas associated to abelian categories with $\operatorname{Ext}$-symmetry and 2-Calabi--Yau triangulated categories with cluster-tilting objects.}


\Classification{17B37; 16G20; 17B20}

\section{Introduction, notation and main results}\label{sec1}
Cluster algebras were introduced by Fomin and Zelevinsky \cite{FZ} in order to find an algebraic framework for understanding total positivity in algebraic groups and canonical bases in quantum groups. They are commutative algebras over $\mathbb{Z}$ generated by certain elements called cluster variables. There is an iterative procedure called mutation to obtain new cluster variables from initial cluster variables. For an cluster algebra of finite type, the cluster variables are in bijection with the almost positive roots of the corresponding simple Lie algebra.

Cluster algebras have close connections to representation theory of algebras via cluster categories and cluster characters. The cluster category was constructed by Buan--Marsh--Reineke--Reiten--Todorov \cite{BMRRT} as a quotient of the bounded derived category of the module category of a~finite-dimensional hereditary algebra. Given an acyclic quiver $Q$, the cluster category $\mathcal{C}_Q$ associated with $Q$ is the orbit category $\mathcal{D}^b(Q)/\tau\circ [-1].$ Indecomposable objects in $\mathcal{C}_Q$ correspond to the almost positive roots in the root system of the Lie algebra $\mathfrak{g}_Q$ of type $Q$. The cluster category has the 2-Calabi--Yau property, i.e.,
\[
\operatorname{Hom}_{\mathcal{C}_Q}(M, N[1])\cong D\operatorname{Hom}_{\mathcal{C}_Q}(N, M[1])
\] for any $M, N\in \mathcal{C}_Q$.

For a Dynkin quiver $Q$ with $n$ vertices, Caldero and Chapoton \cite{CC2006} defined a map $X_{?}$ from the set of objects in the cluster category $\mathcal{C}_Q$ to $\mathbb{Q}(x_1, \dots, x_n)$ such that
\[
 X_M=\sum_{\underline{e}}\chi(\mathrm{Gr}_{\underline{e}}(M))x^{-B\underline{e}-(I_n-R)\underline{m}}
\]
 for a $\mathbb{C}Q$-module $M$ and $X_{P[1]}=x^{\underline{\dim }P/\mathrm{rad}P}$ for a projective $\mathbb{C}Q$-module $P$ (See \cite{CK2005} for the notation in the above formula.)
 In particular, it gives a bijection between the indecomposable rigid objects and cluster variables and then the mutation rule of cluster variables is written as the following property. For indecomposable rigid objects $M$, $N$ in $\mathcal{C}_Q$ with ${\dim _{\mathbb{C}}\operatorname{Hom}_{\mathcal{C}_Q}(M, N[1])=1}$, we have $X_M\cdot X_N=X_{L}+X_{L'}$, where $L$, $L'$ are the middle terms of two non-split triangles induced by $\varepsilon\in \operatorname{Hom}_{\mathcal{C}_Q}(M, N[1])$ and $\varepsilon'\in \operatorname{Hom}_{\mathcal{C}_Q}(N, M[1])$, respectively. This map is called the cluster character. Then the cluster character $X_M$ for an indecomposable rigid object $M$ can be viewed as the general formula of a~cluster variable, computed from the initial variables $X_{P_1[1]}, \dots, X_{P_n[1]}$ using the mutation rule. Caldero and Keller \cite{CK2006} generalized these results to any acyclic quiver and then categorified all cluster algebras with acyclic initial clusters. As a~corollary, the cluster algebra $\mathcal{A}(Q)$ is generated by $X_{M}$ for all rigid objects $M\in \mathcal{C}_Q$. Palu \cite{Palu2008} generalized the definition of cluster character to a 2-Calabi--Yau triangulated category $\mathcal{C}$ with a~cluster tilting object~$T$. His cluster character, denoted by $X^T_?$, is given by the formula
\[
X^T_M=\sum_{\underline{e}}\chi(\mathrm{Gr}_{\underline{e}}(\operatorname{Hom}_{\mathcal{C}}(T, M)))x^{\operatorname{ind}_{T}M-\iota(\underline{e})}
\]
for $M\in \mathcal{C}$ (see Section \ref{index} for more details).

In \cite{CK2005}, Caldero and Keller proved a higher-dimensional multiplication formula between cluster characters for the cluster category $\mathcal{C}_Q$ of a Dynkin quiver $Q$. They showed that for any indecomposable objects $M, N\in \mathcal{C}_Q$,
\begin{align*}
& \chi(\mathbb{P}\operatorname{Hom}_{\mathcal{C}_Q}(M, N[1]))X_M\cdot X_N\\
&\quad=\sum_{[L]}(\chi(\mathbb{P}\operatorname{Hom}_{\mathcal{C}_Q}(M, N[1])_{L[1]})+\chi(\mathbb{P}\operatorname{Hom}_{\mathcal{C}_Q}(N, M[1])_{L[1]}))X_L,
\end{align*}
where $\chi$ is the Euler--Poincar\'{e} characteristic. The proof heavily depends on the 2-Calabi--Yau property of $\mathcal{C}_Q$. This multiplication formula was generalized firstly to acyclic quivers \cite{XiaoXu2007, Xu2010} and then to any 2-Calabi--Yau triangulated category with a cluster-tilting object \cite{Palu2012}.

Inspired by the results in \cite{CK2005}, Geiss, Leclerc and Schr\"{o}er \cite[Theorem 3]{GLS2007} proved an analogous formula for nilpotent module categories of preprojective algebras. Let $\Lambda_Q$ be the preprojective algebra of an acyclic quiver $Q$ and $\operatorname{nil}(\Lambda_Q)$ be the category of finite dimensional nilpotent $\Lambda_Q$-modules. Let $\mathfrak{n}$ be the positive part of $\mathfrak{g}_Q$. Given two sequences $\mathbf{i}=(i_1,\dots,i_m)$, $\mathbf{a}=(a_1,\dots,a_m)$ where $i_j\in\{1,\dots,n\}$ and $a_j\in \{0,1\}$ for $1 \leq j \leq m$, Lusztig introduced the constructible functions $d_{\mathbf{i}, \mathbf{a}}\colon \operatorname{nil}(\Lambda_Q)\rightarrow \mathbb{C}$ defined by
\[
d_{\mathbf{i}, \mathbf{a}}(M)=\chi(\mathcal{F}_{\mathbf{i}, \mathbf{a}}(M)),
\]
where $\mathcal{F}_{\mathbf{i}, \mathbf{a}}(M)$ is the set of flags of type $(\mathbf{i}, \mathbf{a})$ (see Section \ref{type} for the definition). He proved that the enveloping algebra $U(\mathfrak{n})$ is isomorphic to $\mathcal{M}=\bigoplus \mathcal{M}_{\underline{d}}$, where the $\mathcal{M}_{\underline{d}}$ are certain vector spaces generated by $d_{\mathbf{i}, \mathbf{a}}$ for proper pairs $(\mathbf{i}, \mathbf{a}).$ Geiss, Leclerc and Schr\"{o}er showed that the category $\operatorname{nil}(\Lambda_Q)$ has the Ext-symmetry, i.e.,
\[\mathrm{Ext}^1_{\Lambda_Q}(M, N)\cong D\mathrm{Ext}^1_{\Lambda_Q}(N, M)\]
for $M, N\in \operatorname{nil}(\Lambda_Q)$. Let $\mathcal{M}^*$ be the graded dual of $\mathcal{M}$. For $M\in \operatorname{nil}(\Lambda_Q)$, they defined an evaluation form $\delta_M\in \mathcal{M}^*$, which is the map from $\mathcal{M}$ to $\mathbb{Z}$ mapping $d_{\mathbf{i}, \mathbf{a}}$ to $d_{\mathbf{i}, \mathbf{a}}(M)$ for any pair $(\mathbf{i}, \mathbf{a})$. We denote by $\Lambda_{\underline{d}}$ the variety of nilpotent $\Lambda$-modules of dimension vector $\underline{d}$ and set $\langle M \rangle:=\{X\in \Lambda_{\underline{d}}\mid \delta_X=\delta_M\}$. Then there exists a finite subset $S(\underline{d})$ of $\Lambda_{\underline{d}}$ such that $\Lambda_{\underline{d}}=\bigsqcup_{M\in S(\underline{d})}\langle M\rangle$. For $M\in \Lambda_{\underline{d}}$ and $N\in \Lambda_{\underline{d}'}$, they obtained the following multiplication formula:
\[
\chi\bigl(\mathbb{P}\mathrm{Ext}^1_{\Lambda_Q}(M, N)\bigr)\delta_M\cdot \delta_N=\sum_{L\in S(\underline{d}) }\bigl(\chi\bigl(\mathbb{P}\mathrm{Ext}^1_{\Lambda_Q}(M, N)_{\langle L \rangle}\bigr)+\chi\bigl(\mathbb{P}\mathrm{Ext}^1_{\Lambda_Q}(N, M)_{\langle L\rangle}\bigr)\bigr)\delta_L.
\]

The similarity between $\mathcal{C}_Q$ and $\operatorname{nil}(\Lambda_Q)$ was further studied by many authors~\cite{BIRS, GLS2011}. Let~$W$ be the Weyl group of $\mathfrak{g}_Q$. For $w\in W,$ Buan, Iyama, Reiten, and Scott \cite{BIRS} constructed a~2-Calabi--Yau Frobenius subcategory $\mathcal{C}_w$ of $\operatorname{nil}(\Lambda_Q)$. When $w=(s_n\cdots s_1)^2$, the stable category $\underline{\mathcal{C}_w}$ of $\mathcal{C}_w$ is just $\mathcal{C}_Q.$ Given a reduced expression of $w=s_{i_r}\cdots s_{i_1}$ and set $\mathbf{i}:=(i_r, \dots, i_1)$. For~$M\in \mathcal{C}_w$,~$\delta_M$ can be reformulated into the following form:
\[
\Delta_{M}=\sum_{\mathbf{a}\in \mathbb{N}^r}\delta_M(d_{\mathbf{i}, \mathbf{a}})t^{\mathbf{a}}=\sum_{\mathbf{a}\in \mathbb{N}^r}\chi(\mathcal{F}_{\mathbf{i},\mathbf{a}}{(M)})t^{\mathbf{a}}.
\]
Geiss, Leclerc and Schr\"{o}er made explicit correspondences between the two multiplication formulas as follows. They defined a cluster tilting object $W_{\mathbf{i}}$ in $\mathcal{C}_w$ and set $A:=\underline{\mathrm{End}}_{\mathcal{C}_w}(W_{\mathbf{i}})^{\rm op}$. For $X\in \mathcal{C}_w$, $F(X)=\mathrm{Ext}^1_{\Lambda}(W_{\mathbf{i}}, X)$ is an $A$-module. There is a bijiection $d_{\mathbf{i}, M}$ between $\{\mathbf{a}\in \mathbb{N}^r\mid \mathcal{F}_{\mathbf{i}, \mathbf{a}}(M)\neq \varnothing\}$ and $\{\underline{d}\in \mathbb{N}^n\mid \operatorname{Gr}_{\underline{d}}(F(M))\neq \varnothing\}$. Furthermore, there is an isomorphism of algebraic varieties between $\mathcal{F}_{\mathbf{i}, \mathbf{a}}(M)$ and $\operatorname{Gr}_{d_{\mathbf{i}, M}(\mathbf{a})}(F(M))$. These lead to a correspondence between $\Delta_{M}$ and $X^{T}_{\underline{M}}$ for $M\in \mathcal{C}_w$ where $T=\underline{W_{\mathbf{i}}}$ and $\underline{M}$ are respectively the image of $W_{\mathbf{i}}$ and $M$ under the natural functor from $\mathcal{C}_w$ to $\underline{\mathcal{C}_w}.$

Quantum cluster algebras were introduced by Berenstein and Zelevinsky \cite{BZ05} as the quantization of cluster algebras. They are defined to be certain noncommutative algebras over $\mathbb{Z}[q^{\pm}]$ generated by quantum cluster variables for a compatible pair $(B, \Lambda)$ of two matrices $B$ and $\Lambda$. Rupel \cite{Rup2011} defined a quantum cluster character as a quantum analogue of the Caldero--Chapoton map. Let $Q$ be an acyclic valued quiver over a finite field $\mathbb{F}=\mathbb{F}_q$. For $M\in \mathrm{Rep}_{\mathbb{F}}Q$,
\[X^{\mathbb{F}}_M=\sum_{\underline{e}}q^{-\frac{1}{2}\langle \underline{e}, \underline{m}-\underline{e}\rangle}|\mathrm{Gr}_{\underline{e}}(M)|x^{-B\underline{e}-(I-R)\underline{m}}\]
 (see Section \ref{sec5} for the explicit definition). Qin \cite{Qin2} gave an alternative definition for the quantum cluster characters of rigid objects in the cluster category of an acyclic equally valued quiver. The definition does not involve the choice of the ground finite field. Then he proved that the definition is consistent with the mutation rules between quantum cluster variables, i.e., for indecomposable rigid objects $M$, $N$ in $\mathcal{C}_Q$ with $\dim _{\mathbb{C}}\operatorname{Hom}_{\mathcal{C}_Q}(M, N[1])=1$, we have
 \[X^{\mathbb{F}}_M\cdot X^{\mathbb{F}}_N=q^{\frac{1}{2}\Lambda(\operatorname{ind}M, \operatorname{ind}L)}X^{\mathbb{F}}_{L}+q^{\frac{1}{2}\Lambda(\operatorname{ind}M, \operatorname{ind}L')}X^{\mathbb{F}}_{L'}.\]
 The multiplication formula of quantum cluster characters was generalized in \cite{DX}. More generally, we give a quantum analogue of Caldero--Keller's multiplication formula in Section \ref{sec4}. Recently, another quantum analogue was proved in \cite{CDZ}. The exact relation between these two formulas will be discovered in the near future.

In contrast to the case of the cluster character, the general definition of the quantum cluster character is still unknown. The aim of this paper is to define the weighted quantum cluster functions for both an abelian category $\mathcal{A}$ with the Ext-symmetry and a 2-Calabi--Yau triangulated category $\mathcal{C}$ with a cluster tilting object. A weighted quantum cluster function is of the form~$f*_{[\epsilon]}\delta_L$ or~$f*_{\epsilon}X_L$ with a weight function $f$.
So far we do not know whether or not weighted quantum cluster functions can give quantum cluster characters which compute the quantum cluster variables
in the sense of Berenstein--Zelevinsky by taking the proper weights.

In this paper, we use the following notation: given a finite set $S$, and a function $g$ on $S$, we define
\[\int_{x\in S} g(x):=\displaystyle\sum_{x\in S} g(x).\]

The first main result in this paper is the following theorem (cf.\ Theorem~\ref{maintheorem1} for details).
\begin{Theorem}\label{th1}
Let $\mathcal{A}$ be a ${\rm Hom}$-finite, $\Ext$-finite abelian category with $\EExt$-symmetry over a~finite field $k=\mathbb{F}_q$ such that the iso-classes of objects form a set. For any weighted quantum cluster functions
$f*_{[\epsilon']}\delta_M$ and $g*_{[\epsilon'']}\delta_N$ such that $\operatorname{Ext}_{\mathcal{A}}^1(M, N)\neq 0$, we have
\begin{align*}
&\big|\mathbb{P}\EExt_{\mathcal{A}}^1(M,N)\big|(f*_{[\epsilon']}\delta_M)*(g*_{[\epsilon'']}\delta_N) \\
&\quad=\int_{\mathbb{P}\epsilon\in\mathbb{P}\EExt_{\mathcal{A}}^1(M,N)}\bigl(f^+_{\rm ext}*_{[\epsilon]} \mathbb{S}_{fg}\bigr) *_{[\epsilon]} \delta_{{\rm mt}\epsilon}+\int_{\mathbb{P}\eta\in\mathbb{P}\EExt_{\mathcal{A}}^1(N,M)}(f_{\rm hom}*_{[\eta]}
\mathbb{S}_{gf})*_{[\eta]} \delta_{{\rm mt}\eta}\\
&\quad=\int_{\mathbb{P}\epsilon\in\mathbb{P}\EExt_{\mathcal{A}}^1(M,N)} \mathbb{S}_{fg}*_{[\epsilon]}\delta_{{\rm mt}\epsilon}+\int_{\mathbb{P}\eta\in\mathbb{P}\EExt_{\mathcal{A}}^1(N,M)}(f^-_{\rm ext}*_{[\eta]} f_{\rm hom}*_{[\eta]} \mathbb{S}_{gf})*_{[\eta]}\delta_{{\rm mt}\eta}.
\end{align*}
\end{Theorem}
This theorem gives the explicit multiplication between two weighted quantum cluster functions over an abelian category $\mathcal{A}$ with the Ext-symmetry.

The second main result is as follows (cf.\ Theorem~\ref{maintheorem2} for details).
\begin{Theorem}\label{th2}
Let $\mathcal{C}$ be a ${\rm Hom}$-finite, $2$-Calabi--Yau, Krull--Schmidt triangulated category over a finite field $k=\mathbb{F}_q$ with a cluster tilting object $T$. For any weighted quantum cluster functions
 $f*_{\epsilon'}X_M$ and $g*_{\epsilon''}X_N$ such that ${\rm Hom}_{\mathcal{C}}(M,\Sigma N)\neq0$, we have
\begin{gather*}
|\mathbb{P}\operatorname{Hom}_{\mathcal{C}}(M,\Sigma N)|(f*_{\epsilon'}X_M)\cdot(g*_{\epsilon''}X_N)\\
\quad=\int_{\mathbb{P}\epsilon\in\mathbb{P}\operatorname{Hom}_{\mathcal{C}}(M,\Sigma N)}\bigl(g^+_{\rm ext}*_\epsilon {f}_{\rm spec}*_\epsilon \mathbb{T}_{fg}\bigr)*_\epsilon X_{{\rm mt}\epsilon}\\
\phantom{\quad=}{}+\int_{\mathbb{P}\eta\in\mathbb{P}\operatorname{Hom}_{\mathcal{C}}(N,\Sigma M)}(g_{\rm skew}*_\eta {f}_{\rm spec}*_\eta \mathbb{T}_{gf})*_\eta X_{{\rm mt}\eta}\\
\quad=\int_{\mathbb{P}\epsilon\in\mathbb{P}\operatorname{Hom}_{\mathcal{C}}(M,\Sigma N)}({f}_{\rm spec}*_\epsilon \mathbb{T}_{fg})*_\epsilon X_{{\rm mt}\epsilon}\\
\phantom{\quad=}{}+\int_{\mathbb{P}\eta\in\mathbb{P}\operatorname{Hom}_{\mathcal{C}}(N,\Sigma M)}(g^-_{\rm ext}*_\eta g_{\rm skew}*_\eta {f}_{\rm spec}*_\eta \mathbb{T}_{gf})*_\eta X_{{\rm mt}\eta}.
\end{gather*}
\end{Theorem}
This theorem gives the explicit multiplication between two weighted quantum cluster functions over a 2-Calabi--Yau triangulated category $\mathcal{C}$ with a cluster tilting object $T$. It can be viewed as the quantum analogue of Caldero--Keller's multiplication formula.

The third main result is the quantum analogue of the above Geiss--Leclerc--Schr\"{o}er's multiplication formula (cf.\ Theorem~\ref{maintheorem3} for details).
\begin{Theorem}\label{th3}
For any weighted quantum cluster functions $f*_{\epsilon'}\Delta_{\mathbf{i},M}$ and $g*_{\epsilon''}\Delta_{\mathbf{i},N}$ such that $\EExt^1_{\mathcal{C}_{\omega}}(M,N)\neq 0$, in $A_{m,\lambda}$ we have
\begin{gather*}
\big|\mathbb{P}\EExt_{\mathcal{C}_{\omega}}^1(M,N)\big| (f*_{\epsilon'}\Delta_{\mathbf{i},M} )\cdot (g*_{\epsilon''}\Delta_{\mathbf{i},N}) \\
\quad=\int_{\mathbb{P}\epsilon\in\mathbb{P}\EExt_{\mathcal{C}_{\omega}}^1(M,N)}\bigl(f^+_{\rm ext}*_\epsilon {f}_{\rm spec}*_\epsilon \mathbb{S}_{fg}\bigr)*_\epsilon \Delta_{\mathbf{i},{\rm mt}\epsilon}\\
\phantom{\quad=}{}+\int_{\mathbb{P}\eta\in\mathbb{P}\EExt_{\mathcal{C}_{\omega}}^1(N,M)}(f_{\rm skew}*_\eta {f}_{\rm spec}*_\eta \mathbb{S}_{gf})*_\eta \Delta_{\mathbf{i},{\rm mt}\eta}\\
\quad=\int_{\mathbb{P}\epsilon\in\mathbb{P}\EExt_{\mathcal{C}_{\omega}}^1(M,N)}({f}_{\rm spec}*_\epsilon \mathbb{S}_{fg})*_\epsilon \Delta_{\mathbf{i},{\rm mt}\epsilon}\\
\phantom{\quad=}{}+\int_{\mathbb{P}\eta\in\mathbb{P}\EExt_{\mathcal{C}_{\omega}}^1(N,M)}(f^-_{\rm ext}*_\eta f_{\rm skew}*_\eta {f}_{\rm spec}*_\eta \mathbb{S}_{gf})*_\eta \Delta_{\mathbf{i},{\rm mt}\eta}.
\end{gather*}
\end{Theorem}

The fourth main result focuses on weighted quantum cluster functions from hereditary algebras. By choosing a particular weight function, we recover Qin's result~\cite[Proposition 5.4.1]{fanqin} (cf.\ Theorem~\ref{maintheorem4} for details).
\begin{Theorem}\label{th4}
In the cluster category $\mathcal{C}=D^b\bigl(\widetilde{\mathcal{A}}\bigr)/{\tau^{-1}\Sigma}$ of a hereditary algebra $\widetilde{A}$, given two indecomposable coefficient-free rigid
objects $M$, \smash{$N \in \widetilde{\mathcal{A}}$} with
\[\dim_k\Hom_{\mathcal{C}}(M,\Sigma N)=1\]
 and two non-split triangles
\[
N\rightarrow L\rightarrow M\xrightarrow{\epsilon}\Sigma N \qquad \text{and}\qquad
M\rightarrow L'\rightarrow N\xrightarrow{\eta}\Sigma M,
\]
where $L'$ is located in the fundamental domain, then we have
\[
\tilde{X}_M\cdot\tilde{X}_N=q^{\frac{1}{2}\lambda(\operatorname{ind}M,\operatorname{ind}N)-\frac{1}{2}}\cdot\tilde{X}_L+q^{\frac{1}{2}\lambda(\operatorname{ind}M,\operatorname{ind}N)}\cdot\tilde{X}_{L'}.
\]
\end{Theorem}

Recently, Keller, Plamondon and Qin gave a refined multiplication formula for cluster characters over $2$-Calabi--Yau triangulated category with tilting objects \cite{KPQ}. Its quantum analogue has been obtained in \cite{XXY}
via a motivic version of weighted quantum cluster functions.

The paper is organized as follows. In Section~\ref{sec2}, we first define a chain of monomorphisms as a~generalization of flag and then define a set of chains of monomorphisms of certain type. In Section~\ref{sec3}, we introduce the quantum cluster function attached to an abelian category with the Ext-symmetry in a general context. Corollaries~\ref{vectorbundle1} and~\ref{vectorbundle2} play an important role in characterizing the quantum cluster function. Then we define the weight function and the weighted quantum cluster function. By choosing appropriate pair of weight functions, we prove the multiplication formula of weighted quantum cluster functions in Section \ref{sectionofmaintheorem1}. Section~\ref{sec4} is devoted to defining the weighted quantum cluster function for a 2-Calabi--Yau triangulated category with a cluster tilting object. Lemma \ref{vectorbundle3} helps to describe the structure of the quantum cluster function. We prove the multiplication formula of weighted quantum cluster functions in Section~\ref{sectionofmaintheorem2}. Section~\ref{sec5} establishes the explicit connection between the two multiplication formulas as stated in Theorem~\ref{th1} and Theorem~\ref{th2}. We make a correspondence between weighted quantum cluster functions in Section~\ref{sec3} and Section~\ref{sec4} by applying Geiss--Leclerc--Schr\"{o}er's correspondence as above. In Section~\ref{sec6}, we show that Theorem \ref{maintheorem2} recovers the multiplcation formula in \cite{fanqin} by assigning proper weight functions.

\section{Chains of monomorphisms}\label{sec2}

\subsection{Chains of morphisms}

Let $\mathcal{A}$ be a ${\operatorname{Hom}}$-finite, $\operatorname{Ext}$-finite abelian category with $\operatorname{Ext}$-symmetry over a finite field $k=\mathbb{F}_q$ such that the isoclasses of objects form a set.
We assume that $\mathcal{A}$ has finitely many simple objects up to isomorphism and fix a complete set of simple objects $\{S_1,\dots,S_n\}$ up to isomorphism in~$\mathcal{A}$.
For a fixed $m \in \mathbb{N}$, we denote by $\widetilde{\mathcal{F}}^{\rm mor}_m$ the set of all chains of morphisms
\[\mathbf{c}=\bigl(L_m\xrightarrow{\iota_{\mathbf{c}, m}}L_{m-1}\longrightarrow\dots\longrightarrow L_1\xrightarrow{\iota_{\mathbf{c}, 1}}L_0\bigr),\]
where $L_0, \dots, L_m\in\mathcal{A}$ and $\iota_{\mathbf{c},1}, \dots, \iota_{\mathbf{c}, m}\in\rm{Mor}\mathcal{A}$.

Using the isomorphism relations $L\cong L'$ in $\mathcal{A}$, we can induce the isomorphism relations in~$\widetilde{\mathcal{F}}^{\rm mor}_m$.
\begin{Definition}
Two chains of morphisms $\mathbf{c}=\bigl(L_m\xrightarrow{\iota_{\mathbf{c},m}}L_{m-1}\longrightarrow\dots\longrightarrow L_1\xrightarrow{\iota_{\mathbf{c},1}}L_0\bigr)$ and \smash{$\mathbf{c}'=\bigl(L'_m\xrightarrow{\iota_{\mathbf{c}',m}}L'_{m-1}\longrightarrow\dots\longrightarrow L'_1\xrightarrow{\iota_{\mathbf{c}',1}}L'_0\bigr)$} in \smash{$\widetilde{\mathcal{F}}^{\rm mor}_m$} are called isomorphic if there are isomorphisms \smash{$L_i\xrightarrow{f_i}L'_j$} in $\rm{Mor}\mathcal{A}$ such that $f_{i-1}\circ \iota_{\mathbf{c},i}=\iota_{\mathbf{c}',i}\circ f_i$, i.e., we have the following commutative diagram:
\[
\begin{tikzcd}
L_m \arrow[r,"\iota_{\mathbf{c},m}"]\arrow[d,"f_m"] & L_{m-1} \arrow[r,"\iota_{\mathbf{c},m-1}"]\arrow[d,"f_{m-1}"]&\cdots\arrow[r]&L_1\arrow[r,"\iota_{\mathbf{c},1}"]\arrow[d,"f_1"]&L_0\arrow[d,"f_0"] \\
L'_m \arrow[r,"\iota_{\mathbf{c}',m}"] & L'_{m-1} \arrow[r,"\iota_{\mathbf{c}',m-1}"]&\cdots\arrow[r]&L'_1\arrow[r,"\iota_{\mathbf{c}',1}"]&L'_0.
\end{tikzcd}
\]
\end{Definition}
We denote the set of isomorphism classes in $\widetilde{\mathcal{F}}^{\rm mor}_m$ by $\mathcal{F}^{\rm mor}_m$ and we still write $\mathbf{c}$ for its isomorphism class in $\mathcal{F}^{\rm mor}_m$.

\subsection{Exact structure}
Given a short exact sequence
\[
\xymatrix{
\epsilon\colon\ 0\ar[r]&N\ar[r]&L\ar[r]&M\ar[r]&0}
\]
in $\mathcal{A}$, write $[\epsilon]$ for its equivalence class in $\Ext_{\mathcal{A}}^1(M,N)$. Recall that two short exact sequences
\[
\xymatrix{
0\ar[r]&N\ar[r]^{f_1}&L\ar[r]^{g_1}&M\ar[r]&0} \qquad \text{and}\qquad \xymatrix{
0\ar[r]&N\ar[r]^{f_2}&L'\ar[r]^{g_2}&M\ar[r]&0}
\]
correspond to the same element in $\Ext_{\mathcal{A}}^1(M,N)$ if and only if there is an $\mathcal{A}$-isomorphism ${a \colon L \longrightarrow L'}$ with $f_2= a f_1$ and $g_2=g_1 a^{-1}$, i.e., we have the following commutative diagram:
\[
\xymatrix{
&0\ar[r]&N \ar[r]^{f_2}& L'\ar[r]^{g_2} &M \ar[r]&0\,\,\\
&0\ar[r]&N \ar[r]^{f_1}\ar@{=}[u]& L\ar[r]^{g_1}\ar[u]^{a}&M \ar[r]\ar@{=}[u]&0.}
\]
Then any morphism $\lambda\in {\rm Mor}(N,N')$ induces an equivalence class $\lambda\circ[\epsilon]$ in $\Ext_{\mathcal{A}}^1(M,N')$ which is called the pushout of $[\epsilon]$ along $\lambda$. Similarly, for any $\rho\in{\rm Mor}(M'',M)$, there is an equivalence class $[\epsilon]\circ\rho$ in $\Ext_{\mathcal{A}}^1(M'',N)$ which is called the pullback of $[\epsilon]$ along $\rho$. We have the following commutative diagram:
\[\xymatrix{
\lambda\circ[\epsilon]\colon \quad 0\ar[r]&N'\ar[r]& L'\ar[r]&M\ar[r]&0\,\,\\
[\epsilon]\colon\quad &0\ar[r] N\ar[u]^{\lambda}\ar[r] &L\ar[r]\ar[u]&M\ar[r]\ar@{=}[u]&0\,\, \\
[\epsilon]\circ\rho\colon\quad 0\ar[r]&N\ar[r]\ar@{=}[u]& L''\ar[r]\ar[u]&M''\ar[r]\ar[u]_{\rho}&0.}
\]
Conversely, given
\begin{gather*}
[\epsilon']\in\Ext_{\mathcal{A}}^1(M',N'),\qquad [\epsilon'']\in\Ext_{\mathcal{A}}^1(M'',N''),\qquad\lambda\in{\rm Mor}(N'',N'),\\ \rho\in{\rm Mor}(M'',M')
\end{gather*}
as in the diagram
\[
\xymatrix{
[\epsilon']\colon\quad 0\ar[r]&N'\ar[r]& L'\ar[r]&M'\ar[r]&0\,\, \\
[\epsilon'']\colon\quad 0\ar[r]&N''\ar[r]\ar[u]^\lambda& L''\ar[r]&M''\ar[r]\ar[u]_\rho&0,}
\]
we can complete the commutative diagram by adding an appropriate morphism in ${\rm Mor}(L'',L')$ if and only if the pushout and pullback coincide, i.e.,
\[
\lambda\circ[\epsilon'']=[\epsilon']\circ\rho.
\]
Moreover, in this case, the morphism in ${\rm Mor}(L'',L')$ is given by $g\circ f$ as shown in the following diagram:
\[
\xymatrix{
[\epsilon']\colon&0\ar[r]&N'\ar[r]& L'\ar[r]&M'\ar[r]&0\,\, \\
\lambda\circ[\epsilon'']=[\epsilon']\circ\rho\colon&0\ar[r]&N'\ar[r]\ar@{=}[u]&L\ar[r]\ar[u]_g&M''\ar[r]\ar[u]_\rho&0\,\, \\
[\epsilon'']\colon&0\ar[r]&N''\ar[r]\ar[u]^\lambda& L''\ar[r]\ar[u]^f&M''\ar[r]\ar@{=}[u]&0.}
\]

The equality $[\epsilon']\circ \rho= \lambda \circ [\epsilon']$ can be illustrated by the following commutative diagram:
\[
\xymatrix{
[\epsilon']\colon&0\ar[r]&N'\ar[r] & L'\ar[r] &M'\ar[r]&0 \,\, \\
[\epsilon']\circ\rho\colon&0\ar[r]&N'\ar[r]\ar@{=}[u]&\widetilde{L}'\ar[r]\ar[u] \ar[d]^{\cong} &M''\ar[r]\ar[u]_\rho&0\,\, \\
\lambda\circ[\epsilon']\ar@{=}[u]\colon&0\ar[r]&N''\ar[r] \ar@{=}[u] & \widetilde{L}''\ar[r]&M''\ar[r]\ar@{=}[u]&0\,\, \\
[\epsilon'']\colon&0\ar[r]&N''\ar[r]\ar[u]^{\lambda} & L''\ar[r]\ar[u] &M''\ar[r]\ar@{=}[u]&0.}
\]

Since every $[\epsilon]\in\Ext_{\mathcal{A}}^1(M,N)$ can be represented by
\[
\xymatrix{
\epsilon\colon\ 0\ar[r]&N\ar[r]&L\ar[r]&M\ar[r]&0}
\]
with $M$, $N$, $L$ being unique up to isomorphism, we denote ${\rm qt}\epsilon:=M$, ${\rm st}\epsilon:=N$ and ${\rm mt}\epsilon:=L$ in the following sections.

Now, more generally, let $[\epsilon_i]\in\Ext_{\mathcal{A}}^1(M_i,N_i)$ for $0 \leq i \leq m$, and take two isomorphism classes of chains of morphisms
\begin{gather*}
\mathbf{c}'=\bigl(M_m\xrightarrow{\iota_{\mathbf{c}',m}}M_{m-1}\longrightarrow\dots\longrightarrow M_1\xrightarrow{\iota_{\mathbf{c}',1}}M_0=M\bigr),\\
\mathbf{c}''=\bigl(N_m\xrightarrow{\iota_{\mathbf{c}'',m}}N_{m-1}\longrightarrow\dots\longrightarrow N_1\xrightarrow{\iota_{\mathbf{c}'',1}}N_0=N\bigr)
\end{gather*}
in $\mathcal{F}^{\rm mor}_m$ as in the diagram
\begin{equation}\label{complete-diagram}
\begin{split}&
\xymatrix{
[\epsilon_0]\colon&0\ar[r]&N_0\ar[r]& L_0\ar[r]&M_0\ar[r]&0\,\, \\
[\epsilon_1]\colon&0\ar[r]&N_1\ar[r]\ar[u]^{\iota_{\mathbf{c}'',1}}& L_1\ar[r]&M_1\ar[r]\ar[u]_{\iota_{\mathbf{c}',1}}&0\, \, \\
&&\vdots\ar[u]&\vdots&\vdots\ar[u]&\\
[\epsilon_{m-1}]\colon&0\ar[r]&N_{m-1}\ar[r]\ar[u]& L_{m-1}\ar[r]&M_{m-1}\ar[r]\ar[u]&0\,\, \\
[\epsilon_m]\colon&0\ar[r]&N_m\ar[r]\ar[u]^{\iota_{\mathbf{c}'',m}}& L_m\ar[r]&M_m\ar[r]\ar[u]_{\iota_{\mathbf{c}',m}}&0.}
\end{split}
\end{equation}
In order to decide whether there exists a chain of morphisms
\[
\mathbf{c}=\bigl(L_m\xrightarrow{\iota_{\mathbf{c},m}}L_{m-1}\longrightarrow\dots\longrightarrow L_1\xrightarrow{\iota_{\mathbf{c},1}}L_0=L\bigr)
\]
to complete the above commutative diagram, we introduce a linear map 
\begin{align*}
\beta_{\mathbf{c}',\mathbf{c}''}\colon\ \bigoplus_{j=0}^m\Ext^1(M_j,N_j)&\longrightarrow\bigoplus_{j=0}^{m-1}\EExt^1(M_{j+1},N_j),\\
([\epsilon_0],\dots,[\epsilon_m])&\longmapsto([\epsilon_{j-1}]\circ\iota_{\mathbf{c}',j}-\iota_{\mathbf{c}'',j}\circ[\epsilon_j], 1\leq j \leq m).
\end{align*}
The definition and some properties of maps $\beta_{\mathbf{c}',\mathbf{c}''}$ and $\beta'_{\mathbf{c}'',\mathbf{c}'}$, which will be introduced later as the
linear dual of $\beta_{\mathbf{c}',\mathbf{c}''}$, were firstly given in \cite{CK2005} for module categories.
We generalize these two maps to abelian categories and prove some necessary properties here.

\begin{Lemma}\label{chain-mono}
There exists a chain of morphisms
\[
\mathbf{c}=\bigl(L_m\xrightarrow{\iota_{\mathbf{c},m}}L_{m-1}\longrightarrow\dots\longrightarrow L_1\xrightarrow{\iota_{\mathbf{c},1}}L_0=L\bigr),
\]
which can complete the commutative diagram \eqref{complete-diagram} if and only if
\[
([\epsilon_0],\dots,[\epsilon_m])\in\operatorname{Ker}\beta_{\mathbf{c}',\mathbf{c}''}.
\]
\end{Lemma}
\begin{proof}
A family of elements $([\epsilon_0],\dots,[\epsilon_m])\in\operatorname{Ker}\beta_{\mathbf{c}',\mathbf{c}''}$ is exactly one which makes pushout and pullback coincide at every level.
\end{proof}

In this case, we get the morphism $\iota_{L,i}$ by composing two vertical morphisms in the following diagram:
\[
\xymatrix{
[\epsilon_{i-1}]\colon&0\ar[r]&N_{i-1}\ar[r]& L_{i-1}\ar[r]&M_{i-1}\ar[r]&0 \,\, \\
\iota_{\mathbf{c}'',i}\circ[\epsilon_i]=[\epsilon_{i-1}]\circ\iota_{\mathbf{c}',i}\colon&0\ar[r]&N_{i-1}\ar[r]\ar@{=}[u]&L'_i\ar[r]\ar[u]_g&M_i\ar[r]\ar[u]_{\iota_{\mathbf{c}',i}}&0\,\, \\
[\epsilon_i]\colon&0\ar[r]&N_i\ar[r]\ar[u]^{\iota_{\mathbf{c}'',i}}& L_i\ar[r]\ar[u]^f&M_i\ar[r]\ar@{=}[u]&0.}
\]
We denote this assignment by
\begin{align*}
B_{\mathbf{c}',\mathbf{c}''}\colon\ \operatorname{Ker}\beta_{\mathbf{c}',\mathbf{c}''}& \longrightarrow\mathcal{F}^{\rm mor}_m,\\
([\epsilon_0],\dots,[\epsilon_m])& \longmapsto(L_m\xrightarrow{\iota_{\mathbf{c},m}}L_{m-1}\longrightarrow\dots\longrightarrow L_1\xrightarrow{\iota_{\mathbf{c},1}}L_0=L).
\end{align*}

\begin{Lemma}\label{chain-mono-inj}
For any chains of morphisms $\mathbf{c}'$, $\mathbf{c}''$, the map $B_{\mathbf{c}',\mathbf{c}''}$ is well defined and injective.
\end{Lemma}
\begin{proof}
The proof is similar to the discussion in the case when $m=1$. The choice of short exact sequences representing $([\epsilon_0],\dots,[\epsilon_m])$ is not unique, but is unique up to the following commutative diagram:
\[
\xymatrix{
\epsilon_i\colon\quad 0\ar[r]&N_i\ar[r]& L_i\ar[r]&M_i\ar[r]&0 \, \, \\
\epsilon'_i\colon\quad 0\ar[r]&N_i\ar[r]\ar@{=}[u]& L'_i\ar[r]\ar[u]^{a_i}&M_i\ar[r]\ar@{=}[u]&0,}
\]
where $a_i$ is an isomorphism. Then we get two chains of morphisms $B_{\mathbf{c}',\mathbf{c}''}([\epsilon_0],\dots,[\epsilon_m])$ and $B_{\mathbf{c}',\mathbf{c}''}([\epsilon'_0],\dots,[\epsilon'_m])$. One can check they are isomorphic in \smash{$\widetilde{\mathcal{F}}^{\rm mor}_m$} through the family of isomorphisms $a_i$.

Conversely, given two choices of exact sequences which induce the same chain of morphisms in~$\mathcal{F}^{\rm mor}_m$, then the isomorphism $L_i \rightarrow L'_i$ also give the equivalence between the two exact sequences:
\[
\begin{tikzcd}
 & N_{i-1} \arrow[rr] & & L'_{i-1} \arrow[rr] & & M_{i-1} \,\, \\
N_i \arrow[ru] \arrow[rr] & & L'_i \arrow[ru] \arrow[rr, crossing over] & & M_i \arrow[ru] & \\
  & N_{i-1} \arrow[rr] \arrow[uu,equal] & & L_{i-1} \arrow[rr] \arrow[uu, "\simeq"] & & M_{i-1} \arrow[uu,equal] . \qedhere \\
N_i \arrow[ru] \arrow[rr] \arrow[uu,equal] & & L_i \arrow[ru] \arrow[rr] \arrow[uu, "\simeq", crossing over] & & M_i \arrow[ru] \arrow[uu,equal] &
\end{tikzcd}\tag*{\qed}
\]\renewcommand{\qed}{}
\end{proof}

Since all extension groups considered here are finite dimensional over a finite field, we have
\begin{Corollary}
Given chains of morphisms $\mathbf{c}'$, $\mathbf{c}''$, the map $B_{\mathbf{c}',\mathbf{c}''}\colon\Ker \beta_{\mathbf{c}',\mathbf{c}''}\longrightarrow\operatorname{Im}B_{\mathbf{c}',\mathbf{c}''}$ is a~bijection and $|\operatorname{Im}B_{\mathbf{c}',\mathbf{c}''}|=|\KKer \beta_{\mathbf{c}',\mathbf{c}''}|$.
\end{Corollary}

\subsection{Chains of monomorphisms}
\begin{Definition}
Given a chain of morphisms
\[
\mathbf{c}=\bigl(L_m\xrightarrow{\iota_{\mathbf{c},m}}L_{m-1}\longrightarrow\dots\longrightarrow L_1\xrightarrow{\iota_{\mathbf{c},1}}L_0\bigr)
\]
in $\widetilde{\mathcal{F}}^{\rm mor}_m$, it is called a chain of monomorphisms if all $\iota_{\mathbf{c},i}$ are monomorphisms and $L_m=0$.
\end{Definition}

One can easily check this definition is independent of choice of chain in an isomorphism class and we denote by $\mathcal{F}^{\rm mono}_m$ the subset of $\mathcal{F}^{\rm mor}_m$ consisting of all isomorphism classes of chains of monomorphisms. We also denote by $\mathcal{F}^{\rm mono}_{m, L}$ the subset of $\mathcal{F}^{\rm mono}_m$ consisting of all isomorphism classes of chains of monomorphisms with $L_0\cong L$.
\begin{Lemma}
If $\mathbf{c}',\mathbf{c}''\in\mathcal{F}^{\rm mono}_m$, then $\operatorname{Im}B_{\mathbf{c}',\mathbf{c}''}\subseteq\mathcal{F}^{\rm mono}_m$.
\end{Lemma}
\begin{proof}
By the snake lemma, all rows in the following commutative diagram are exact:
\[
\begin{tikzcd}
0 \arrow[r] & N_{j-1} \arrow[r] & L_{j-1} \arrow[r] & M_{j-1} \arrow[r] & 0\,\, \\
0 \arrow[r] & N_j \arrow[r] \arrow[u] & L_j \arrow[r] \arrow[u] & M_j \arrow[u] \arrow[r] & 0. \\
 & {\operatorname{Ker}\iota_{\mathbf{c}'',j}} \arrow[u] \arrow[r] & {\operatorname{Ker}\iota_{\mathbf{c},j}} \arrow[u] \arrow[r] & {\operatorname{Ker}\iota_{\mathbf{c}',j}} \arrow[u] &
\end{tikzcd}
\]
Since $\mathbf{c}'$, $\mathbf{c}''$ are chains of monomorphisms, $\operatorname{Ker}\iota_{\mathbf{c}'',j}=0$ and $\operatorname{Ker}\iota_{\mathbf{c}',j}=0$. So is $\operatorname{Ker}\iota_{\mathbf{c},j}$.
\end{proof}

\subsection{Type}\label{type}

Recall that $\mathcal{A}$ admits a complete set simple objects $\{S_1,\dots,S_n\}$ up to isomorphism, so we can consider composition factors of an object in $\mathcal{A}$ with finite length.
\begin{Definition}
Given a chain of monomorphisms
\[
\mathbf{c}=\bigl(L_m\xrightarrow{\iota_{\mathbf{c},m}}L_{m-1}\longrightarrow\dots\longrightarrow L_1\xrightarrow{\iota_{\mathbf{c},1}}L_0\bigr)
\]
and two sequences $\mathbf{i}=(i_1,\dots,i_m)$, $\mathbf{a}=(a_1,\dots,a_m)$, where $i_j\in\{1,\dots,n\}$, $a_j\in\{0,1\}$, $\mathbf{c}$ is called of type $(\mathbf{i},\mathbf{a})$ if $\operatorname{Coker}\iota_{\mathbf{c},j} \cong S_{i_j}$ when $a_j=1$ and $0$ otherwise for $1\leq j \leq m$.
\end{Definition}

We denote the set of all chains of monomorphisms of type $(\mathbf{i},\mathbf{a})$ by $\mathcal{F}^{\rm mono}_{\mathbf{i},\mathbf{a}}$ and the set of all chains of monomorphisms of type $(\mathbf{i},\mathbf{a})$ with $L_0\cong L$ by $\mathcal{F}^{\rm mono}_{\mathbf{i},\mathbf{a},L}$.
\begin{Lemma}
If $(\mathbf{c}',\mathbf{c}'')\in\mathcal{F}^{\rm mono}_{\mathbf{i},\mathbf{a}',M}\times\mathcal{F}^{\rm mono}_{\mathbf{i},\mathbf{a}'',N}$ and $a'_j+a''_j \leq 1$ for $1 \leq j \leq m$, then \[\operatorname{Im}B_{\mathbf{c}',\mathbf{c}''}\subseteq\mathcal{F}^{\rm mono}_{\mathbf{i},\mathbf{a}'+\mathbf{a}''}.\]
\end{Lemma}
\begin{proof}
Since $a'_j+a''_j \leq 1$ for $1 \leq j \leq m$, $\mathcal{F}^{\rm mono}_{\mathbf{i},\mathbf{a}'+\mathbf{a}''}$ is well defined.
By the snake lemma, we have an exact sequence
\[
\operatorname{Ker}\iota_{\mathbf{c}',j}\rightarrow\operatorname{Coker}\iota_{\mathbf{c}'',j}\rightarrow\operatorname{Coker}\iota_{\mathbf{c},j}\rightarrow\operatorname{Coker}\iota_{\mathbf{c}',j}\rightarrow 0.
\]
Since $\mathbf{c}'$ is a chain of monomorphisms, $\operatorname{Ker}\iota_{\mathbf{c}',j}=0$. So the composition factors of $\operatorname{Coker}\iota_{\mathbf{c},j}$ are the union of composition factors of $\operatorname{Coker}\iota_{\mathbf{c}',j}$ and $\operatorname{Coker}\iota_{\mathbf{c}'',j}$, which are one copy of $S_{i_j}$.
\end{proof}

If $\mathbf{a}=(1, \dots, 1)$, we simply write $\mathcal{F}^{\rm mono}_{\mathbf{i}}$ instead of $\mathcal{F}^{\rm mono}_{\mathbf{i},\mathbf{a}}$.

\section[Abelian categories with Ext-symmetry and the multiplication formula]{Abelian categories with $\boldsymbol{\operatorname{Ext}}$-symmetry \\
and the multiplication formula}\label{sec3}

\subsection{Quantum cluster function}

Now we introduce the concept of the quantum cluster functions over $\mathcal{A}$.

Let $d_{\mathbf{i},\mathbf{a}}$ be a formal notation representing a type of chains of monomorphisms and
\[
\mathcal{M}_q:=\bigoplus_{(\mathbf{i},\mathbf{a})}\mathbb{Q}d_{\mathbf{i},\mathbf{a}}
\]
be the $\mathbb{Q}$-space spanned by all $d_{\mathbf{i},\mathbf{a}}$.

For each object $L$ in $\mathcal{A}$, we define a $\mathbb{Q}$-valued linear function $\delta_L$ on $\mathcal{M}_q$ by $\delta_L(d_{\mathbf{i},\mathbf{a}})=|\mathcal{F}^{\rm mono}_{\mathbf{i},\mathbf{a},L}|$ and call such a function the quantum cluster function of $L$.

The core purpose of this section is to study the relationship between quantum cluster functions of objects in $\mathcal{A}$ related by a short exact sequence.

\subsection{Mappings with affine fibers}

We denote by
\[
\EF_{\mathbf{i},\mathbf{a}}(M,N):=\big\{([\epsilon],\mathbf{c})\mid [\epsilon]\in\operatorname{Ext}_{\mathcal{A}}^1(M,N),\, \mathbf{c}\in\mathcal{F}^{\rm mono}_{\mathbf{i},\mathbf{a},{\rm mt}\epsilon}\big\}
\]
the set of all pairs of extensions of $M$, $N$ and chain of monomorphisms of the middle term of type $(\mathbf{i},\mathbf{a})$.

From $([\epsilon],\mathbf{c})\in \EF_{\mathbf{i},\mathbf{a}}(M,N)$, we can induce two chains of monomorphisms ending with $M$, $N$ respectively. More precisely, we have
\[
\begin{tikzcd}
0 \arrow[r] & N \arrow[r, "i"] & L \arrow[r, "p"] & M \arrow[r] & 0 \, \, \\
0 \arrow[r] & N_1 \arrow[u, "{\iota_{\mathbf{c}'',1}}", hook] \arrow[r, "i_1"] & L_1 \arrow[u, "{\iota_{\mathbf{c},1}}"] \arrow[r, "p_1"] & M_1 \arrow[u, "{\iota_{\mathbf{c}',1}}", hook] \arrow[r] & 0 \, \, \\
 & \vdots \arrow[u, hook] & \vdots \arrow[u] & \vdots \arrow[u, hook] & \\
0 \arrow[r] & N_{m-1} \arrow[r, "i_{m-1}"] \arrow[u, hook] & L_{m-1} \arrow[r, "p_{m-1}"] \arrow[u] & M_{m-1} \arrow[r] \arrow[u, hook] & 0 \,\,\\
0 \arrow[r] & N_m \arrow[r, "i_m"] \arrow[u, "{\iota_{\mathbf{c}'',m}}", hook] & L_m \arrow[r, "p_m"] \arrow[u, "{\iota_{\mathbf{c},m}}"] & M_m \arrow[r] \arrow[u, "{\iota_{\mathbf{c}',m}}", hook] & 0,
\end{tikzcd}
\]
where
\[
N_j=i^{-1}(\operatorname{Im}\iota_{\mathbf{c},1}\circ\dots\circ\iota_{\mathbf{c},j}), \qquad M_j=p(\operatorname{Im}\iota_{\mathbf{c},1}\circ\dots\circ\iota_{\mathbf{c},j}),\]
$\iota_{\mathbf{c}'', j}$, $\iota_{\mathbf{c}',j}$ are natural embeddings and $i_j$, $p_j$ are naturally induced by $i$ and $p$ respectively for~${1\leq j \leq m}$.

We write this assignment as
\[
\phi_{MN}\colon\ \EF_{\mathbf{i},\mathbf{a}}(M,N) \rightarrow\coprod_{\mathbf{a}'+\mathbf{a}''=\mathbf{a}}\mathcal{F}^{\rm mono}_{\mathbf{i},\mathbf{a}',M}\times\mathcal{F}^{\rm mono}_{\mathbf{i},\mathbf{a}'',N}.
\]
The coproduct runs over all pairs $(\mathbf{a}',\mathbf{a}'')$ with $\mathbf{a}'+\mathbf{a}''=\mathbf{a}$ since the composition factor at each level is fixed by $\mathbf{a}$.

Given two chains of monomorphisms $\mathbf{c}'$, $\mathbf{c}''$ ending with $M$ and $N$ respectively, we are interested in the fiber $\phi^{-1}_{MN}(\mathbf{c}',\mathbf{c}'')$.

\begin{Lemma}\label{fiber-related}
Consider chains of monomorphisms $\mathbf{c}'$, $\mathbf{c}''$ ending with $M$, $N$ respectively.
There exists a bijection between $\operatorname{Ker}\beta_{\mathbf{c}',\mathbf{c}''}$ and $\phi^{-1}_{MN}(\mathbf{c}',\mathbf{c}'')$, given by
\[([\epsilon_0],\dots,[\epsilon_m]) \longmapsto ([\epsilon_0], B_{\mathbf{c}',\mathbf{c}''}([\epsilon_0],\dots,[\epsilon_m])).\]
\end{Lemma}
\begin{proof}
Consider chains of monomorphisms $\mathbf{c}'$, $\mathbf{c}''$ ending with $M$, $N$ respectively and $[\epsilon_0]\in\EExt_{\mathcal{A}}^1(M,N)$ with ${\rm mt}\epsilon_0=L$ as in the
diagram
\[
\begin{tikzcd}
\epsilon_0 \colon \quad 0 \arrow[r] & N \arrow[r] & L \arrow[r] & M \arrow[r] & 0. \\
 & & N_1 \arrow[u, "{\iota_{\mathbf{c}'',1}}", hook] & & M_1 \arrow[u, "{\iota_{\mathbf{c}',1}}", hook] & \\
 & & \vdots \arrow[u, hook] & & \vdots \arrow[u, hook] & \\
 & & N_{m-1} \arrow[u, hook] & & M_{m-1} \arrow[u, hook] & \\
 & & N_m \arrow[u, "{\iota_{\mathbf{c}'',m}}", hook] & & M_m \arrow[u, "{\iota_{\mathbf{c}',m}}", hook] &
\end{tikzcd}
\]
For any $\mathbf{c}\in\mathcal{F}^{\rm mono}_{m,L}$, the following two statements are equivalent:
\begin{enumerate}\itemsep=0pt
\item[(1)] there are extensions $([\epsilon_1],\dots,[\epsilon_m])$ satisfying $([\epsilon_0],[\epsilon_1], \dots,[\epsilon_m])\in\KKer\beta_{\mathbf{c}',\mathbf{c}''}$ such that $B_{\mathbf{c}',\mathbf{c}''}([\epsilon_0],[\epsilon_1],\dots,[\epsilon_m])=\mathbf{c}$. By definition of $B_{\mathbf{c}',\mathbf{c}''}$, this means that there is a commutative diagram
\[
\begin{tikzcd}
\epsilon_0\colon & 0 \arrow[r] & N \arrow[r, "i"] & L \arrow[r, "p"] & M \arrow[r] & 0 \,\, \\
\epsilon_1\colon & 0 \arrow[r] & N_1 \arrow[u, hook] \arrow[r, "i_1"] & L_1 \arrow[u, "{\iota_{\mathbf{c},1}}", dashed, hook] \arrow[r, "p_1"] & M_1 \arrow[u, hook] \arrow[r] & 0 \,\, \\
 & & \vdots \arrow[u, hook] & \vdots \arrow[u, dashed, hook] & \vdots \arrow[u, hook] & \\
\epsilon_{m-1}\colon & 0 \arrow[r] & N_{m-1} \arrow[u, hook] \arrow[r, "i_{m-1}"] & L_{m-1} \arrow[u, dashed, hook] \arrow[r, "p_{m-1}"] & M_{m-1} \arrow[u, hook] \arrow[r] & 0 \,\, \\
\epsilon_m\colon & 0 \arrow[r] & N_m \arrow[r, "i_m"] \arrow[u, hook] & L_m \arrow[r, "p_m"] \arrow[u, "{\iota_{\mathbf{c},m}}", dashed, hook] & M_m \arrow[r] \arrow[u, hook] & 0,
\end{tikzcd}
\]
where $\iota_{\mathbf{c},i}$ are given by composing the middle vertical morphisms in diagrams of the form
\[
\xymatrix{
0\ar[r]&N_{i-1}\ar[r]& L_{i-1}\ar[r]&M_{i-1}\ar[r]&0 \,\, \\
0\ar[r]&N_{i-1}\ar[r]\ar@{=}[u]&L'_i\ar[r]\ar[u]_g&M_i\ar[r]\ar[u]_{\iota_{\mathbf{c}',i}}&0\,\, \\
0\ar[r]&N_i\ar[r]\ar[u]^{\iota_{\mathbf{c}'',i}}& L_i\ar[r]\ar[u]^f&M_i\ar[r]\ar@{=}[u]&0;}
\]
\item[(2)] $\phi_{MN}([\epsilon_0],\mathbf{c})=(\mathbf{c}',\mathbf{c}'')$. By definition of $\phi_{MN}$, this means that there is a commutative diagram
\[
\begin{tikzcd}
\epsilon_0\colon & 0 \arrow[r] & N \arrow[r, "i"] & L \arrow[r, "p"] & M \arrow[r] & 0,\\
 & & N_1 \arrow[u, hook] \arrow[r, "i_1", dashed] & L_1 \arrow[u, "{\iota_{\mathbf{c},1}}", hook] \arrow[r, "p_1", dashed] & M_1 \arrow[u, hook] & \\
 & & \vdots \arrow[u, hook] & \vdots \arrow[u, hook] & \vdots \arrow[u, hook] & \\
 & & N_{m-1} \arrow[u, hook] \arrow[r, "i_{m-1}", dashed] & L_{m-1} \arrow[u, hook] \arrow[r, "p_{m-1}", dashed] & M_{m-1} \arrow[u, hook] & \\
 & & N_m \arrow[r, "i_m", dashed] \arrow[u, hook] & L_m \arrow[r, "p_m", dashed] \arrow[u, "{\iota_{\mathbf{c},m}}", hook] & M_m \arrow[u, hook] &
\end{tikzcd}
\]
where $N_j=i_j^{-1}(L_j)$, $M_j=p_j(L_j)$, and $i_j$ and $p_j$ are restrictions of $i$ and $p$.
\end{enumerate}

Assume (1) is true. Since all rows in $(1)$ are short exact sequences, we have $N_j=i_j^{-1}(L_j)$ and $M_j=p_j(L_j)$, i.e., $(2)$ holds.

Conversely, given $(2)$, by definitions of $N_j$ and $M_j$, we know that $N_j \stackrel{i_j}{\longrightarrow} L_j \stackrel{p_j}{\longrightarrow} M_j$ is a~short exact sequence for $1 \leq j \leq m$.
Then we obtain a commutative diagram
\[
\begin{tikzcd}
\epsilon_0\colon & 0 \arrow[r] & N \arrow[r, "i"] & L \arrow[r, "p"] & M \arrow[r] & 0 \,\, \\
\epsilon_1\colon & 0 \arrow[r] & N_1 \arrow[u, hook] \arrow[r, "i_1", dashed] & L_1 \arrow[u, "{\iota_{\mathbf{c},1}}", hook] \arrow[r, "p_1", dashed] & M_1 \arrow[u, hook] \arrow[r] & 0 \,\, \\
 & & \vdots \arrow[u, hook] & \vdots \arrow[u, hook] & \vdots \arrow[u, hook] & \\
\epsilon_{m-1} & 0 \arrow[r] & N_{m-1} \arrow[u, hook] \arrow[r, "i_{m-1}", dashed] & L_{m-1} \arrow[u, hook] \arrow[r, "p_{m-1}", dashed] & M_{m-1} \arrow[u, hook] \arrow[r] & 0 \,\, \\
\epsilon_m\colon & 0 \arrow[r] & N_m \arrow[r, "i_m", dashed] \arrow[u, hook] & L_m \arrow[r, "p_m", dashed] \arrow[u, "{\iota_{\mathbf{c},m}}", hook] & M_m \arrow[u, hook] \arrow[r] & 0,
\end{tikzcd}
\]
which consists of chains of monomorphisms as columns and short exact sequences as rows. From the uniqueness of pushouts and pullbacks, the middle column $\mathbf{c}$ must be isomorphic to $B_{\mathbf{c}',\mathbf{c}''}([\epsilon_0],[\epsilon_1],\dots,[\epsilon_m])$ where $([\epsilon_1],\dots,[\epsilon_m])$ is as shown in the diagram.
\end{proof}

Let $[\epsilon]\in \EExt_{\mathcal{A}}^1(M,N)$ with ${\rm mt} \epsilon=L$, and let $\mathbf{c}' \in \mathcal{F}^{\rm mono}_{m,M}$, $\mathbf{c}''\in \mathcal{F}^{\rm mono}_{m,N}$. Then the chains $\mathbf{c}\in \mathcal{F}^{\rm mono}_{m,L}$
such that $\phi_{MN}([\epsilon],\mathbf{c})=(\mathbf{c}',\mathbf{c}'')$ are precisely those of the form $B_{\mathbf{c}',\mathbf{c}''}([\epsilon],\dots,[\epsilon_m])$ with $([\epsilon],\dots,[\epsilon_m])\in\KKer\beta_{\mathbf{c}',\mathbf{c}''}$.

On the other hand, note that $\phi_{MN}$ describes the relationship among $\mathcal{F}^{\rm mono}_{m,L}$, $\mathcal{F}^{\rm mono}_{m,M}$ and $\mathcal{F}^{\rm mono}_{m,N}$, and hence the relationship among characters $\delta_M$, $\delta_N$ and $\delta_L$. In order to describe the structure more clearly, we need to refine this map.

Fix $[\epsilon] \in \operatorname{Ext}^1_{\mathcal{A}} (M,N)$, and define
\[
\phi_{[\epsilon]}:=\phi_{MN}([\epsilon],-)\colon\ \mathcal{F}^{\rm mono}_{\mathbf{i},\mathbf{a},L}\rightarrow\coprod_{\mathbf{a}'+\mathbf{a}''=\mathbf{a}}\mathcal{F}^{\rm mono}_{\mathbf{i},\mathbf{a}',M}\times\mathcal{F}^{\rm mono}_{\mathbf{i},\mathbf{a}'',N},
\]
where $L={\rm mt}\epsilon$ and write the preimage as
\[
\mathcal{F}^{\rm mono}_{\mathbf{i},\mathbf{a},L}([\epsilon],\mathbf{a}',\mathbf{a}'')=\phi_{[\epsilon]}^{-1}(\mathcal{F}^{\rm mono}_{\mathbf{i},\mathbf{a}',M}\times\mathcal{F}^{\rm mono}_{\mathbf{i},\mathbf{a}'',N}).
\]
Since the coproduct
\[
\coprod_{\mathbf{a}'+\mathbf{a}''=\mathbf{a}}\mathcal{F}^{\rm mono}_{\mathbf{i},\mathbf{a}',M}\times\mathcal{F}^{\rm mono}_{\mathbf{i},\mathbf{a}'',N}
\]
is a disjoint union, we naturally have
\[
\mathcal{F}^{\rm mono}_{\mathbf{i},\mathbf{a},L}=\coprod_{\mathbf{a}'+\mathbf{a}''=\mathbf{a}}\mathcal{F}^{\rm mono}_{\mathbf{i},\mathbf{a},L}([\epsilon],\mathbf{a}',\mathbf{a}'').
\]
However, the structure of a fiber of
\[
\phi_{[\epsilon]}\colon\ \mathcal{F}^{\rm mono}_{\mathbf{i},\mathbf{a},L}([\epsilon],\mathbf{a}',\mathbf{a}'')\rightarrow\mathcal{F}^{\rm mono}_{\mathbf{i},\mathbf{a}',M}\times\mathcal{F}^{\rm mono}_{\mathbf{i},\mathbf{a}'',N}
\]
is heavily dependent on the relation between $[\epsilon]$ and specific chains $(\mathbf{c}',\mathbf{c}'')$.

The following corollary describes the image of $\phi_{[\epsilon]}$ restricted to $\mathcal{F}^{\rm mono}_{\mathbf{i},\mathbf{a},L}([\epsilon],\mathbf{a}',\mathbf{a}'')$.

\begin{Corollary}\label{image-related}
Given $(\mathbf{c}',\mathbf{c}'')$ in $ \mathcal{F}^{\rm mono}_{\mathbf{i},\mathbf{a}',M}\times\mathcal{F}^{\rm mono}_{\mathbf{i},\mathbf{a}'',N}$,
\[
(\mathbf{c}',\mathbf{c}'')\in \phi_{[\epsilon]}(\mathcal{F}^{\rm mono}_{\mathbf{i},\mathbf{a},L}([\epsilon],\mathbf{a}',\mathbf{a}''))\qquad \text{if and only if} \quad [\epsilon]\in p_0\KKer\beta_{\mathbf{c}',\mathbf{c}''},
\]
where $p_0$ is the projection from $\bigoplus_{j=0}^m\EExt_{\mathcal{A}}^1(M_j,N_j)$ to $\Ext_{\mathcal{A}}^1(M,N)$.
\end{Corollary}

Now, in the case that $[\epsilon]\in p_0\KKer\beta_{\mathbf{c}',\mathbf{c}''}$, we calculate the fiber $\phi^{-1}_{[\epsilon]}(\mathbf{c}',\mathbf{c}'')$.

From Lemma~\ref{chain-mono-inj}, we know that the assignment $B_{\mathbf{c}',\mathbf{c}''}$ which maps a family of extensions $([\epsilon_0],\dots,[\epsilon_m])$ to a chain of monomorphisms ending with $L={\rm mt}\epsilon_0$ is injective. We denote
\[
K(\mathbf{c}',\mathbf{c}'',[\epsilon]):=\{([\epsilon_0],\dots,[\epsilon_m])\in \KKer\beta_{\mathbf{c}',\mathbf{c}''}|[\epsilon_0]=[\epsilon]\}.
\]
Then we have
\[
\KKer\beta_{\mathbf{c}',\mathbf{c}''}=\coprod_{[\epsilon]\in p_0\KKer\beta_{\mathbf{c}',\mathbf{c}''}}K(\mathbf{c}',\mathbf{c}'',[\epsilon]).
\]
Then by Lemma~\ref{fiber-related}, we have
$\phi^{-1}_{[\epsilon]}(\mathbf{c}',\mathbf{c}'') \neq \varnothing$ only if $[\epsilon]\in p_0\KKer\beta_{\mathbf{c}',\mathbf{c}''}$ and obtain the immediate corollary

\begin{Corollary}\label{vectorbundle1}
Given $(\mathbf{c}',\mathbf{c}'')$ in $ \mathcal{F}^{\rm mono}_{\mathbf{i},\mathbf{a}',M}\times\mathcal{F}^{\rm mono}_{\mathbf{i},\mathbf{a}'',N}$,
\begin{enumerate}\itemsep=0pt
\item[$(1)$] $B_{\mathbf{c}',\mathbf{c}''}\colon K(\mathbf{c}',\mathbf{c}'',[\epsilon])\rightarrow\phi^{-1}_{[\epsilon]}(\mathbf{c}',\mathbf{c}'')$ is bijective;
\item[$(2)$] If $[\epsilon]\notin p_0\KKer\beta_{\mathbf{c}',\mathbf{c}''}$, $\phi^{-1}_{[\epsilon]}(\mathbf{c}',\mathbf{c}'')=\varnothing$.
\end{enumerate}
\end{Corollary}

\subsection{Dual case}
Recall that the abelian category $\mathcal{A}$ considered here admits $\operatorname{Ext}$-symmetry. That is to say, for any two objects $M$, $N$ in $\mathcal{A}$, there is a natural isomorphism
\[
E_{M,N}\colon\ \Ext_{\mathcal{A}}^1(M,N)\cong D\Ext_{\mathcal{A}}^1(N,M),
\]
where $D={\rm Hom}_k(-,k)$ is the linear dual. That is
\begin{Proposition}\label{nature-iso}
Given $\lambda\in\Hom(N,N')$, $\rho\in\Hom(M',M)$, $[\epsilon]\in\Ext_{\mathcal{A}}^1(M,N)$ and $[\eta]\in\Ext_{\mathcal{A}}^1(N',M')$, we have
\[
E_{M',N'}(\lambda\circ[\epsilon]\circ\rho)([\eta])=E_{M,N}([\epsilon])(\rho\circ[\eta]\circ\lambda).
\]
\end{Proposition}

Consider the linear map
\begin{align*}
\beta_{\mathbf{c}',\mathbf{c}''}\colon\ \bigoplus_{j=0}^m\Ext^1(M_j,N_j)&\longrightarrow\bigoplus_{j=0}^{m-1}\Ext^1(M_{j+1},N_j),\\
([\epsilon_0],\dots,[\epsilon_m])&\longmapsto([\epsilon_{j-1}]\circ\iota_{\mathbf{c}',j}-\iota_{\mathbf{c}'',j}\circ[\epsilon_j], 1\leq j \leq m).
\end{align*}

Since it is a linear map between finite dimensional spaces, we can calculate the dual of $\beta_{\mathbf{c}',\mathbf{c}''}$ explicitly using Proposition~\ref{nature-iso}
\begin{equation*}
\begin{aligned}
\beta'_{\mathbf{c}'',\mathbf{c}'}=D\beta_{\mathbf{c}',\mathbf{c}''}\colon\ \bigoplus_{j=0}^{m-1}\EExt^1(N_j,M_{j+1})&\longrightarrow\bigoplus_{j=0}^{m}\EExt^1(N_{j},M_j),\\
([\eta_0],\dots,[\eta_{m-1}])&\longmapsto (\iota_{\mathbf{c}',j+1}\circ[\eta_j] -[\eta_{j-1}]\circ\iota_{\mathbf{c}'',j}, 0\leq j \leq m),
\end{aligned}
\end{equation*}
where $[\eta_{-1}]:=\mathbf{0}$ and $[\eta_{m}]:=\mathbf{0}$.

To decide whether a pair of chains $(\mathbf{c}'',\mathbf{c}')$ is located in the image of $\phi_{[\eta]}$, we need
\begin{Lemma}\label{dual-1}
Given $(\mathbf{c}'',\mathbf{c}')$ in $ \mathcal{F}^{\rm mono}_{\mathbf{i},\mathbf{a}'',N}\times\mathcal{F}^{\rm mono}_{\mathbf{i},\mathbf{a}',M}$,
\[
(\mathbf{c}'',\mathbf{c}')\in \phi_{[\eta]}(\mathcal{F}^{\rm mono}_{\mathbf{i},\mathbf{a},L}([\eta],\mathbf{a}'',\mathbf{a}'))\qquad \text{if and only if}\quad [\eta]\in \operatorname{Im}\beta'_{\mathbf{c}'',\mathbf{c}'}\cap\Ext_{\mathcal{A}}^1(N,M),
\]
where the intersection is realized through regarding $\EExt_{\mathcal{A}}^1(N,M)$ as a linear subspace of $\bigoplus_{j=0}^{m}$ $\Ext_{\mathcal{A}}^1 (N_{j},M_j)$.
\end{Lemma}
\begin{proof}
By definition of $\beta'_{\mathbf{c}'',\mathbf{c}'}$, $[\eta]\in \operatorname{Im}\beta'_{\mathbf{c}'',\mathbf{c}'}\cap\EExt_{\mathcal{A}}^1(N,M)$ if and only if there exists a family of extensions $([\eta_0],\dots,[\eta_{m-1}])$ such that
\begin{enumerate}\itemsep=0pt
\item[(1)] $\iota_{\mathbf{c}',1}\circ[\eta_0]=[\eta]$;
\item[(2)] $\iota_{\mathbf{c}',j+1}\circ[\eta_j]=[\eta_{j-1}]\circ\iota_{\mathbf{c}'',j}$ for $1\leq j \leq m-1$;
\item[(3)] $[\eta_{m-1}]\circ\iota_{\mathbf{c}'',m}=0$.
\end{enumerate}

The third condition is always satisfied since $\EExt_{\mathcal{A}}^1(N_m,M_m)=\EExt_{\mathcal{A}}^1(0,0)=0$.

Since $\mathbf{a}'+\mathbf{a}''=\mathbf{a}$ is a $0$-$1$ sequence, either $\iota_{\mathbf{c}',j}$ or $\iota_{\mathbf{c}'',j}$ is an isomorphism. In either case, one can check that the first and second condition are equivalent to the condition that the pushout of $[\eta_j]$ along $\iota_{\mathbf{c}',j}$ coincides with the pullback of $[\eta_{j-1}]$ along $\iota_{\mathbf{c}'',j}$ and
the pushout of extension between $M_j$ and $N_j$ along $\iota_{M,j}$ coincides with the pullback of extension between $M_{j-1}$ and $N_{j-1}$ along $\iota_{N,j}$
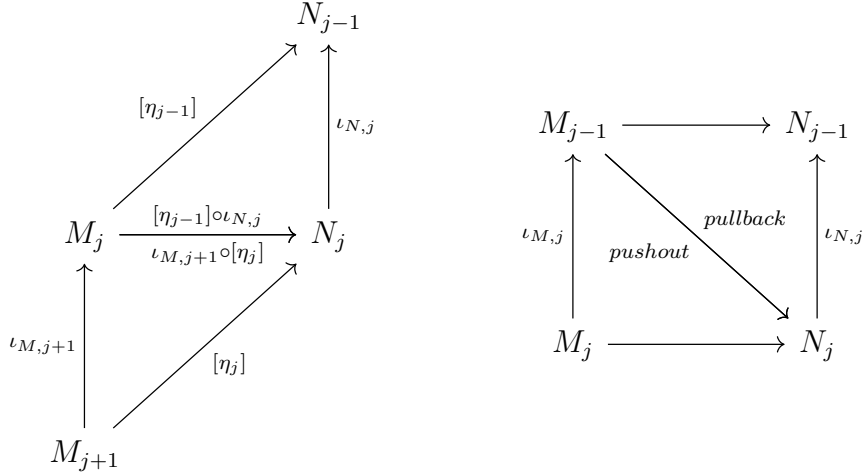
\[
\begin{tikzcd}
    & & N_{j-1} & & & & \\
                                                                                                                    & & & & M_{j-1} \arrow[rr] \arrow[rrdd, "{\rm pullback}"] \arrow[rrdd, "{\rm pushout}"'] & & N_{j-1} \,\, \\
M_j \arrow[rruu, "{[\eta_{j-1}]}"] \arrow[rr, "{\iota_{\mathbf{c}',j+1}\circ[\eta_j]}"'] \arrow[rr, "{[\eta_{j-1}]\circ\iota_{\mathbf{c}'',j}}"] & & N_j , \arrow[uu, "{\iota_{\mathbf{c}'',j}}"'] & & & & \\
                                                                                                                     & & & & M_j \arrow[rr] \arrow[uu, "{\iota_{\mathbf{c}',j}}"] & & N_j. \arrow[uu, "{\iota_{\mathbf{c}'',j}}"'] \\
M_{j+1} \arrow[rruu, "{[\eta_j]}"'] \arrow[uu, "{\iota_{\mathbf{c}',j+1}}"] & & & & & &
\end{tikzcd}
\]
According to Lemma~\ref{fiber-related}, the condition that the pullback and the pushout coincide ensures the preimage $\phi_{[\eta]}^{-1}(\mathbf{c}'',\mathbf{c}')$ is non-empty.
\end{proof}

\begin{Remark}
In fact, we can check that
\[
(\mathbf{c}'',\mathbf{c}')\in \phi_{[\eta]}(\mathcal{F}^{\rm mono}_{\mathbf{i},\mathbf{a},L}([\eta],\mathbf{a}'',\mathbf{a}')) \qquad \text{ if and only if }\quad [\eta]\in p_0\KKer \beta_{\mathbf{c}'',\mathbf{c}'}
\]
as in Corollary~\ref{image-related}.

\end{Remark}

By Lemmas~\ref{fiber-related} and~\ref{dual-1}, we obtain the dual of Corollary~\ref{vectorbundle1}.
\begin{Corollary}\label{vectorbundle2}
Given $(\mathbf{c}'',\mathbf{c}')$ in $ \mathcal{F}^{\rm mono}_{\mathbf{i},\mathbf{a}'',N}\times\mathcal{F}^{\rm mono}_{\mathbf{i},\mathbf{a}',M}$,
\begin{enumerate}\itemsep=0pt
\item[$(1)$] $B_{\mathbf{c}'',\mathbf{c}'}\colon K(\mathbf{c}'',\mathbf{c}',[\eta])\rightarrow\phi^{-1}_{[\eta]}(\mathbf{c}'',\mathbf{c}')$ is bijective;
\item[$(2)$] If $[\eta]\notin \operatorname{Im}\beta'_{\mathbf{c}'',\mathbf{c}'}\cap\EExt_{\mathcal{A}}^1(N,M)$, $\phi^{-1}_{[\eta]}(\mathbf{c}_N,\mathbf{c}_M)=\varnothing$.
\end{enumerate}
\end{Corollary}

\subsection{Cardinality}
Based on several lemmas above, we can refine the calculation of $\delta_L$. Recall that all extension groups are finite dimensional over a finite field. Since
\[
\mathcal{F}^{\rm mono}_{\mathbf{i},\mathbf{a},L}=\coprod_{\mathbf{a}'+\mathbf{a}''=\mathbf{a}}\mathcal{F}^{\rm mono}_{\mathbf{i},\mathbf{a},L}([\epsilon],\mathbf{a}',\mathbf{a}''),
\]
we can write
\[
|\mathcal{F}^{\rm mono}_{\mathbf{i},\mathbf{a},L}|=\int_{\mathbf{a}'+\mathbf{a}''=\mathbf{a}}|\mathcal{F}^{\rm mono}_{\mathbf{i},\mathbf{a},L}([\epsilon],\mathbf{a}',\mathbf{a}'')|.
\]
Note that the decomposition of $\mathcal{F}^{\rm mono}_{\mathbf{i},\mathbf{a},L}$ depends on the choice of the short exact sequence $\epsilon$.

By linear algebra, $|K(\mathbf{c}',\mathbf{c}'',[\epsilon])|=|K(\mathbf{c}',\mathbf{c}'',{\bf 0})|$ if $[\epsilon]\in p_0\KKer \beta_{\mathbf{c}',\mathbf{c}''}$. Note that $K(\mathbf{c}',\mathbf{c}'',{\bf 0})$ is a vector space over $k$ so we can
define $k(\mathbf{c}',\mathbf{c}''):=\operatorname{dim}_k K(\mathbf{c}',\mathbf{c}'',{\bf 0})$.

Moreover, Corollary~ \ref{vectorbundle1} shows that if $[\epsilon]\notin p_0\KKer \beta_{\mathbf{c}',\mathbf{c}''}$, then $|\phi^{-1}_{[\epsilon]}(\mathbf{c}',\mathbf{c}'')|=0$, and otherwise
 $|\phi^{-1}_{[\epsilon]}(\mathbf{c}',\mathbf{c}')|=|K(\mathbf{c}',\mathbf{c}'',[\epsilon])|=q^{k(\mathbf{c}',\mathbf{c}'')}$. So if ${\rm mt}\epsilon=L$, we have
\begin{align*}
\delta_L(d_{\mathbf{i},\mathbf{a}})&=|\mathcal{F}^{\rm mono}_{\mathbf{i},\mathbf{a},L}|=\int_{\mathbf{a}'+\mathbf{a}''=\mathbf{a}}|\mathcal{F}^{\rm mono}_{\mathbf{i},\mathbf{a},L}([\epsilon],\mathbf{a}',\mathbf{a}'')|\\
&=\iint_{\mathbf{a}'+\mathbf{a}''=\mathbf{a}, (\mathbf{c}',\mathbf{c}'')\in\mathcal{F}^{\rm mono}_{\mathbf{i},\mathbf{a}',M}\times\mathcal{F}^{\rm mono}_{\mathbf{i},\mathbf{a}'',N}}\big|\phi^{-1}_{[\epsilon]}(\mathbf{c}',\mathbf{c}'')\big|\\
&=\iint_{\mathbf{a}'+\mathbf{a}''=\mathbf{a}, (\mathbf{c}',\mathbf{c}'')\in\phi_{[\epsilon]}(\mathcal{F}^{\rm mono}_{\mathbf{i},\mathbf{a},L}([\epsilon],\mathbf{a}',\mathbf{a}''))}q^{k(\mathbf{c}',\mathbf{c}'')},
\end{align*}
where all integrals are finite sums.

\begin{Remark}
Since different $[\epsilon]$, $[\epsilon']$ in $\EExt_{\mathcal{A}}^1(M,N)$ may admit the same middle term ${\rm mt}\epsilon={\rm mt}\epsilon'=L$, the refinement of $\delta_L(d_{\mathbf{i},\mathbf{a}})$ in the above equation depends on
the choice of $[\epsilon]$.
\end{Remark}

Moreover, Corollary~\ref{vectorbundle2} shows if $[\eta]\notin \operatorname{Im}\beta'_{\mathbf{c}'', \mathbf{c}'}\cap\EExt_{\mathcal{A}}^1(N,M)$, $|\phi^{-1}_{[\eta]}(\mathbf{c}'',\mathbf{c}'')|=0$ and if $[\eta]\in \operatorname{Im}\beta'_{\mathbf{c}'', \mathbf{c}'}\cap\EExt_{\mathcal{A}}^1(N,M)$,
$|\phi^{-1}_{[\eta]}(\mathbf{c}'',\mathbf{c}')|=|K(\mathbf{c}'',\mathbf{c}',[\eta])|=q^{k(\mathbf{c}'',\mathbf{c}')}$. So if ${\rm mt}\eta=L$, we have
\begin{align*}
\delta_L(d_{\mathbf{i},\mathbf{a}})&=|\mathcal{F}^{\rm mono}_{\mathbf{i},\mathbf{a},L}|=\int_{\mathbf{a}''+\mathbf{a}'=\mathbf{a}}|\mathcal{F}^{\rm mono}_{\mathbf{i},\mathbf{a},L}([\eta],\mathbf{a}'',\mathbf{a}')|\\
&=\iint_{\mathbf{a}''+\mathbf{a}'=\mathbf{a}, (\mathbf{c}'',\mathbf{c}')\in\mathcal{F}^{\rm mono}_{\mathbf{i},\mathbf{a}'',N}\times\mathcal{F}^{\rm mono}_{\mathbf{i},\mathbf{a}',M}}\big|\phi^{-1}_{[\eta]}(\mathbf{c}'',\mathbf{c}')\big|\\
&=\iint_{\mathbf{a}''+\mathbf{a}'=\mathbf{a}, (\mathbf{c}'',\mathbf{c}'')\in\phi_{[\eta]}(\mathcal{F}^{\rm mono}_{\mathbf{i},\mathbf{a},L}([\eta],\mathbf{a}'',\mathbf{a}'))}q^{k(\mathbf{c}'',\mathbf{c}')}.
\end{align*}

\subsection{Weight}
Based on the calculation of $\delta_L$, we introduce the notion of weight functions and weighted quantum cluster functions. 
Denote $\mathcal{F}^{\rm mono}:=\bigcup_{m\in \mathbb{N}} \mathcal{F}^{\rm mono}_m$.

We define
\[
\operatorname{MF}:=\mathcal{F}^{\rm mono} \times \mathcal{F}^{\rm mono}=\{(\mathbf{c}',\mathbf{c}'')|\mathbf{c}'\in\mathcal{F}^{\rm mono}_{m,M},\,\mathbf{c}''\in\mathcal{F}^{\rm mono}_{m, N},\, M, N\in\mathcal{A},\, m \in \mathbb{N}\}
\]
and set
\[
\mathbb{Z}_{\operatorname{MF}}:=\big\{ f\colon \operatorname{MF} \times \EExt_{\mathcal{A}}^1 \rightarrow Z| f(\mathbf{c}', \mathbf{c}'', [\epsilon])=0 \text{ unless }\mathbf{c}'_0={\rm qt}\epsilon, \, \mathbf{c}''_0={\rm st}\epsilon \big\},
\]
where \smash{$\EExt_{\mathcal{A}}^1=\coprod_{M, N \in \mathcal{A}}\EExt_{\mathcal{A}}^1 (M,N)$} and $Z=\big\{\frac{n}{2}\mid n\in\mathbb{Z}\big\}$ is the set of all half integers. The functions in $\mathbb{Z}_{\operatorname{MF}}$ are called weight functions.
Given $[\epsilon] \in \Ext_{\mathcal{A}}^1 (M, N)$, we define
\[\mathbb{Z}_{\operatorname{MF}}[\epsilon]:= \{f\in \mathbb{Z}_{\operatorname{MF}}\mid f(\mathbf{c}', \mathbf{c}'', [\rho])=0 \text{ if } [\rho]\neq [\epsilon]\}.\]
Given $f\in \mathbb{Z}_{\operatorname{MF}}[\epsilon]$, we write $f(\mathbf{c}', \mathbf{c}'', [\epsilon])$ instead as $f(\mathbf{c}', \mathbf{c}'')$.

\begin{Definition}[weighted quantum cluster function]
Given a weight function $f \in\mathbb{Z}_{\operatorname{MF}}[\epsilon]$, the weighted quantum cluster function $f*_{[\epsilon]} \delta_L$ is the linear function on $\mathcal{M}_q$ defined by
\begin{align*}
f*_{[\epsilon]}\delta_L(d_{\mathbf{i},\mathbf{a}})&=\iint_{\mathbf{a}'+\mathbf{a}''=\mathbf{a}, (\mathbf{c}',\mathbf{c}'')\in\mathcal{F}^{\rm mono}_{\mathbf{i},\mathbf{a}',M}\times\mathcal{F}^{\rm mono}_{\mathbf{i},\mathbf{a}'',N}}\big|\phi^{-1}_{[\epsilon]}(\mathbf{c}',\mathbf{c}'')\big|\cdot q^{f(\mathbf{c}',\mathbf{c}'')}\\
&=\iint_{\mathbf{a}'+\mathbf{a}''=\mathbf{a}, (\mathbf{c}',\mathbf{c}'')\in\phi_{[\epsilon]}(\mathcal{F}^{\rm mono}_{\mathbf{i},\mathbf{a},L}([\epsilon],\mathbf{a}',\mathbf{a}''))}q^{k(\mathbf{c}',\mathbf{c}'')}\cdot q^{f (\mathbf{c}',\mathbf{c}'')},
\end{align*}
where $M={\rm qt}\epsilon$, $N={\rm st}\epsilon$ and $L={\rm mt}\epsilon$.
\end{Definition}
If $f\in \mathbb{Z}_{\MF}$ is the zero function, we have $f *_{[\epsilon]}\delta_L=\delta_L$.
So weight functions provide $q$-deformations of $\delta_L$.

\subsection{Multiplication}
We now introduce the multiplication of weighted quantum cluster functions.
\begin{Definition}[multiplication of quantum cluster functions]
Given $M$, $N$ in $\mathcal{A}$, we define the multiplication of quantum cluster functions as
\[
\delta_M*\delta_N=\delta_{M\oplus N}.
\]
\end{Definition}
By definition, we have
\begin{Proposition}
The multiplication $*$ is associative and commutative.
\end{Proposition}

We denote the zero element in $\EExt_{\mathcal{A}}^1(M,N)$ by $\mathbf{0}_{MN}$.
\begin{Remark}
From the refinement of the quantum cluster function, we have
\[\delta_{M\oplus N}(d_{\mathbf{i},\mathbf{a}})={\hskip -0.2cm}\int_{\mathbf{a}'+\mathbf{a}''=\mathbf{a}}|\mathcal{F}^{\rm mono}_{\mathbf{i},\mathbf{a},M\oplus N}(\mathbf{0}_{MN},\mathbf{a}',\mathbf{a}'')|
={\hskip -0.2cm}\iint_{\mathbf{a}'+\mathbf{a}''=\mathbf{a}, (\mathbf{c}',\mathbf{c}'')\in\mathcal{F}^{\rm mono}_{\mathbf{i},\mathbf{a}',M}\times\mathcal{F}^{\rm mono}_{\mathbf{i},\mathbf{a}'',N}}q^{k(\mathbf{c}',\mathbf{c}'')}.\]

Notice that the image of $\phi_{\mathbf{0}_{MN}}$ is the whole of $\mathcal{F}^{\rm mono}_{\mathbf{i},\mathbf{a}',M}\times\mathcal{F}^{\rm mono}_{\mathbf{i},\mathbf{a}'',N}$ since any two chains $\mathbf{c}'$, $\mathbf{c}''$ can be assembled into a chain ending with $M\oplus N$ through direct sum.

On the other hand, the convolution of the two functions $\delta_M$ and $\delta_N$ can be formally defined~as
\begin{align*}
\int_{\mathbf{a}'+\mathbf{a}''=\mathbf{a}}\delta_M(d_{\mathbf{i},\mathbf{a}'})\cdot\delta_N(d_{\mathbf{i},\mathbf{a}''})
&=\int_{\mathbf{a}'+\mathbf{a}''=\mathbf{a}}|\mathcal{F}^{\rm mono}_{\mathbf{i},\mathbf{a}',M}|\cdot|\mathcal{F}^{\rm mono}_{\mathbf{i},\mathbf{a}'',N}|\\
&=\iint_{\mathbf{a}'+\mathbf{a}''=\mathbf{a}, (\mathbf{c}',\mathbf{c}'')\in\mathcal{F}^{\rm mono}_{\mathbf{i},\mathbf{a}',M}\times\mathcal{F}^{\rm mono}_{\mathbf{i},\mathbf{a}'',N}}1.
\end{align*}
So the definition of multiplication is a $q$-deformation of conventional convolution.
\end{Remark}

\begin{Definition}[multiplication of weighted quantum cluster functions]\label{associativity}
Given weighted quantum cluster functions $f*_{[\epsilon']}\delta_M$ and $g*_{[\epsilon'']}\delta_N$, we define the multiplication as
\[
(f*_{[\epsilon']}\delta_M)*(g*_{[\epsilon'']}\delta_N)=h*_{\mathbf{0}_{MN}}\delta_{M\oplus N},
\]
where $h(\mathbf{c}', \mathbf{c}'', [\epsilon])=0$ unless $[\epsilon]={\mathbf{0}_{MN}}$ and $h( \mathbf{c}', \mathbf{c}'') =f (\phi_{[\epsilon']} (\mathbf{c}'))+g (\phi_{[\epsilon'']} (\mathbf{c}''))$.
\end{Definition}

\begin{Proposition}
The multiplication operation from Definition~{\rm \ref{associativity}} is associative.
\end{Proposition}
\begin{proof}
Given $f_i \in\mathbb{Z}_{\MF}[{\epsilon_i}]$ such that $\operatorname{mt}\epsilon_i=M_i$ for $i=1,2,3$, we set
\begin{gather*}
((f_1 *_{[\epsilon_1]}\delta_{M_1})*(f_2 *_{[\epsilon_2]}\delta_{M_2}))*(f_3 *_{[\epsilon_3]}\delta_{M_3})=f_{(12)3}*_{\mathbf{0}_{M_1\oplus M_2,M3}}\delta_{M_1\oplus M_2\oplus M_3},\\
(f_1 *_{[\epsilon_1]}\delta_{M_1})*((f_2 *_{[\epsilon_2]}\delta_{M_2})*(f_3 *_{[\epsilon_3]}\delta_{M_3}))=f_{1(23)}*_{\mathbf{0}_{M_1, M_2\oplus M3}}\delta_{M_1\oplus M_2\oplus M_3},
\end{gather*}
where $f_{(12)3} \in \mathbb{Z}_{\MF}\mathbf{0}_{M_1\oplus M_2,M3}$, $f_{1(23)} \in\mathbb{Z}_{\MF}\mathbf{0}_{M_1, M_2\oplus M3}$.

Note that the non-weighted parts of the two compositions are the same.
We have \[(f_1 *_{[\epsilon_1]}\delta_{M_1})*(f_2 *_{[\epsilon_2]}\delta_{M_2})=f_{(12)} *_{\mathbf{0}_{M_1\oplus M_2}} \delta_{M_1\oplus M_2},\]
where $f_{(12)}\in \mathbb{Z}_{\MF}\mathbf{0}_{MN}$ and
\[f_{(12)} (\mathbf{c}_1, \mathbf{c}_2)=f_1 (\phi_{[\epsilon_1]} (\mathbf{c}_1))+f_2 (\phi_{[\epsilon_2]} (\mathbf{c}_2)).\]

Then
\begin{align*}
f_{(12)3} (\mathbf{c}_1 \oplus \mathbf{c}_2, \mathbf{c}_3)&= f_{(12)} (\mathbf{c}_1, \mathbf{c}_2) + f_3 (\phi_{[\epsilon_3]} (\mathbf{c}_3))\\
&=f_1 (\phi_{[\epsilon_1]} (\mathbf{c}_1))+f_2 (\phi_{[\epsilon_2]} (\mathbf{c}_2)) + f_3 (\phi_{[\epsilon_3]} (\mathbf{c}_3)).
\end{align*}
Similarly, we have that
\begin{align*}
f_{1(23)}(\mathbf{c}_1, \mathbf{c}_2 \oplus \mathbf{c}_3)&= f_1 (\phi_{[\epsilon_1]} (\mathbf{c}_1)) + f_{(23)} (\mathbf{c}_2, \mathbf{c}_3)\\
&=f_1 (\phi_{[\epsilon_1]} (\mathbf{c}_1))+f_2 (\phi_{[\epsilon_2]} (\mathbf{c}_2)) + f_3 (\phi_{[\epsilon_3]} (\mathbf{c}_3)).
\end{align*}

Then we have
\begin{gather*}
f_{(12)3}*_{\mathbf{0}_{M_1 \oplus M_2 \oplus M_3}}\delta_{M_1\oplus M_2\oplus M_3}(d_{\mathbf{i},\mathbf{a}})\\
\quad=\iint_{\mathbf{a}'+\mathbf{a}''+\mathbf{a}'''=\mathbf{a}, (\mathbf{c}_{1},\mathbf{c}_{2},\mathbf{c}_{3})}q^{k_{(12)3}(\mathbf{c}_{1},\mathbf{c}_{2},\mathbf{c}_{3})}\cdot q^{
f_1 (\phi_{[\epsilon_1]} (\mathbf{c}_1))+f_2 (\phi_{[\epsilon_2]} (\mathbf{c}_2)) + f_3 (\phi_{[\epsilon_3]} (\mathbf{c}_3))}
\end{gather*}
and
\begin{gather*}
f_{1(23)}*_{\mathbf{0}_{M_1 \oplus M_2 \oplus M_3}} \delta_{M_1\oplus M_2\oplus M_3}(d_{\mathbf{i},\mathbf{a}})\\
\quad=\iint_{\mathbf{a}'+\mathbf{a}''+\mathbf{a}'''=\mathbf{a}, (\mathbf{c}_{1},\mathbf{c}_{2},\mathbf{c}_{3})}q^{k_{1(23)}(\mathbf{c}_{1},\mathbf{c}_{2},\mathbf{c}_{3})}\cdot q^{
f_1 (\phi_{[\epsilon_1]} (\mathbf{c}_1))+f_2 (\phi_{[\epsilon_2]} (\mathbf{c}_2)) + f_3 (\phi_{[\epsilon_3]} (\mathbf{c}_3))},
\end{gather*}
where the second integral in both equations runs over \[\mathcal{F}^{\rm mono}_{\mathbf{i},\mathbf{a}',M_1}\times\mathcal{F}^{\rm mono}_{\mathbf{i},\mathbf{a}'',M_2}\times\mathcal{F}^{\rm mono}_{\mathbf{i},\mathbf{a}''',M_3}\]
 and the functions $k_{(12)3}$ and $k_{1(23)}$
compute the dimensions of fibers of the diagonal maps $\phi_{(12)3}$ and~$\phi_{1(23)}$ respectively in the diagram
\[
\begin{tikzcd}[sep=1.4cm]
{\mathcal{F}^{\rm mono}_{\mathbf{i},\mathbf{a},M_1\oplus M_2\oplus M_3}} \arrow[dd, "{\phi_{M_1,M_2\oplus M_3,\mathbf{0}_{M_1,M_2\oplus M_3}}}"'] \arrow[rr, "{\phi_{M_1\oplus M_2,M_3,\mathbf{0}_{M_1\oplus M_2,M_3}}}"] \arrow[rrdd, "\phi_{(12)3}", dashed] \arrow[rrdd, "\phi_{1(23)}"', dashed] & & {\coprod\mathcal{F}^{\rm mono}_{\mathbf{i},\mathbf{a}'+\mathbf{a}'',M_1\oplus M_2}\times\mathcal{F}^{\rm mono}_{\mathbf{i},\mathbf{a}'',M_3}} \arrow[dd, "{\coprod \phi_{M_1,M_2,\mathbf{0}_{M_1,M_2}}\times {\rm id}}"] \\
                                                                                                                                                                                                                                                                                    & &  \\
{\coprod\mathcal{F}^{\rm mono}_{\mathbf{i},\mathbf{a}',M_1}\times\mathcal{F}^{\rm mono}_{\mathbf{i},\mathbf{a}''+\mathbf{a}''',M_2\oplus M_3}} \arrow[rr, "{\coprod{\rm id}\times\phi_{M_2,M_3,\mathbf{0}_{M_2,M_3}}}"'] & & {\coprod\mathcal{F}^{\rm mono}_{\mathbf{i},\mathbf{a}',M_1}\times\mathcal{F}^{\rm mono}_{\mathbf{i},\mathbf{a}'',M_2}\times\mathcal{F}^{\rm mono}_{\mathbf{i},\mathbf{a}''',M_3}}.
\end{tikzcd}
\]
From the associativity of $\delta_M*\delta_N$, we know that the diagram commutes. So $k_{(12)3}=k_{1(23)}$.
\end{proof}

\begin{Remark}
By definition, $\delta_M*\delta_N=\delta_N*\delta_M$. This implies
\[\iint_{\mathbf{a}'+\mathbf{a}''=\mathbf{a}, (\mathbf{c}',\mathbf{c}'')\in\mathcal{F}^{\rm mono}_{\mathbf{i},\mathbf{a}',M}\times\mathcal{F}^{\rm mono}_{\mathbf{i},\mathbf{a}'',N}}q^{k(\mathbf{c}',\mathbf{c}'')}
=\iint_{\mathbf{a}'+\mathbf{a}''=\mathbf{a}, (\mathbf{c}',\mathbf{c}'')\in\mathcal{F}^{\rm mono}_{\mathbf{i},\mathbf{a}',M}\times\mathcal{F}^{\rm mono}_{\mathbf{i},\mathbf{a}'',N}}q^{k(\mathbf{c}'',\mathbf{c}')}.\]
However, it need not be the case
that $k(\mathbf{c}',\mathbf{c}'')=k(\mathbf{c}'',\mathbf{c}')$, since one value is computed using~$\mathbf{0}_{MN}$ and the other using
$\mathbf{0}_{NM}$.

So even though the weight functions $h$ and $h'$ in
\begin{gather*}
(f *_{[\epsilon']} \delta_M)*(g*_{[\epsilon'']}\delta_N)=h *_{\mathbf{0}_{MN}}\delta_{M\oplus N},\\
(g*_{[\epsilon'']}\delta_N) *(f *_{[\epsilon']} \delta_M) =h'*_{\mathbf{0}_{NM}}*\delta_{M\oplus N}
\end{gather*}
satisfy $h(\mathbf{c}',\mathbf{c}'')=h'(\mathbf{c}'',\mathbf{c}')$ by definition, we still can not deduce that
\begin{gather*}
\iint_{\mathbf{a}'+\mathbf{a}''=\mathbf{a}, (\mathbf{c}',\mathbf{c}'')\in\mathcal{F}^{\rm mono}_{\mathbf{i},\mathbf{a}',M}\times\mathcal{F}^{\rm mono}_{\mathbf{i},\mathbf{a}'',N}}q^{k(\mathbf{c}',\mathbf{c}'')}\cdot q^{h (\mathbf{c}',\mathbf{c}'')}\\
\quad=\iint_{\mathbf{a}'+\mathbf{a}''=\mathbf{a}, (\mathbf{c}',\mathbf{c}'')\in\mathcal{F}^{\rm mono}_{\mathbf{i},\mathbf{a}',M}\times\mathcal{F}^{\rm mono}_{\mathbf{i},\mathbf{a}'',N}}q^{k(\mathbf{c}'',\mathbf{c}')}\cdot q^{h'(\mathbf{c}'',\mathbf{c}')}.
\end{gather*}
\end{Remark}

\subsection[The projectivization of Ext\_A\^1(M,N)]{The projectivization of $\boldsymbol{\Ext_{\mathcal{A}}^1(M,N)}$}

Since $\Ext_{\mathcal{A}}^1(M,N)$ is a finite dimensional vector space, we can consider $\mathbb{P}\Ext_{\mathcal{A}}^1(M,N)$. We denote the equivalence class of $[\epsilon]$ in $\mathbb{P}\Ext_{\mathcal{A}}^1(M,N)$ by $\mathbb{P}[\epsilon]$.

The core aim in this subsection is to check that the multiplication of weighted quantum cluster functions is still well defined if we replace $[\epsilon]$ by $\mathbb{P}[\epsilon]$.
In this subsection, we fix a non-zero parameter $\lambda$ in the field $k$.
Recall the mapping
\[
\phi_{MN}\colon\ \EF_{\mathbf{i},\mathbf{a}}(M,N)\rightarrow\coprod_{\mathbf{a}'+\mathbf{a}''=\mathbf{a}}\mathcal{F}^{\rm mono}_{\mathbf{i},\mathbf{a}',M}\times\mathcal{F}^{\rm mono}_{\mathbf{i},\mathbf{a}'',N}
\]
with affine fibers. By Lemma~\ref{fiber-related} and the linearity of $\beta_{\mathbf{c}',\mathbf{c}''}$, $[\epsilon]\in p_0\KKer\beta_{\mathbf{c}',\mathbf{c}''}$ if and only if~$\lambda[\epsilon]\in p_0\operatorname{Ker}\beta_{\mathbf{c}',\mathbf{c}''}$. So in this case,
\[
\big|\phi^{-1}_{[\epsilon]}(\mathbf{c}',\mathbf{c}'')\big|=|K(\mathbf{c}',\mathbf{c}'',[\epsilon])|=|K(\mathbf{c}',\mathbf{c}'',\lambda[\epsilon])|=\big|\phi^{-1}_{\lambda[\epsilon]}(\mathbf{c}',\mathbf{c}'')\big|.
\]
Otherwise, they are all zero. So we have
\begin{Proposition}
If $f, f'\in\mathbb{Z}_{\MF}$ satisfy
\[
f(\mathbf{c}', \mathbf{c}'', [\epsilon])=f'(\mathbf{c}', \mathbf{c}'', \lambda[\epsilon])\qquad \text{for all}\quad (\mathbf{c}',\mathbf{c}'') \in\mathcal{F}^{\rm mono}_{\mathbf{i},\mathbf{a}',M}\times\mathcal{F}^{\rm mono}_{\mathbf{i},\mathbf{a}'',N},\]
then $f*_{[\epsilon]}\delta_L=f*_{\lambda[\epsilon]}\delta_L$.
\end{Proposition}
\begin{proof}
\begin{align*}
f*_{[\epsilon]}\delta_L(d_{\mathbf{i},\mathbf{a}})&=\iint_{\mathbf{a}'+\mathbf{a}''=\mathbf{a}, (\mathbf{c}',\mathbf{c}'')\in\mathcal{F}^{\rm mono}_{\mathbf{i},\mathbf{a}',M}\times\mathcal{F}^{\rm mono}_{\mathbf{i},\mathbf{a}'',N}}|\phi^{-1}_{[\epsilon]}(\mathbf{c}',\mathbf{c}'')|\cdot q^{f(\mathbf{c}',\mathbf{c}'', [\epsilon])}\\
&=\iint_{\mathbf{a}'+\mathbf{a}''=\mathbf{a}, (\mathbf{c}',\mathbf{c}'')\in\mathcal{F}^{\rm mono}_{\mathbf{i},\mathbf{a}',M}\times\mathcal{F}^{\rm mono}_{\mathbf{i},\mathbf{a}'',N}}|\phi^{-1}_{\lambda[\epsilon]}(\mathbf{c}',\mathbf{c}'')|\cdot q^{f(\mathbf{c}',\mathbf{c}'', \lambda[\epsilon])}\\
&=f*_{\lambda[\epsilon]}\delta_L(d_{\mathbf{i},\mathbf{a}}). \tag*{\qed}
\end{align*}\renewcommand{\qed}{}
\end{proof}

Recall that given $f*_{[\epsilon']}\delta_M$ and $g*_{[\epsilon'']}\delta_N$, we define
\[
(f*_{[\epsilon']}\delta_M)*(g*_{[\epsilon'']}\delta_N)=h*_{\mathbf{0}_{MN}}\delta_{M\oplus N},
\]
where $h( \mathbf{c}', \mathbf{c}'') =f (\phi_{[\epsilon']} (\mathbf{c}'))+g (\phi_{[\epsilon'']} (\mathbf{c}''))$. Since $\phi_{[\epsilon]} (\mathbf{c})=\phi_{\lambda [\epsilon]} (\mathbf{c})$ for the corresponding $[\epsilon]$ and~$\mathbf{c}$, we have
\[
(f*_{[\epsilon']}\delta_M)*(g*_{[\epsilon'']}\delta_N)=(f*_{\lambda[\epsilon']}\delta_M)*(g*_{\mu[\epsilon'']}\delta_N)\]
for any non-zero $\lambda$ and $\mu$ in $k$.

\subsection{Multiplication of weight functions}
\begin{Definition}
Given weight functions
 $f, g \in \mathbb{Z}_{\operatorname{MF}}$, $f*_{[\eta]}g$ is defined by
 \begin{equation*}
f*_{[\eta]}g (\mathbf{c}' ,\mathbf{c}'', [\epsilon]) =\left\{
 \begin{aligned}
& f(\mathbf{c}',\mathbf{c}'', [\epsilon])+g(\mathbf{c}',\mathbf{c}'', [\epsilon]), &~\text{if}~[\epsilon]=[\eta],\\
& 0 , &~\text{otherwise}.
 \end{aligned}
 \right.
\end{equation*}
 In particular, $f*_{[\eta]}g \in \mathbb{Z}_{\operatorname{MF}} [\eta]$.
\end{Definition}

\begin{Definition}
Given $(f*_{[\epsilon']}\delta_M)*(g *_{[\epsilon'']}\delta_N)=h*_{\mathbf{0}_{M\oplus N}}\delta_{M\oplus N}$ and $[\epsilon]\in\EExt_{\mathcal{A}}^1(M,N)$, we define $\mathbb{S}_{fg}\in \mathbb{Z}_{\MF} [\epsilon]$ by
\[\mathbb{S}_{fg} (\mathbf{c}',\mathbf{c}'', [\epsilon])=h( \mathbf{c}', \mathbf{c}'') =f (\phi_{[\epsilon']} (\mathbf{c}'))+g (\phi_{[\epsilon'']} (\mathbf{c}'')).\]
\end{Definition}
\begin{Corollary}
Given $f *_{[\epsilon]'}\delta_M$ and $g*_{[\epsilon]''}\delta_N$, for any $[\epsilon]\in\EExt_{\mathcal{A}}^1(M,N)$ and $[\eta]\in\EExt_{\mathcal{A}}^1(N,M)$, we have
\[
\mathbb{S}_{fg} (\mathbf{c}',\mathbf{c}'', [\epsilon])=\mathbb{S}_{gf}(\mathbf{c}'',\mathbf{c}', [\eta]).
\]
\end{Corollary}
\begin{proof}
Both sides are equal to $f (\phi_{[\epsilon']} (\mathbf{c}'))+g (\phi_{[\epsilon'']} (\mathbf{c}''))$.
\end{proof}

\subsection{Multiplication formula and balanced pairs}\label{sectionofmaintheorem1}
Recall the linear map
\begin{align*}
\beta_{\mathbf{c}',\mathbf{c}''}\colon\ \bigoplus_{j=0}^m\Ext^1(M_j,N_j)&\longrightarrow\bigoplus_{j=0}^{m-1}\EExt^1(M_{j+1},N_j),\\
([\epsilon_0],\dots,[\epsilon_m])&\longmapsto([\epsilon_{j-1}]\circ\iota_{\mathbf{c}',j}-\iota_{\mathbf{c}'',j}\circ[\epsilon_j], 1\leq j \leq m)
\end{align*}
and its dual
\begin{align*}
\beta'_{\mathbf{c}'',\mathbf{c}'}=D\beta_{\mathbf{c}',\mathbf{c}''}\colon\ \bigoplus_{j=0}^{m-1}\EExt^1(N_j,M_{j+1})&\longrightarrow\bigoplus_{j=0}^{m}\EExt^1(N_{j},M_j)\\
([\eta_0],\dots,[\eta_{m-1}])&\longmapsto (\iota_{\mathbf{c}',j+1}\circ[\eta_j] -[\eta_{j-1}]\circ\iota_{\mathbf{c}'',j}, 0\leq j \leq m),
\end{align*}
where $[\eta_{-1}]:=\mathbf{0}$ and $[\eta_{m}]:=\mathbf{0}$.

So we have
\begin{Lemma}
For any chains of morphisms $\mathbf{c}'$, $\mathbf{c}''$ ending with $M$, $N$ respectively, $\Ker\beta_{\mathbf{c}',\mathbf{c}''}=(\operatorname{Im}\beta'_{\mathbf{c}'',\mathbf{c}'})^{\perp}$. In particular,
\[
\dim_k p_0\KKer\beta_{\mathbf{c}',\mathbf{c}''}+\dim_k \bigl(\operatorname{Im}\beta'_{\mathbf{c}'',\mathbf{c}'}\cap\Ext_{\mathcal{A}}^1(N,M)\bigr)=\dim_k \EExt_{\mathcal{A}}^1(M,N).
\]
\end{Lemma}

Based on this lemma, we introduce several special weight functions.
\begin{Definition}
There are three weight functions in $\mathbb{Z}_{\MF}$ defined as
\begin{enumerate}\itemsep=0pt
\item[(1)] $f_{\rm hom}(\mathbf{c}'',\mathbf{c}', [\eta])=k(\mathbf{c}',\mathbf{c}'')-k(\mathbf{c}'',\mathbf{c}')$;
\item[(2)] $f^+_{\rm ext}(\mathbf{c}',\mathbf{c}'', [\epsilon])=\dim_k \bigl(\operatorname{Im}\beta'_{\mathbf{c}'',\mathbf{c}'}\cap\EExt_{\mathcal{A}}^1(N,M)\bigr)$;
\item[(3)] $f^-_{\rm ext}(\mathbf{c}'',\mathbf{c}', [\eta])=\dim_k p_0\KKer\beta_{\mathbf{c}',\mathbf{c}''}$
\end{enumerate}
for any $M$, $N$ in $\mathcal{A}$, $\mathbf{c}' \in\mathcal{F}^{\rm mono}_{m,M}$, $\mathbf{c}''\in \mathcal{F}^{\rm mono}_{m,N}$, $[\epsilon] \in \EExt_{\mathcal{A}}^1(M,N)$ and
$[\eta] \in \EExt_{\mathcal{A}}^1(N,M)$.
Note that all of these functions are constant in the extension.
\end{Definition}

\begin{Definition}\label{forproof-balance}
Given a pair of weight functions $(f^+,f^-)$, set \begin{gather*}
\sigma_1\bigl(f^+\bigr):=\int_{\mathbb{P}\epsilon\in\mathbb{P}p_0\operatorname{Ker}\beta_{\mathbf{c}',\mathbf{c}''}}q^{k(\mathbf{c}',\mathbf{c}'')+f^+(\mathbf{c}',\mathbf{c}'', [\epsilon])},\\
\sigma_2(f^-):=\int_{\mathbb{P}\eta\in\mathbb{P}\operatorname{Im}\beta'_{\mathbf{c}',\mathbf{c}''}\cap\operatorname{Ext}_{\mathcal{A}}^1(N,M)}q^{k(\mathbf{c}'',\mathbf{c}')+f^- (\mathbf{c}'',\mathbf{c}', [\eta])}.\end{gather*}
This pair is called pointwise balanced if \[\frac{q^{\operatorname{dim}_k \operatorname{Ext}_{\mathcal{A}}^1(M,N)}-1}{q-1}\cdot q^{k(\mathbf{c}',\mathbf{c}'')}=\sigma_1\bigl(f^+\bigr)+\sigma_2(f^-)\]
holds for any $M,N\in\mathcal{A}$ and $(\mathbf{c}',\mathbf{c}'')\in\mathcal{F}^{\rm mono}_{\mathbf{i},\mathbf{a}',M}\times\mathcal{F}^{\rm mono}_{\mathbf{i},\mathbf{a}'',N}$.
\end{Definition}

\begin{Proposition}\label{mainlemma11}
The following two pairs of weighted functions
\begin{enumerate}\itemsep=0pt
\item[$(1)$] $\bigl(f^+_{\rm ext},f_{\rm hom}\bigr)$;
\item[$(2)$] $(0,f^-_{\rm ext}+f_{\rm hom})$
\end{enumerate}
are pointwise balanced.
\end{Proposition}
\begin{proof}
We first consider
\begin{align*}
\sigma_2(f_{\rm hom})&=\int_{\mathbb{P}\eta\in\mathbb{P}\operatorname{Im}\beta'_{\mathbf{c}'',\mathbf{c}'}\cap\operatorname{Ext}_{\mathcal{A}}^1(N,M)}q^{k(\mathbf{c}'',\mathbf{c}')+f_{\rm hom}((\mathbf{c}'',\mathbf{c}'), [\eta])}\\
&=\int_{\mathbb{P}\eta\in\mathbb{P}\operatorname{Im}\beta'_{\mathbf{c}'',\mathbf{c}'}\cap\operatorname{Ext}_{\mathcal{A}}^1(N,M)}q^{k(\mathbf{c}',\mathbf{c}'')}.
\end{align*}
Notice that $k(\mathbf{c}',\mathbf{c}'')$ and $k(\mathbf{c}'',\mathbf{c}')$ are independent of $[\epsilon]$ and $[\eta]$ as long as they are located in the domain of integration as follows.
Since
\[
\dim_k p_0\operatorname{Ker}\beta_{\mathbf{c}',\mathbf{c}''}+\dim_k \bigl(\operatorname{Im}\beta'_{\mathbf{c}'',\mathbf{c}'}\cap\operatorname{Ext}_{\mathcal{A}}^1(N,M)\bigr)=\dim_k \operatorname{Ext}_{\mathcal{A}}^1(M,N),
\]
we have two equalities
\begin{gather*}
\frac{q^{\dim_k \operatorname{Ext}_{\mathcal{A}}^1(M,N)}-1}{q-1}\\
\quad=q^{\dim_k (\operatorname{Im}\beta'_{\mathbf{c}'',\mathbf{c}'}\cap\operatorname{Ext}_{\mathcal{A}}^1(N,M))}\cdot\frac{q^{\dim_k p_0\operatorname{Ker}\beta_{\mathbf{c}',\mathbf{c}''}}-1}{q-1}+\frac{q^{\dim_k \operatorname{Im}\beta'_{\mathbf{c}'',\mathbf{c}'}\cap\operatorname{Ext}_{\mathcal{A}}^1(N,M)}-1}{q-1}\\
\quad=\frac{q^{\dim_k p_0\operatorname{Ker}\beta_{\mathbf{c}',\mathbf{c}''}}-1}{q-1}+q^{\dim_k p_0\operatorname{Ker}\beta_{\mathbf{c}',\mathbf{c}''}}\cdot\frac{q^{\dim_k \operatorname{Im}\beta'_{\mathbf{c}'',\mathbf{c}'}\cap\operatorname{Ext}_{\mathcal{A}}^1(N,M)}-1}{q-1}.
\end{gather*}
So by definition,
\begin{gather*}
\sigma_1\bigl(f^+_{\rm ext}\bigr)+\sigma_2({f_{\rm hom}})\\
\quad=\int_{\mathbb{P}\epsilon\in\mathbb{P}p_0\operatorname{Ker}\beta_{\mathbf{c}',\mathbf{c}''}}q^{k(\mathbf{c}',\mathbf{c}'')}\cdot q^{\dim_k (\operatorname{Im}\beta'_{\mathbf{c}'',\mathbf{c}'}\cap\operatorname{Ext}_{\mathcal{A}}^1(N,M))}
+\int_{\mathbb{P}\eta\in\mathbb{P}\operatorname{Im}\beta'_{\mathbf{c}'',\mathbf{c}'}\cap\operatorname{Ext}_{\mathcal{A}}^1(N,M)}q^{k(\mathbf{c}',\mathbf{c}'')}\\
\quad=q^{k(\mathbf{c}',\mathbf{c}'')}\cdot\bigg(\frac{q^{\dim_k p_0\operatorname{Ker}\beta_{\mathbf{c}',\mathbf{c}''}}-1}{q-1}\cdot q^{\dim_k (\operatorname{Im}\beta'_{\mathbf{c}'',\mathbf{c}'}\cap\operatorname{Ext}_{\mathcal{A}}^1(N,M))}\\
\phantom{\quad=}{}+\frac{q^{\dim_k \operatorname{Im}\beta'_{\mathbf{c}'',\mathbf{c}'}\cap\operatorname{Ext}_{\mathcal{A}}^1(N,M)}-1}{q-1}\bigg)\\
\quad=q^{k(\mathbf{c}',\mathbf{c}'')}\cdot \frac{q^{\dim_k \operatorname{Ext}_{\mathcal{A}}^1(M,N)}-1}{q-1}
\end{gather*}
and
\begin{gather*}
\sigma_1(0)+\sigma_2({f^-_{\rm ext}+f_{\rm hom}})\\
\quad=\int_{\mathbb{P}\epsilon\in\mathbb{P}p_0\operatorname{Ker}\beta_{\mathbf{c}',\mathbf{c}''}}q^{k(\mathbf{c}',\mathbf{c}'')}
+\int_{\mathbb{P}\eta\in\mathbb{P}\operatorname{Im}\beta'_{\mathbf{c}'',\mathbf{c}'}\cap\operatorname{Ext}_{\mathcal{A}}^1(N,M)}q^{k(\mathbf{c}',\mathbf{c}'')}\cdot q^{\dim_k p_0\operatorname{Ker}\beta_{\mathbf{c}',\mathbf{c}''}}\\
\quad=q^{k(\mathbf{c}',\mathbf{c}'')}\cdot\bigg(\frac{q^{\dim_k p_0\operatorname{Ker}\beta_{\mathbf{c}',\mathbf{c}''}}-1}{q-1}+\frac{q^{\dim_k \operatorname{Im}\beta'_{\mathbf{c}'',\mathbf{c}'}\cap\operatorname{Ext}_{\mathcal{A}}^1(N,M)}-1}{q-1}\cdot q^{\dim_k p_0\operatorname{Ker}\beta_{\mathbf{c}',\mathbf{c}''}}\bigg)\\
\quad=q^{k(\mathbf{c}',\mathbf{c}'')}\cdot \frac{q^{\dim_k \operatorname{Ext}_{\mathcal{A}}^1(M,N)}-1}{q-1}. \tag*{\qed}
\end{gather*} \renewcommand{\qed}{}
\end{proof}

\begin{Theorem}\label{mainlemma12}
If a pair of weight functions $\bigl(f^+,f^-\bigr)$ in $\mathbb{Z}_{\MF}$ is pointwise balanced, then for any weighted quantum cluster functions $f*_{[\epsilon']}\delta_M$ and $g*_{[\epsilon'']}\delta_N$,
we have
\begin{gather*}
 \big|\mathbb{P}\EExt_{\mathcal{A}}^1(M,N)\big|(f*_{[\epsilon']}\delta_M)*(g*_{[\epsilon'']}\delta_N)\\
\quad= \int_{\mathbb{P}\epsilon\in\mathbb{P}\EExt_{\mathcal{A}}^1(M,N)}\bigl(f^+*_{ [\epsilon]}\mathbb{S}_{fg}\bigr)*_{ [\epsilon]}\delta_{{\rm mt}\epsilon}+\int_{\mathbb{P}\eta\in\mathbb{P}\EExt_{\mathcal{A}}^1(N,M)}
(f^-*_{[\eta]}\mathbb{S}_{gf})*_{[\eta]}\delta_{{\rm mt}\eta}.
\end{gather*}
\end{Theorem}
\begin{proof}
We simplify the equality in the theorem to ${\rm l.h.s.}=\Sigma_1\bigl(f^+\bigr)+\Sigma_2(f^-)$.

Just for simplicity, we omit some independent variables of functions which are obvious in following calculation.
For example, $\mathbb{S}_{fg}(\mathbf{c}',\mathbf{c}'', \mathbf{0}_{MN})$ is simplified to
$\mathbb{S}_{fg}$.

Direct calculation:
\begin{gather*}
{\rm l.h.s.}(d_{\mathbf{i},\mathbf{a}})=\big|\mathbb{P}\EExt_{\mathcal{A}}^1(M,N)\big|(f*_{[\epsilon']}\delta_M)*(g*_{[\epsilon'']}\delta_N)(d_{\mathbf{i},\mathbf{a}})\\
\phantom{{\rm l.h.s.}(d_{\mathbf{i},\mathbf{a}})}{}=\big|\mathbb{P}\EExt_{\mathcal{A}}^1(M,N)\big|\iint_{\mathbf{a}'+\mathbf{a}''=\mathbf{a}, (\mathbf{c}',\mathbf{c}'')\in\mathcal{F}^{\rm mono}_{\mathbf{i},\mathbf{a}',M}\times\mathcal{F}^{\rm mono}_{\mathbf{i},\mathbf{a}'',N}}q^{k(\mathbf{c}',\mathbf{c}'')+\mathbb{S}_{fg}}\\
\phantom{{\rm l.h.s.}(d_{\mathbf{i},\mathbf{a}})}{}=\iint_{\mathbf{a}'+\mathbf{a}''=\mathbf{a}, (\mathbf{c}',\mathbf{c}'')\in\mathcal{F}^{\rm mono}_{\mathbf{i},\mathbf{a}',M}\times\mathcal{F}^{\rm mono}_{\mathbf{i},\mathbf{a}'',N}}q^{\mathbb{S}_{fg}}\cdot\frac{q^{\dim_k \EExt_{\mathcal{A}}^1(M,N)}-1}{q-1}\cdot q^{k(\mathbf{c}',\mathbf{c}'')},\\
\Sigma_1\bigl(f^+\bigr)(d_{\mathbf{i},\mathbf{a}})\\
\quad=\int_{\mathbb{P}\epsilon\in\mathbb{P}\EExt_{\mathcal{A}}^1(M,N)}\bigl(f^+*_{[\epsilon]} \mathbb{S}_{fg}\bigr)*_{[\epsilon]}\delta_{{\rm mt}\epsilon}(d_{\mathbf{i},\mathbf{a}})\\
\quad=\iiint_{\mathbb{P}\epsilon\in\mathbb{P}\EExt_{\mathcal{A}}^1(M,N), \mathbf{a}'+\mathbf{a}''=\mathbf{a}, (\mathbf{c}',\mathbf{c}'')\in\phi_{[\epsilon]}(\mathcal{F}^{\rm mono}_{\mathbf{i},\mathbf{a},L}([\epsilon],\mathbf{a}',\mathbf{a}''))}q^{k(\mathbf{c}',\mathbf{c}'')+ \mathbb{S}_{fg} +f^+(\mathbf{c}',\mathbf{c}'', [\epsilon])}\\
\quad=\iint_{\mathbf{a}'+\mathbf{a}''=\mathbf{a}, (\mathbf{c}',\mathbf{c}'')\in\mathcal{F}^{\rm mono}_{\mathbf{i},\mathbf{a}',M}\times\mathcal{F}^{\rm mono}_{\mathbf{i},\mathbf{a}'',N}}q^{\mathbb{S}_{fg}}\int_{\mathbb{P}\epsilon\in\mathbb{P}p_0\KKer\beta_{\mathbf{c}',\mathbf{c}''}}q^{k(\mathbf{c}',\mathbf{c}'')+f^+(\mathbf{c}',\mathbf{c}'', [\epsilon])}\\
\quad=\iint_{\mathbf{a}'+\mathbf{a}''=\mathbf{a}, (\mathbf{c}',\mathbf{c}'')\in\mathcal{F}^{\rm mono}_{\mathbf{i},\mathbf{a}',M}\times\mathcal{F}^{\rm mono}_{\mathbf{i},\mathbf{a}'',N}}q^{\mathbb{S}_{fg}}\cdot\sigma_1\bigl(f^+\bigr),\\
\Sigma_2(f^-)(d_{\mathbf{i},\mathbf{a}})\\
\quad=\int_{\mathbb{P}\eta\in\mathbb{P}\EExt_{\mathcal{A}}^1(N,M)}(f^-*_{[\eta]} \mathbb{S}_{gf})*_{[\eta]}\delta_{{\rm mt}\eta}(d_{\mathbf{i},\mathbf{a}})\\
\quad=\iiint_{\mathbb{P}\eta\in\mathbb{P}\EExt_{\mathcal{A}}^1(N,M), \mathbf{a}'+\mathbf{a}''=\mathbf{a}, (\mathbf{c}'',\mathbf{c}')\in\phi_{[\eta]}(\mathcal{F}^{\rm mono}_{\mathbf{i},\mathbf{a},L}([\eta],\mathbf{a}'',\mathbf{a}'))}q^{k(\mathbf{c}'',\mathbf{c}')+\mathbb{S}_{gf}+f^-(\mathbf{c}'', \mathbf{c}', [\eta])}\\
\quad=\iint_{\mathbf{a}'+\mathbf{a}''=\mathbf{a}, (\mathbf{c}',\mathbf{c}'')\in\mathcal{F}^{\rm mono}_{\mathbf{i},\mathbf{a}',M}\times\mathcal{F}^{\rm mono}_{\mathbf{i},\mathbf{a}'',N}}q^{
\mathbb{S}_{fg}}\int_{\mathbb{P}\eta\in\mathbb{P}\operatorname{Im}\beta'_{\mathbf{c}'',\mathbf{c}'}\cap\EExt_{\mathcal{A}}^1(N,M)}q^{k(\mathbf{c}'',\mathbf{c}')+f^-(\mathbf{c}'', \mathbf{c}', [\eta])}\\
\quad=\iint_{\mathbf{a}'+\mathbf{a}''=\mathbf{a}, (\mathbf{c}',\mathbf{c}'')\in\mathcal{F}^{\rm mono}_{\mathbf{i},\mathbf{a}',M}\times\mathcal{F}^{\rm mono}_{\mathbf{i},\mathbf{a}'',N}}q^{
\mathbb{S}_{fg}}\cdot\sigma_2(f^-).
\end{gather*}
By Definition~\ref{forproof-balance}, we have ${\rm l.h.s.}=\Sigma_1(f^+)+\Sigma_2(f^-)$.
\end{proof}

From Proposition \ref{mainlemma11} and Theorem \ref{mainlemma12}, we have
\begin{Theorem}[multiplication formula of weighted quantum cluster functions]\label{maintheorem1}
Let $\mathcal{A}$ be a~${\rm Hom}$-finite, $\Ext$-finite abelian category with $\EExt$-symmetry over a finite field $k=\mathbb{F}_q$ such that the iso-classes of objects form a set. For any weighted quantum cluster functions
$f*_{[\epsilon']}\delta_M$ and $g*_{[\epsilon'']}\delta_N$ such that $\operatorname{Ext}_{\mathcal{A}}^1(M, N)\neq 0$, we have
\begin{gather*}
\big|\mathbb{P}\EExt_{\mathcal{A}}^1(M,N)\big|(f*_{[\epsilon']}\delta_M)*(g*_{[\epsilon'']}\delta_N) \\
\quad=\int_{\mathbb{P}\epsilon\in\mathbb{P}\EExt_{\mathcal{A}}^1(M,N)}\bigl(f^+_{\rm ext}*_{[\epsilon]} \mathbb{S}_{fg}\bigr) *_{[\epsilon]} \delta_{{\rm mt}\epsilon}+\int_{\mathbb{P}\eta\in\mathbb{P}\EExt_{\mathcal{A}}^1(N,M)}(f_{\rm hom}*_{[\eta]}
\mathbb{S}_{gf})*_{[\eta]} \delta_{{\rm mt}\eta}\\
\quad=\int_{\mathbb{P}\epsilon\in\mathbb{P}\EExt_{\mathcal{A}}^1(M,N)} \mathbb{S}_{fg}*_{[\epsilon]}\delta_{{\rm mt}\epsilon}+\int_{\mathbb{P}\eta\in\mathbb{P}\EExt_{\mathcal{A}}^1(N,M)}(f^-_{\rm ext}*_{[\eta]} f_{\rm hom}*_{[\eta]} \mathbb{S}_{gf})*_{[\eta]}\delta_{{\rm mt}\eta}.
\end{gather*}
\end{Theorem}

Theorem~\ref{mainlemma12} shows that every pointwise balanced pair is a balanced pair as in the following definition.
\begin{Definition}
A pair of weight functions $\bigl(f^+,f^-\bigr)$ in $\mathbb{Z}_{\MF}$ is called a balanced pair if for any weighted quantum cluster functions $f*_{[\epsilon']}\delta_M$ and $g*_{[\epsilon'']}\delta_N$ such that $\EExt_{\mathcal{A}}^1(M,N)\neq 0$,
\begin{gather*}
\big|\mathbb{P}\EExt_{\mathcal{A}}^1(M,N)\big|(f*_{[\epsilon']}\delta_M)*(g*_{[\epsilon'']}\delta_N)\\
\quad= \int_{\mathbb{P}\epsilon\in\mathbb{P}\EExt_{\mathcal{A}}^1(M,N)}\bigl(f^+*_{[\epsilon]} \mathbb{S}_{fg}\bigr)*_{[\epsilon]}\delta_{{\rm mt}\epsilon}+\int_{\mathbb{P}\eta\in\mathbb{P}\EExt_{\mathcal{A}}^1(N,M)}
(f^-*_{[\eta]} \mathbb{S}_{gf})*_{[\eta]} \delta_{{\rm mt}\eta}
\end{gather*}
holds.
\end{Definition}
Thus by definition, any balanced pair satisfies a multiplication formula as in Theorem~\ref{maintheorem1}.

\subsection[The case dim\_k EExt\_A\^1(M,N)=1]{The case $\boldsymbol{\dim_k\EExt_{\mathcal{A}}^1(M,N)=1}$}
Finally, we simplify the formula in the case that $\dim_k\EExt_{\mathcal{A}}^1(M,N)=\dim_k\EExt_{\mathcal{A}}^1(N,M)=1$, so that $\mathbb{P}\EExt_{\mathcal{A}}^1(M,N)$ and $\mathbb{P}\EExt_{\mathcal{A}}^1(N,M)$ are both singleton sets. Moreover, since $p_0\KKer\beta_{\mathbf{c}',\mathbf{c}''}$ and $\operatorname{Im}\beta'_{\mathbf{c}'',\mathbf{c}'}\cap\EExt_{\mathcal{A}}^1(N,M)$ are orthogonal, one of them is of dimension $1$ and the other is zero.

\begin{Lemma}\label{onedim1}
Assume $\dim_k \EExt_{\mathcal{A}}^1(M,N)=1$, and that $[\epsilon]$ and $[\eta]$ are non-zero in $\EExt_{\mathcal{A}}^1(M,N)$ and $\EExt_{\mathcal{A}}^1(N,M)$, respectively. We have
\begin{enumerate}\itemsep=0pt
\item[$(1)$] If $(\mathbf{c}',\mathbf{c}'')\in \phi_{[\epsilon]}(\mathcal{F}^{\rm mono}_{\mathbf{i},\mathbf{a},L}([\epsilon],\mathbf{a}',\mathbf{a}''))$, then $f^+_{\rm ext}(\mathbf{c}',\mathbf{c}'', [\epsilon])=0$;
\item[$(2)$] If $(\mathbf{c}'',\mathbf{c}')\in \phi_{[\eta]}(\mathcal{F}^{\rm mono}_{\mathbf{i},\mathbf{a},L}([\eta],\mathbf{a}'',\mathbf{a}'))$, then $f^-_{\rm ext}(\mathbf{c}'',\mathbf{c}',[\eta])=0$.
\end{enumerate}
\end{Lemma}
\begin{proof}
By Corollary~\ref{image-related}, $(\mathbf{c}',\mathbf{c}'')\in \phi_{[\epsilon]}(\mathcal{F}^{\rm mono}_{\mathbf{i},\mathbf{a},L}([\epsilon],\mathbf{a}',\mathbf{a}''))$ if and only if $[\epsilon]\in p_0\KKer\beta_{\mathbf{c}',\mathbf{c}''}.$
In this case, $p_0\KKer\beta_{\mathbf{c}',\mathbf{c}''}$ contains a non-zero element, so $\dim_k p_0\Ker\beta_{\mathbf{c}',\mathbf{c}''}=1$.
By definition, $f^+_{\rm ext}((\mathbf{c}',\mathbf{c}''), [\epsilon])=\dim_k \bigl(\operatorname{Im}\beta'_{\mathbf{c}'',\mathbf{c}'}\cap\EExt_{\mathcal{A}}^1(N,M)\bigr)=0$.

Dually, by Lemma~\ref{dual-1},
\[(\mathbf{c}'',\mathbf{c}')\in \phi_{[\eta]}(\mathcal{F}^{\rm mono}_{\mathbf{i},\mathbf{a},L}([\eta],\mathbf{a}'',\mathbf{a}'))\quad \text{if and only if}\quad [\eta]\in \operatorname{Im}\beta'_{\mathbf{c}'',\mathbf{c}'}\cap\EExt_{\mathcal{A}}^1(N,M).\]
In this case, $\dim_k p_0\KKer\beta_{\mathbf{c}',\mathbf{c}''}=0$.
\end{proof}

\begin{Corollary}
Assume that $\dim_k \EExt_{\mathcal{A}}^1(M,N)=1$, and that $[\epsilon]$ and $[\eta]$ are non-zero in $\EExt_{\mathcal{A}}^1(M, N)$ and $\EExt_{\mathcal{A}}^1(N,M)$, respectively. Then
\[(\mathbf{c}',\mathbf{c}'')\in \phi_{[\epsilon]}\bigl(\mathcal{F}^{\rm mono}_{\mathbf{i},\mathbf{a},L}([\epsilon],\mathbf{a}',\mathbf{a}'')\bigr) \qquad\text{if and only if} \quad (\mathbf{c}'',\mathbf{c}')\notin \phi_{[\eta]}\bigl(\mathcal{F}^{\rm mono}_{\mathbf{i},\mathbf{a},L}([\eta],\mathbf{a}'',\mathbf{a}')\bigr).\]
\end{Corollary}
\begin{proof}
Since $\dim_k \EExt_{\mathcal{A}}^1(M,N)=1$, $\dim_k p_0\KKer\beta_{\mathbf{c}',\mathbf{c}''}$ is either $0$ or $1$. Both conditions are equivalent to $\dim_k p_0\KKer\beta_{\mathbf{c}',\mathbf{c}''}=1$.
\end{proof}

\begin{Theorem}
Assume $\dim_k \EExt_{\mathcal{A}}^1(M,N)=1$ and there are nonsplit short exact sequences
\[
\epsilon\colon\ 0\rightarrow N\rightarrow L\rightarrow M\rightarrow0, \quad\text{ and } \quad \eta\colon\ 0\rightarrow M\rightarrow L'\rightarrow N\rightarrow0.
\]

Then we have $(f*_{[\epsilon']}\delta_M)*(g*_{[\epsilon'']}\delta_N)=\mathbb{S}_{fg}*{[\epsilon]}\delta_L+(f_{\rm hom}*_{[\eta]} \mathbb{S}_{gf})*_{[\eta]}\delta_{L'}$. In particular, $\delta_M*\delta_N=\delta_L+(f_{\rm hom}*_{[\eta]} f)*_{[\eta]}\delta_{L'}$, where $f$ is the zero function in $\mathbb{Z}_{\MF}$.
\end{Theorem}
\begin{proof}
Recall the multiplication formula
\begin{gather*}
\big|\mathbb{P}\EExt_{\mathcal{A}}^1(M,N)\big|(f*_{[\epsilon']}\delta_M)*(g*_{[\epsilon'']}\delta_N) \\
\quad=\int_{\mathbb{P}\epsilon\in\mathbb{P}\EExt_{\mathcal{A}}^1(M,N)}\bigl(f^+_{\rm ext}*_{[\epsilon]} \mathbb{S}_{fg}\bigr) *_{[\epsilon]} \delta_{{\rm mt}\epsilon}+\int_{\mathbb{P}\eta\in\mathbb{P}\EExt_{\mathcal{A}}^1(N,M)}(f_{\rm hom}*_{[\eta]}\mathbb{S}_{gf})*_{[\eta]}\delta_{{\rm mt}\eta}\\
\quad=\int_{\mathbb{P}\epsilon\in\mathbb{P}\EExt_{\mathcal{A}}^1(M,N)} \mathbb{S}_{fg}*_{[\epsilon]} \delta_{{\rm mt}\epsilon}+\int_{\mathbb{P}\eta\in\mathbb{P}\EExt_{\mathcal{A}}^1(N,M)}(f^-_{\rm ext}*_{[\eta]}f_{\rm hom}*_{[\eta]} \mathbb{S}_{gf})*_{[\eta]}\delta_{{\rm mt}\eta}
\end{gather*}
from Theorem~\ref{maintheorem1}.

Since $\dim_k\EExt_{\mathcal{A}}^1(M,N)=1$, $\big|\mathbb{P}\EExt_{\mathcal{A}}^1(M,N)\big|=1$ and both integrals in the formula degenerate into a pair of terms indexed by $[\epsilon]$ and $[\eta]$ respectively. Then we have
\begin{align*}
(f*_{[\epsilon']}\delta_M)*(g*_{[\epsilon'']}\delta_N) ={}&\bigl(f^+_{\rm ext}*_{[\epsilon]} \mathbb{S}_{fg}\bigr)*_{\epsilon]} \delta_L+(f_{\rm hom}*_{[\eta]}\mathbb{S}_{gf})*_{[\eta]}\delta_{L'}\\
={}& \mathbb{S}_{fg}*_{[\epsilon]}\delta_L+(f^-_{\rm ext}*_{[\eta]}f_{\rm hom}*_{[\eta]} \mathbb{S}_{gf})*_{[\eta]}\delta_{L'}.
\end{align*}

From Lemma \ref{onedim1},
\begin{align*}
\bigl(f^+_{\rm ext}*_{[\epsilon]} \mathbb{S}_{fg}\bigr)*_{[\epsilon]}\delta_L(d_{\mathbf{i},\mathbf{a}})
&=\iint_{\mathbf{a}'+\mathbf{a}''=\mathbf{a}, (\mathbf{c}',\mathbf{c}'')\in\phi_{[\epsilon]}(\mathcal{F}^{\rm mono}_{\mathbf{i},\mathbf{a},L}([\epsilon],\mathbf{a}',\mathbf{a}''))}q^{k(\mathbf{c}',\mathbf{c}'')+\mathbb{S}_{fg}
+f^+_{\rm ext}(\mathbf{c}',\mathbf{c}'', [\epsilon])}\\
&=\iint_{\mathbf{a}'+\mathbf{a}''=\mathbf{a}, (\mathbf{c}',\mathbf{c}'')\in\phi_{[\epsilon]}(\mathcal{F}^{\rm mono}_{\mathbf{i},\mathbf{a},L}([\epsilon],\mathbf{a}',\mathbf{a}''))}q^{k(\mathbf{c}',\mathbf{c}'')+\mathbb{S}_{fg}}\\
&=\mathbb{S}_{fg}*_{[\epsilon]}\delta_L(d_{\mathbf{i},\mathbf{a}})
\end{align*}
and
\begin{gather*}
(f^-_{\rm ext}*_{[\eta]}f_{\rm hom}*_{[\eta]} \mathbb{S}_{gf})*_{[\eta]}\delta_{L'}(d_{\mathbf{i},\mathbf{a}})\\
\quad=\iint_{\mathbf{a}'+\mathbf{a}''=\mathbf{a}, (\mathbf{c}'',\mathbf{c}'')\in\phi_{[\eta]}(\mathcal{F}^{\rm mono}_{\mathbf{i},\mathbf{a},L}([\eta],\mathbf{a}'',\mathbf{a}'))}q^{k(\mathbf{c}'',\mathbf{c}')+
\mathbb{S}_{gf}+f^-_{\rm ext}(\mathbf{c}'', \mathbf{c}',[\eta])+f_{\rm hom}(\mathbf{c}'', \mathbf{c}',[\eta])}\\
\quad=\iint_{\mathbf{a}'+\mathbf{a}''=\mathbf{a}, (\mathbf{c}'',\mathbf{c}')\in\phi_{[\eta]}(\mathcal{F}^{\rm mono}_{\mathbf{i},\mathbf{a},L}([\eta],\mathbf{a}'',\mathbf{a}'))}q^{k(\mathbf{c}'',\mathbf{c}')+\mathbb{S}_{gf}+f_{\rm hom}
(\mathbf{c}'', \mathbf{c}',[\eta])}\\
\quad=(f_{\rm hom}*_{[\eta]} \mathbb{S}_{gf})*_{[\eta]}\delta_{L'}(d_{\mathbf{i},\mathbf{a}}). \tag*{\qed}
\end{gather*} \renewcommand{\qed}{}
\end{proof}

\section[2-Calabi--Yau triangulated categories and multiplication formula]{2-Calabi--Yau triangulated categories \\ and multiplication formula}\label{sec4}
\subsection{2-Calabi--Yau triangulated categories}\label{index}
Let $\mathcal{C}$ be a $\Hom$-finite, 2-Calabi--Yau, Krull--Schmidt triangulated category over a finite field $k=\mathbb{F}_q$, which admits a cluster tilting object $T$. The shift functor on $\mathcal{C}$ is denoted by $\Sigma$. Let $B$ be the endomorphism algebra of $T$, so there is a functor
\begin{align*}
F:=\Hom(T,-)\colon\ \mathcal{C}\longrightarrow&\mod B,\\
X\longmapsto &\Hom(T,X)
\end{align*}
from $\mathcal{C}$ to the abelian category $\mod B$. This induces an equivalence of categories
 \[\mathcal{C}/(\Sigma T) \stackrel{\simeq}{\longrightarrow} \operatorname{mod}B,\]
where $(\Sigma T)$ denotes the ideal of morphisms of $\mathcal{C}$ which factor through a direct sum of copies of $\Sigma T$
and $\mathcal{C}/(\Sigma T)$ is the corresponding quotient category (see~\cite{Palu2008}).
Let $\{S_1,\dots,S_n\}$ be a~complete set of isomorphism classes of simple objects in $\mod B$.

For any two objects $M$, $N$ in $\mod B$, define
\begin{gather*} \langle M,N\rangle=\dim_k\Hom_B(M,N)-\dim_k\Ext_B^1(M,N), \qquad
\langle M,N\rangle_a=\langle M,N\rangle-\langle N,M\rangle.
\end{gather*}
The form $\langle -,-\rangle_a$ can be reduced to the Grothendieck group $K_0(\operatorname{mod}B)$ (see~\cite[Lemma 1.3]{Palu2008}).

For any object $X$ in $\mathcal{C}$, there are two triangles
\[T_2\rightarrow T_1\rightarrow X\rightarrow \Sigma T_2 \qquad \text{and}\qquad \Sigma T_4\rightarrow X\rightarrow \Sigma^2 T_3\rightarrow \Sigma^2 T_4\]
with $T_1$, $T_2$, $T_3$, $T_4$ in ${\rm add}T$. So we can define index and coindex of $X$
\[
\operatorname{ind}X:=[FT_1]-[FT_2] \qquad \text{and}\qquad \operatorname{coind}X:=[FT_3]-[FT_4],\]
where $[-]$ represents the equivalence class of an object in the Grothendieck group $K_0(\operatorname{proj}B)$.

\subsection{Exact structure}
An element $\epsilon\in\Hom_{\mathcal{C}}(M,\Sigma N)$ induces a triangle
\[
N\xrightarrow{i} L\xrightarrow{p} M\xrightarrow{\epsilon} \Sigma N
\]
in $\mathcal{C}$. In this case, we denote ${\rm qt}\epsilon=M$, ${\rm st}\epsilon=N$ and ${\rm mt}\epsilon=L$.
Applying $F$, we get an exact sequence in $\mod B$
\[
FN\xrightarrow{Fi} FL\xrightarrow{Fp} FM\xrightarrow{F\epsilon} F\Sigma N.
\]
Assume $M_0\subseteq FM$ and $N_0\subseteq FN$ are submodules, and consider the diagram
\[
\begin{tikzcd}
FN \arrow[r] & FL \arrow[r] & FM \arrow[r] & F\Sigma N.\\
N_0 \arrow[u, hook] \arrow[r, dashed] & L_0 \arrow[r, dashed] \arrow[u, dashed, hook] & M_0 \arrow[u, hook] &
\end{tikzcd}
\]
We denote \[{\rm Gr}^{\epsilon}_{M_0,N_0}(FL)=\big\{L_0\subseteq FL|Fi^{-1}(L_0)=N_0,Fp(L_0)=M_0\big\}.\]

We now calculate $|{\rm Gr}^{\epsilon}_{M_0,N_0}(FL)|$. Recall that $\mathcal{C}$ is 2-Calabi--Yau, so for any two objects $M$,~$N$ in $\mathcal{C}$, there are natural isomorphisms
\[\Hom_{\mathcal{C}}(M,\Sigma N)\rightarrow D\Hom_{\mathcal{C}}(N,\Sigma M)\qquad \text{and}\qquad \Hom_{\mathcal{C}}\bigl(\Sigma^{-1}M,N\bigr)\rightarrow \Hom_{\mathcal{C}}(M,\Sigma N),\]
 which induce an isomorphism
\[D_{M,N}\colon\ \Hom_{\mathcal{C}}\bigl(\Sigma^{-1}M,N\bigr)\longrightarrow D\Hom_{\mathcal{C}}(N,\Sigma M)\]
and then a bilinear form
\[\widetilde{D}_{M,N}\colon\ \Hom_{\mathcal{C}}(\Sigma^{-1}M,N)\times\Hom_{\mathcal{C}}(N,\Sigma M)\longrightarrow k,\,\,(a,b)\longmapsto D_{M,N}(a)(b).\]

Via the composition of functors $\mathcal{C} \longrightarrow \mathcal{C}/(\Sigma T) \stackrel{\simeq}{\longrightarrow} \operatorname{mod} B$, any two inclusions $M_0\subseteq FM$ and~$N_0\subseteq FN$
in $\operatorname{mod} B$ can be lifted to two morphisms $\tilde{M}_0\xrightarrow{\iota_M}M$ and $\tilde{N}_0\xrightarrow{\iota_N}N$ in $\mathcal{C}$, respectively.
Then we can define two linear maps
\begin{align*}
\alpha_{M_0,N_0}\colon&\Hom_{\mathcal{C}}\bigl(\Sigma^{-1}M,\tilde{N}_0\bigr)\oplus\Hom_{\mathcal{C}}\bigl(\Sigma^{-1}M,N\bigr)\\
&\quad\longrightarrow\Hom_{\mathcal{C}/(T)}\bigl(\Sigma^{-1}\tilde{M}_0,\tilde{N}_0\bigr)\oplus\Hom_{\mathcal{C}}\bigl(\Sigma^{-1}\tilde{M}_0,N\bigr)\oplus\Hom_{\mathcal{C}/(\Sigma T)}\bigl(\Sigma^{-1}M,N\bigr),\\
&(a,b)\longmapsto \bigl(a\circ\Sigma^{-1}\iota_M,\iota_N\circ a\circ\Sigma^{-1}\iota_M-b\circ\Sigma^{-1}\iota_M,\iota_N\circ a-b\bigr)
\end{align*}
and
\begin{align*}
\alpha'_{N_0,M_0}\colon&\Hom_{\Sigma T}\bigl(\tilde{N}_0,\Sigma\tilde{M}_0\bigr)\oplus\Hom_{\mathcal{C}}\bigl(N,\Sigma \tilde{M}_0\bigr)\oplus\Hom_{\Sigma^2T}(N,\Sigma M)\\
&\quad\longrightarrow\Hom_{\mathcal{C}}\bigl(\tilde{N}_0,\Sigma M\bigr)\oplus\Hom_{\mathcal{C}}(N,\Sigma M),\\
&(a,b,c)\longmapsto (\Sigma\iota_M\circ a+c\circ\iota_N+\Sigma\iota_M\circ b\circ\iota_N,-c-\Sigma\iota_M\circ b),
\end{align*}
where $\Hom_{\Sigma T} (-,-)$ and $\Hom_{\Sigma^2 T} (-,-)$ represent the homomorphisms of $\mathcal{C}$ which factor through a direct sum of copies of $\Sigma T$ and $\Sigma^2 T$, respectively.

The following properties of $\alpha_{M_0,N_0}$ and $\alpha'_{N_0,M_0}$ are given in \cite[Lemma 4.2, Proposition 4.3]{Palu2008}.
\begin{Lemma}
For any $M$, $N$ in $\mathcal{C}$, let $\epsilon\in\Hom_{\mathcal{C}}(M,\Sigma N)$ and $\eta\in\Hom_{\mathcal{C}}(N,\Sigma M)$ and choose submodules $M_0\subseteq FM$, $N_0\subseteq FN$.
\begin{enumerate}\itemsep=0pt
\item[$(1)$] ${\rm Gr}^{\epsilon}_{M_0,N_0}(F({\rm mt}\epsilon))$ is non-empty if and only if $\Sigma^{-1}\epsilon\in p\Ker\alpha_{M_0,N_0}$ where $p$ is the projection from \[\Hom_{\mathcal{C}}\bigl(\Sigma^{-1}M,\tilde{N}_0\bigr)\oplus\Hom_{\mathcal{C}}\bigl(\Sigma^{-1}M,N\bigr)\qquad \text{to} \quad\Hom_{\mathcal{C}}\bigl(\Sigma^{-1}M,N\bigr);\]
\item[$(2)$] ${\rm Gr}^{\eta}_{N_0,M_0}(F({\rm mt}\eta))$ is non-empty if and only if $\eta\in\operatorname{Im}\alpha'_{N_0,M_0}\cap\Hom_{\mathcal{C}}(N,\Sigma M)$ where the intersection is realized through regarding $\Hom_{\mathcal{C}}(N,\Sigma M)$ as a linear subspace of \[\Hom_{\mathcal{C}}\bigl(\tilde{N}_0,\Sigma M\bigr)\oplus\Hom_{\mathcal{C}}(N,\Sigma M).
 \]
\end{enumerate}
\end{Lemma}

Moreover, note that $\alpha_{M_0,N_0}$ is the dual of $\alpha'_{N_0,M_0}$ with respect to the pairing $\widetilde{D}_{M,N}$ \cite[Lemma 7.3.1]{GLS2007}, so
\begin{Lemma}\label{orthogonal2}
For any $M$, $N$ in $\mathcal{C}$ and submodules $M_0\subseteq FM$, $N_0\subseteq FN$, we have
\[
\Ker\alpha_{M_0,N_0}=(\operatorname{Im}\alpha'_{N_0,M_0})^\perp.
\]
In particular,
\[
\dim_k p\Ker\alpha_{M_0,N_0}+\dim_k (\operatorname{Im}\alpha'_{N_0,M_0}\cap\Hom_{\mathcal{C}}(N,\Sigma M))=\dim_k\Hom_{\mathcal{C}}(M,\Sigma N).
\]
\end{Lemma}

Then, to calculate $|{\rm Gr}^{\epsilon}_{M_0,N_0}(FL)|$ when it is non-zero, we use the following result.
\begin{Lemma}\label{same-dim}
If ${\rm Gr}^{\epsilon}_{M_0,N_0}(F({\rm mt}\epsilon))$ is non-empty, then $\Hom_B(M_0,FN/N_0)$ acts
freely and transitively on ${\rm Gr}^{\epsilon}_{M_0,N_0}(F({\rm mt}\epsilon)).$
Moreover, if we set
\[
l(M,N,M_0,N_0):=\dim_k\Hom_B(M_0,FN/N_0),
\]
then $|{\rm Gr}^{\epsilon}_{M_0,N_0}(FL)|=q^{l(M,N,M_0,N_0)}$.
\end{Lemma}
\begin{proof}
We sketch the proof.

Consider the commutative diagram in $\mod B$
\[
\begin{tikzcd}
FN/N_0 & & & \\
FN \arrow[r, "Fi"] \arrow[u, "\pi", two heads] & FL \arrow[r, "Fp"] & FM \arrow[r, "F\epsilon"] & F\Sigma N \\
N_0 \arrow[u, "i_N", hook] \arrow[r, "Fi"] & L_0 \arrow[u, "i_L", hook] \arrow[r, "Fp"] &\, M_0 \arrow[u, "i_M", hook], &
\end{tikzcd}
\]
where ${\rm mt}\epsilon=L$ and $M_0$, $N_0$, $L_0$ are submodules with $Fi^{-1}(L_0)=N_0$, $Fp(L_0)=M_0$.

We can define an action of $\Hom_B(M_0,FN/N_0)$ on ${\rm Gr}^{\epsilon}_{M_0,N_0}(FL)$ as follows. For any
\begin{gather*}
f\in\Hom_B(M_0,FN/N_0) \qquad \text{and}\qquad L_0\in{\rm Gr}^{\epsilon}_{M_0,N_0}(FL),\\
L_0^f:=\{Fi(n)+x\mid n\in FN,\, x\in L_0,\, f(Fp(x))=\pi(n)\}
\end{gather*}
is a linear subspace of $FL$. Since $Fi$, $Fp$, $f$, $\pi$ are $B$-module homomorphisms,
\[
B(Fi(n)+x)=B(Fi(n))+B(x)=Fi(B(n))+B(x)
\]
and
\[
f(Fp(B(x)))=B(f(Fp(x)))=B(\pi(n))=\pi(B(n)).
\]
So $L_0^f$ is a submodule of $FL$. For $Fi(n)+x \in L_0^f$, since $Fp(Fi(n)+x)=Fp(x) \in M_0$ and $Fi^{-1}(Fi(n))=Fi^{-1}(n_0)$ for some $n_0 \in N_0$,
we have $Fi^{-1}(Fi(n)+x)=Fi^{-1}(n_0) +Fi^{-1}(x) \in N_0$ and then \smash{$Fi^{-1}\bigl(L^f_0\bigr)=N_0$}, \smash{$Fp\bigl(L^f_0\bigr)=M_0$}. Thus \smash{$L^f_0\in{\rm Gr}^{\epsilon}_{M_0,N_0}(FL)$}.

If $f=0$, then $f(Fp(x))=0$ and $\pi(n)=0$ implies $n\in N_0$, so $L^0_0=L_0$. Moreover, \smash{$\bigl(L_0^f\bigr)^g=L_0^{f+g}$}. Since $Fi^{-1}(L_0)=N_0$, $Fi(n)\in L_0$ if and only if $n\in N_0$ and in this case, $\pi(n)=0$. So if \smash{$L^f_0=L_0$}, then $f(Fp(L_0))=f(M_0)=0$. That is to say, the action is free. For any $L'_0\in{\rm Gr}^{\epsilon}_{M_0,N_0}(FL)$, define $f$ as $f(m)=\pi(x'-x)$ where $x\in L_0$ and $x'\in L'_0$ with $Fp(x)=Fp(x')=m$. Then \smash{$L^f_0=L'_0$}, so the action is transitive.
\end{proof}

Moreover, according to the proof, we can define a linear structure on ${\rm Gr}^{\epsilon}_{M_0,N_0}(FL)$ with respect to a fixed $L_0$ in ${\rm Gr}^{\epsilon}_{M_0,N_0}(FL)$ by
\[
\lambda L_0^f =L_0^{\lambda f},\qquad L_0^f+L_0^g =L_0^{f+g}.
\]
Then the map $\Hom_B(M_0,FN/N_0) \rightarrow {\rm Gr}^{\epsilon}_{M_0,N_0}(F({\rm mt}\epsilon))$ which sends $f$ to $L_0^f$ is a linear isomorphism.

\begin{Remark}
Notice that $\Hom_B(M_0,FN/N_0)$ just depends on $M$, $N$, $M_0$, $N_0$ and is independent of ${\rm mt}\epsilon$. So as long as $\Sigma^{-1}\epsilon$ is chosen from $p\KKer\alpha_{M_0,N_0}$, the dimension $l(M,N,M_0,N_0)$ is invariant. That is to say, all $\epsilon$ in $\Sigma p\KKer\alpha_{M_0,N_0}$ lead to the same quantity $|{\rm Gr}^{\epsilon}_{M_0,N_0}(FL)|$.
\end{Remark}

\begin{Corollary}\label{commondimension}
All submodules in ${\rm Gr}^{\epsilon}_{M_0,N_0}(F({\rm mt}\epsilon))$ have the same dimension vector.
\end{Corollary}
\begin{proof}


By Lemma~\ref{same-dim}, any submodule in ${\rm Gr}^{\epsilon}_{M_0,N_0}(F({\rm mt}\epsilon))$ has the form $L_0^f$ for some $f\in \Hom_B(M_0,FN/N_0)$.
Since $L^f_0\in{\rm Gr}^{\epsilon}_{M_0,N_0}(FL)$, we have the commutative diagram
\[
\begin{tikzcd}
FN \arrow[r, "Fi"] & FL \arrow[r, "Fp"] & FM \\
N_0 \arrow[u, "i_N", hook] \arrow[r, "Fi"] & L^f_0 \arrow[u, "i_L", hook] \arrow[r, "Fp"] & M_0 \arrow[u, "i_M", hook].
\end{tikzcd}
\]
Both rows are exact at the middle term, therefore $L_0^f/Fi(N_0)\cong M_0$ and then
\[
\operatorname{dim} L_0^f=\operatorname{dim} Fi(N_0)+ \operatorname{dim} M_0.\tag*{\qed}
\]\renewcommand{\qed}{}
\end{proof}

Summarizing all the analysis in this subsection, we have
\begin{Proposition}\label{vectorbundle20}
For any $M$, $N$ in $\mathcal{C}$, let $\epsilon \in\Hom_{\mathcal{C}}(M,\Sigma N)$ and choose submodules $M_0\subseteq FM$, $N_0\subseteq FN$.
\begin{enumerate}\itemsep=0pt
\item[$(1)$] The set ${\rm Gr}^{\epsilon}_{M_0,N_0}(F({\rm mt}\epsilon))$ is non-empty if and only if $\epsilon\in \Sigma p\KKer\alpha_{M_0,N_0}$ if and only if
$\epsilon \in \operatorname{Im}\alpha'_{M_0,N_0}\cap\Hom_{\mathcal{C}}(M,\Sigma N)$;
\item[$(2)$] If ${\rm Gr}^{\epsilon}_{M_0,N_0}(F({\rm mt}\epsilon))$ is non-empty, then there is a bijection \[\Hom_{B}(M_0,FN/N_0) \longrightarrow {\rm Gr}^{\epsilon}_{M_0,N_0}(F({\rm mt}\epsilon))\] and
$|{\rm Gr}^{\epsilon}_{M_0,N_0}(F({\rm mt}\epsilon))|=q^{l(M,N,M_0,N_0)}$.
\end{enumerate}
\end{Proposition}

\subsection{Quantum cluster functions}\label{sectionofskewpolynomial}
Let $A_{n,\lambda}$ be the $\mathbb{Q}$-algebra generated by $x_1^\pm,\dots,x_n^\pm$ with defining relations: for any $\e=(e_1,\dots,\allowbreak e_n)$, $\f=(f_1,\dots,f_n)\in\mathbb{Z}^n$,
\[X^\e\cdot X^\f=q^{\frac{1}{2}\lambda(\e,\f)}X^{\e+\f},\]
where $X^\e=x_1^{e_1}\dots x_n^{e_n}$ and $X^\f=x_1^{f_1}\dots x_n^{f_n}$ are monomials in $A_{n,\lambda}$ and $\lambda(-,-)$ is a skew-symmetric bilinear form. Such $A_{n,\lambda}$ is called a skew-polynomial algebra and elements in $A_{n,\lambda}$ are called skew-polynomials.

On the other hand, for any $L\in\mathcal{C}$ and $\g\in\mathbb{N}^n$, the set of all submodules $F_0$ of $FL$ with dimension vector $\g$ is denoted by ${\rm Gr}_\g(FL)$, called a quiver Grassmannian.

\begin{Definition}\label{pld}
For each object $L$ in $\mathcal{C}$, we assign a skew-polynomial in $A_{n,\lambda}$ as
\[X_L:=\int_{\g}|{\rm Gr}_\g(FL)|\cdot X^{p(L,\g)},\]

where the integral runs over all $\g$ in $\mathbb{N}^n$ and $p(L,\g)\in\mathbb{N}^n$ with \[p(L,\g)_i=-({\rm coind}L)_i+\langle \underline{\dim}_k S_i,\g\rangle_a.\] The skew-polynomial is called a quantum cluster function of $L$.
\end{Definition}

Notice that the set of all $\g$ such that ${\rm Gr}_\g(FL)$ is non-empty is finite, so the integral is just a finite sum.

\subsection{Mappings with affine fibers}\label{sec4.4}
Given $M,N\in\mathcal{C}$, define
\[\operatorname{EG}(M,N):=\biggl\{(\epsilon,L_0)\mid \epsilon\in\Hom_{\mathcal{C}}(M,\Sigma N),\, L_0\in\coprod_{\g}{\rm Gr}_\g(F({\rm mt}\epsilon))\biggr\}.\]
Then from $(\epsilon,L_0)\in \operatorname{EG}(M,N)$, we can induce two submodules $M_0$, $N_0$ of $FM$ and $FN$ respectively as $N_0=i^{-1}(L_0)$ and $M_0=p(L_0)$. Let $L:={\rm mt}\epsilon$.
Then we naturally have the commutative diagram
\[
\begin{tikzcd}
FN \arrow[r, "i"] & FL \arrow[r, "p"] & FM \arrow[r, "F\epsilon"] & F\Sigma N.\\
N_0 \arrow[r] \arrow[u, hook] & L_0 \arrow[r] \arrow[u, hook] & M_0 \arrow[u, hook] &
\end{tikzcd}
\]
We denote this assignment by
\[\psi_{MN}\colon\ \operatorname{EG}(M,N)\longrightarrow\coprod_{\e,\f}{\rm Gr}_\e(FM)\times{\rm Gr}_\f(FN),\]
where the coproduct runs over $\mathbb{N}^n\times\mathbb{N}^n$.

Fix $\epsilon$, and set
\[\psi_{\epsilon}:=\psi_{MN}(\epsilon,-)\colon\coprod_{\g}{\rm Gr}_\g(FL)\longrightarrow\coprod_{\e,\f}{\rm Gr}_\e(FM)\times{\rm Gr}_\f(FN).\]
To separate all pieces indexed by $\g$, we set
\[\psi_{\epsilon,\g}:=\psi_{\epsilon}|_{{\rm Gr}_\g(FL)}\colon{\rm Gr}_\g(FL)\longrightarrow\coprod_{\e,\f}{\rm Gr}_\e(FM)\times{\rm Gr}_\f(FN).\]
Moreover, denote
\[{\rm Gr}^\epsilon_\g(FL,\e,\f):=\psi_{\epsilon,\g}^{-1}({\rm Gr}_\e(FM)\times{\rm Gr}_\f(FN))\]
and
\[{\rm Gr}_{\e,\f}^\epsilon(FM,FN,\g):=\psi_{\epsilon,\g}({\rm Gr}^\epsilon_\g(FL,\e,\f)).\]

Notice that ${\rm Gr}_{\e,\f}^\epsilon(FM,FN,\g)$ may not be the whole of ${\rm Gr}_\e(FM)\times{\rm Gr}_\f(FN)$ because $\psi_{\epsilon,\g}$ is not surjective in general.

Rewriting Proposition \ref{vectorbundle20} in this notation, we get the following.
\begin{Lemma}\label{vectorbundle3}
Given $(M_0,N_0)\in{\rm Gr}_\e(FM)\times{\rm Gr}_\f(FN)$ and $\epsilon\in\Hom_{\mathcal{C}}(M,\Sigma N)$,
If $(M_0,N_0)\in{\rm Gr}_{\e,\f}^\epsilon(FM,FN,\g)$, then $\Hom_{B}(M_0,FN/N_0)$ acts
freely and transitively on $\psi_{\epsilon,\g}^{-1}(M_0,N_0)$.
\end{Lemma}




\subsection{Cardinality}
In this subsection, we calculate cardinalities of sets involved in several lemmas above and refine the calculation of $X_L$.

Considering $\psi_{MN}$, since ${\rm Gr}^\epsilon_\g(FL,\e,\f)=\psi_{\epsilon,\g}^{-1}({\rm Gr}_\e(FM)\times{\rm Gr}_\f(FN))$ and the Grassmannian is a finite set here, we have
\[|{\rm Gr}_\g(FL)|=\int_{\e,\f}|{\rm Gr}^\epsilon_\g(FL,\e,\f)|\]
for any $\epsilon\in\Hom_{\mathcal{C}}(M,\Sigma N)$ with ${\rm mt}\epsilon=L$, and
\[|{\rm Gr}^\epsilon_\g(FL,\e,\f)|=\int_{(M_0,N_0)\in{\rm Gr}_\e(FM)\times{\rm Gr}_\f(FN)}\big|\psi_{\epsilon,\g}^{-1}(M_0,N_0)\big|.\]
Moreover, based on Lemma \ref{vectorbundle3}, we have the following.
\begin{Lemma}\label{calculate-cardinality}
Let $(M_0,N_0)\in{\rm Gr}_\e(FM)\times{\rm Gr}_\f(FN)$ and $\epsilon\in\Hom_{\mathcal{C}}(M,\Sigma N)$.
\begin{enumerate}
\item[$(1)$] If $(M_0,N_0)\notin{\rm Gr}_{\e,\f}^\epsilon(FM,FN,\g)$, $\big|\psi_{\epsilon,\g}^{-1}(M_0,N_0)\big|=0$;
\item[$(2)$] If $(M_0,N_0)\in{\rm Gr}_{\e,\f}^\epsilon(FM,FN,\g)$, $\big|\psi_{\epsilon,\g}^{-1}(M_0,N_0)\big|=|\Hom_{B}(M_0,FN/N_0)|=q^{l(M,N,M_0,N_0)}$.
\end{enumerate}
\end{Lemma}

Given these calculations, if ${\rm mt}\epsilon=L$, we can refine the formula for $X_L$ to
\begin{align*}
X_L&=\int_{\g}|{\rm Gr}_\g(FL)|\cdot X^{p(L,\g)}=\int_{\g}\int_{\e,\f}|{\rm Gr}^\epsilon_\g(FL,\e,\f)|\cdot X^{p(L,\g)}\\
&=\int_{\g} \int_{\e,\f}\int_{(M_0,N_0)\in{\rm Gr}_\e(FM)\times{\rm Gr}_\f(FN)}\big|\psi_{\epsilon,\g}^{-1}(M_0,N_0)\big|\cdot X^{p(L,\g)}\\
&=\int_{\g}\int_{\e,\f}\int_{(M_0,N_0)\in{\rm Gr}_{\e,\f}^\epsilon(FM,FN,\g)}q^{l(M,N,M_0,N_0)}\cdot X^{p(L,\g)}.
\end{align*}

\begin{Remark}
Similar to the case in Section~\ref{sec3}, the refinement of $X_L$ still depends on the choice of $\epsilon$ which is not unique. But, although different $\epsilon$ lead to different refinements, they all evaluate
to $X_L$ in the end.
\end{Remark}

\subsection{Weight}
Based on the refinements of $X_L$, we introduce weight functions and weighted quantum cluster functions.

Define
\[\MG:=\biggl\{(M_0,N_0)\mid (M_0,N_0)\in\coprod_{\e,\f}{\rm Gr}_\e(FM)\times{\rm Gr}_\f(FN),\, M,N\in\mathcal{C} \biggr\}\]
and set
\[\mathbb{Z}_{\operatorname{MG}}:=\big\{f\colon\operatorname{MG} \times \EExt_{\mathcal{C}}^1 \rightarrow Z\mid  f (M_0, N_0, \epsilon)=0 \text{ unless }M={\rm qt}\epsilon, \, N={\rm st}\epsilon \big\},\]
where $\EExt_{\mathcal{C}}^1=\coprod_{M, N \in \mathcal{C}}\Hom_{\mathcal{C}} (M, \Sigma N)$ and $Z=\big\{\frac{n}{2}\mid n\in\mathbb{Z}\big\}$ is the set of all half integers.
The functions in $\mathbb{Z}_{\operatorname{MG}}$ are called weight functions.
Given $\epsilon \in \Hom_{\mathcal{C}}^1 (M, \Sigma N)$, we define
\[ \mathbb{Z}_{\MG}(\epsilon):=\{f\in \mathbb{Z}_{\MG}\mid  f(M_0, N_0, \rho)=0 \text{ if } \rho\neq \epsilon\}.\]
Given $f \in \mathbb{Z}_{\MG}(\epsilon)$, we write $f(M_0, N_0, \epsilon)$ instead as $f(M_0, N_0)$.

\begin{Definition}[weighted quantum cluster function]\label{weight-qcc}
Given a weight function $f \in \mathbb{Z}_{\MG}({\epsilon})$, the weighted quantum cluster function $f*_{\epsilon}X_L$ is a skew-polynomial in $A_{n,\lambda}$ defined as
\begin{align*}
f*_{\epsilon}X_L&=\int_{\g}\int_{\e,\f}\int_{(M_0,N_0)\in{\rm Gr}_\e(FM)\times{\rm Gr}_\f(FN)}\big|\psi_{\epsilon,\g}^{-1}(M_0,N_0)\big|\cdot q^{f(M_0,N_0)}\cdot X^{p(L,\g)}\\
&=\int_{\g}\int_{\e,\f}\int_{(M_0,N_0)\in{\rm Gr}_{\e,\f}^{\epsilon}(FM,FN,\g)}q^{l(M,N,M_0,N_0)}\cdot q^{f(M_0,N_0)}\cdot X^{p(L,\g)},
\end{align*}
where $\epsilon \in \Hom_{\mathcal{C}}(M,\Sigma N)$ with $M={\rm qt}\epsilon$, $N={\rm st}\epsilon$ and $L={\rm mt}\epsilon$.
\end{Definition}

Similarly, we can take $f$ to be the zero function to check that this definition is a $q$-deformation of the quantum cluster function.

For the zero function $f\in \mathbb{Z}_{\MG}(\epsilon)$ with ${\rm mt}\epsilon=L$, $f*_{\epsilon}X_L=X_L$.

\subsection{Multiplication}
We denote the zero morphism in $\Hom_{\mathcal{C}}(M,\Sigma N)$ by $\mathbf{0}_{MN}$. First, we consider the multiplication of quantum cluster functions.
Recall that the multiplication of skew-polynomials follows the generating relation $X^\e\cdot X^\f=q^{\frac{1}{2}\lambda(\e,\f)}X^{\e+\f}$.

So given $M$, $N$ in $\mathcal{C}$, we can multiply the quantum cluster functions as
\begin{align*}
X_M\cdot X_N&=\int_{\e}|{\rm Gr}_\e(FM)|\cdot X^{p(M,\e)}\cdot\int_{\f}|{\rm Gr}_\f(FN)|\cdot X^{p(N,\f)}\\
&=\int_{\e,\f}|{\rm Gr}_\e(FM)|\cdot|{\rm Gr}_\f(FN)|\cdot q^{\frac{1}{2}\lambda(p(M,\e),p(N,\f))}\cdot X^{p(M,\e)+p(N,\f)}\\
&=\int_{\e,\f}\int_{(M_0,N_0)\in{\rm Gr}_\e(FM)\times{\rm Gr}_\f(FN)}q^{\frac{1}{2}\lambda(p(M,\e),p(N,\f))}\cdot X^{p(M,\e)+p(N,\f)}.
\end{align*}

On the other hand, we calculate $X_{M\oplus N}$ as
\begin{align*}
X_{M\oplus N}&=\int_{\g}\int_{\e,\f}\int_{(M_0,N_0)\in{\rm Gr}_\e(FM)\times{\rm Gr}_\f(FN)}\big|\psi_{\mathbf{0}_{MN},\g}^{-1}(M_0,N_0)\big|\cdot X^{p(M\oplus N,\g)}\\
&=\int_{\g}\int_{\e,\f}\int_{(M_0,N_0)\in{\rm Gr}_{\e,\f}^{\mathbf{0}_{MN}}(FM,FN,\g)}q^{l(M,N,M_0,N_0)}\cdot X^{p(M\oplus N,\g)}.
\end{align*}

Comparing these two equalities, there are three differences: the power of $X$, the domain of integration, and the power of $q$. We analyze them successively.

Firstly, given $\epsilon\in\Hom_{\mathcal{C}}(M,\Sigma N)$ with ${\rm mt}\epsilon=L$ and submodules
\[
\begin{tikzcd}
FN \arrow[r, "Fi"] & FL \arrow[r, "Fp"] & FM \arrow[r, "F\epsilon"] & F\Sigma N, \\
N_0 \arrow[u, hook] \arrow[r] & L_0 \arrow[u, hook] \arrow[r] & M_0 \arrow[u, hook] &
\end{tikzcd}
\]
where $N_0=Fi^{-1}(L_0)$, $M_0=Fp(L_0)$, we denote the dimension vectors
of $M_0$, $N_0$ and $L_0$ by $\e$,~$\f$ and $\g$ respectively.
\begin{Lemma}[{\cite[Lemma 5.1]{Palu2008}}]
With the notation above, $p(M,\e)+p(N,\f)=p(L,\g)$.
\end{Lemma}

This lemma allows us to identify the powers of $X$ in $X_M\cdot X_N$ and $X_{M\oplus N}$, providing we are careful about the indices $\e$, $\f$ and $\g$ of the terms in these sums.

Secondly, recall the map
\[\psi_{\epsilon}=\psi_{MN}(\epsilon,-)\colon\ \coprod_{\g}{\rm Gr}_\g(FL)\longrightarrow\coprod_{\e,\f}{\rm Gr}_\e(FM)\times{\rm Gr}_\f(FN)\]
and pieces
\[\psi_{\epsilon,\g}=\psi_{\epsilon}|_{{\rm Gr}_\g(FL)}\colon\ {\rm Gr}_\g(FL)\longrightarrow\coprod_{\e,\f}{\rm Gr}_\e(FM)\times{\rm Gr}_\f(FN).\]
Although ${\rm Gr}^{\epsilon}_{M_0,N_0}(FL)$ can contain many submodules, they all have
the same dimension vector by Corollary~\ref{commondimension}. That is, the sets ${\rm Gr}_{\e,\f}^{\epsilon}(FM,FN,\g)$ do not intersect for different values of $\g$.
Thus the data $(\e,\f,M_0,N_0)$ indexing a term in the expression for $X_M\cdot X_N$ uniquely determines the additional datum $\g$ needed to index a term in the expression for $X_{M\oplus N}$. Moreover, for~$\mathbf{0}_{MN}$,
\[\psi_{\mathbf{0}_{MN}}\colon \ \coprod_{\g}{\rm Gr}_\g(FL)\longrightarrow\coprod_{\e,\f}{\rm Gr}_\e(FM)\times{\rm Gr}_\f(FN)\]
is surjective because any two submodules of $FM$ and $FN$ can be assembled into a submodule of $F(M\oplus N)$ through direct sum. So in this case, for any $\e$, $\f$,
\[{\rm Gr}_\e(FM)\times{\rm Gr}_\f(FN)=\coprod_{\g}{\rm Gr}_{\e,\f}^{\mathbf{0}_{MN}}(FM,FN,\g).\]

Then we have
\begin{align*}
X_M\cdot X_N&=\int_{\e,\f}\int_{(M_0,N_0)\in{\rm Gr}_\e(FM)\times{\rm Gr}_\f(FN)}q^{\frac{1}{2}\lambda(p(M,\e),p(N,\f))}\cdot X^{p(M,\e)+p(N,\f)}\\
&=\int_{\g}\int_{\e,\f}\int_{(M_0,N_0)\in{\rm Gr}_{\e,\f}^{\mathbf{0}_{MN}}(FM,FN,\g)}q^{\frac{1}{2}\lambda(p(M,\e),p(N,\f))}\cdot X^{p(L,\g)}.
\end{align*}

Finally, we introduce a special family of weight functions.
\begin{Definition}\label{special-weight}
For any $M$, $N$ in $\mathcal{C}$, define a weight function ${f}_{\rm spec} \in \mathbb{Z}_{\MG}$ by
\[{f}_{\rm spec}(M_0,N_0, \epsilon)=\frac{1}{2}\lambda(p(M,\e),p(N,\f))-l(M,N,M_0,N_0)\]
if $M_0\in\coprod_\e{\rm Gr}_\e(FM)$, $N_0\in\coprod_\f{\rm Gr}_\f(FN)$, $\epsilon \in \Hom_{\mathcal{C}} (M, \Sigma N)$, and $0$ otherwise.

\end{Definition}
Note that ${f}_{\rm spec}$ is constant in $\epsilon$. The following is then immediate.
\begin{Proposition}\label{direct-sum}
For any $M$, $N$ in $\mathcal{C}$, in $A_{n,\lambda}$
\[X_M\cdot X_N={f}_{\rm spec}*_{ \mathbf{0}_{MN}}X_{M\oplus N}.\]
\end{Proposition}

\begin{Definition}
Given weight functions
 $f, g \in \mathbb{Z}_{\operatorname{MG}}$, $f*_\eta g$ is defined by
 \begin{gather*}
f*_\eta g (M_0, N_0, \epsilon)=
 \begin{cases}
f(M_0,N_0, \epsilon)+g(M_0,N_0, \epsilon), &\text{if}\ \epsilon=\eta,\\
 0 , &\text{otherwise},
 \end{cases}
\end{gather*}
 for $M_0\in\coprod_{\e}{\rm Gr}_\e(FM)$ and $N_0\in\coprod_{\f}{\rm Gr}_\f(FN)$. That is to say, $f*_\eta g \in \mathbb{Z}_{\operatorname{MG}}(\eta).$
\end{Definition}
Now we consider the multiplication $(f*_{\epsilon'}X_M)\cdot(g*_{\epsilon''}X_N)$. Obviously, the product must contain information about the weight functions $f \in \mathbb{Z}_{\operatorname{MG}}(\epsilon')$ and $g \in \mathbb{Z}_{\operatorname{MG}}(\epsilon'')$. To record these, we define a corresponding weight function for the middle term.

\begin{Definition}\label{weight-func}
Given weighted quantum cluster functions $f*_{\epsilon'}X_M$ and $g*_{\epsilon''}X_N$, and a~morphism $\epsilon \in\Hom_{\mathcal{C}}(M,\Sigma N)$, define a weight function $\mathbb{T}_{fg} \in\mathbb{Z}_{\MG}(\epsilon)$ by
\[ \mathbb{T}_{fg} (M_0,N_0,\epsilon)=f(\psi_{\epsilon'}(M_0)) +g(\psi_{\epsilon''}(N_0)).\]
\end{Definition}

By Definition~\ref{weight-func}, given weighted quantum cluster functions $f*_{\epsilon'}X_M$ and $g*_{\epsilon''}X_N$, for any $\epsilon\in\Hom_{\mathcal{C}}(M,\Sigma N)$ and $\eta\in\Hom_{\mathcal{C}}(N,\Sigma M)$, we have
\[\mathbb{T}_{fg} (M_0,N_0, \epsilon)=\mathbb{T}_{gf}(N_0,M_0, \eta)\]
since both sides are equal to $f(\psi_{\epsilon'}(M_0)) +g(\psi_{\epsilon''}(N_0)).$

\begin{Proposition}\label{weight-mult-sum}
For any weighted quantum cluster functions $f*_{\epsilon'}X_M$ and $g*_{\epsilon''}X_N$, in $A_{n,\lambda}$ we have
\[(f*_{\epsilon'}X_M) \cdot(g*_{\epsilon''}X_N)=({f}_{\rm spec}*_{\mathbf{0}_{MN}}\mathbb{T}_{fg})*_{\mathbf{0}_{MN}}X_{M\oplus N}.\]
\end{Proposition}
\begin{proof}
For simplicity, without causing ambiguity, we omit some variables of weight functions in following calculation. For example,
$\mathbb{T}_{fg}(M_0,N_0, \mathbf{0}_{MN})$ is simplified to
$\mathbb{T}_{fg}$.

The key step is to calculate the fibers of the following composition of mappings
\[
\begin{tikzcd}
 & & \coprod_{\e_1}{\rm Gr}_{\e_1}(F({\rm qt}\epsilon')) \\
 & \coprod_{\e}{\rm Gr}_{\e}(FM) \arrow[rd, "p_2\circ\psi_{\epsilon'}"'] \arrow[ru, "p_1\circ\psi_{\epsilon'}"] & \times \\
 & & \coprod_{\e_2}{\rm Gr}_{\e_2}(F({\rm st}\epsilon')) \\
\coprod_{\g}{\rm Gr}_{\g}(F(M\oplus N)) \arrow[rdd, "p_2\circ\psi_{\mathbf{0}_{MN}}"'] \arrow[ruu, "p_1\circ\psi_{\mathbf{0}_{MN}}"] & \times & \\
 & & \coprod_{\f_1}{\rm Gr}_{\f_1}(F({\rm qt}\epsilon'')) \\
 & \coprod_{\f}{\rm Gr}_{\f}(FN) \arrow[rd, "p_2\circ\psi_{\epsilon''}"'] \arrow[ru, "p_1\circ\psi_{\epsilon''}"] & \times \\
 & & \coprod_{\f_1}{\rm Gr}_{\f_1}(F({\rm st}\epsilon'').
\end{tikzcd}
\]

Notice that $\psi_{\mathbf{0}_{MN}}$ is surjective, but $\psi_{\epsilon'}$ and $\psi_{\epsilon''}$ may not be surjective. We have
\begin{gather*}
\coprod_{\g}{\rm Gr}^{\mathbf{0}_{MN}}_{\e,\f}(FM,FN,\g)={\rm Gr}_{\e}(FM)\times{\rm Gr}_{\f}(FN),\\
\coprod_{\e}{\rm Gr}^{\epsilon'}_{\e_1,\e_2}(F({\rm qt}\epsilon'),F({\rm st}\epsilon''),\e)\subseteq{\rm Gr}_{\e_1}(F({\rm qt}\epsilon'))\times{\rm Gr}_{\e_2}(F({\rm st}\epsilon'')),\\
\coprod_{\f}{\rm Gr}^{\epsilon''}_{\f_1,\f_2}(F({\rm qt}\epsilon_N),F({\rm st}\epsilon_N),\f)\subseteq{\rm Gr}_{\f_1}(F({\rm qt}\epsilon''))\times{\rm Gr}_{\f_2}(F({\rm st}\epsilon'')).
\end{gather*}

By Definitions~\ref{weight-qcc} and~\ref{special-weight}, and Proposition~\ref{direct-sum}, we obtain
\begin{gather*}
({f}_{\rm spec}*_{\mathbf{0}_{MN}} \mathbb{T}_{fg} )*_{\mathbf{0}_{MN}}X_{M\oplus N}\\
\quad=\int_{\g}\int_{\e,\f}\int_{(M_0,N_0)\in{\rm Gr}^{\mathbf{0}_{MN}}_{\e,\f}(FM,FN,\g)}q^{l(M,N,M_0,N_0)+\mathbb{T}_{fg}+{f}_{\rm spec}(M_0, N_0)}\cdot X^{p(M\oplus N,\g)}\\
\quad=\int_{\e,\f}\int_{(M_0,N_0)\in{\rm Gr}_{\e}(FM)\times{\rm Gr}_{\f}(FN)}q^{\mathbb{T}_{fg}}\cdot X^{p(M,\e)}\cdot X^{p(M,\f)}\\
\quad=\int_{\e}\int_{M_0\in{\rm Gr}_{\e}(FM)}X^{p(M,\e)}\int_{\f}\int_{N_0\in{\rm Gr}_{\f}(FN)}X^{p(M,\f)}\cdot q^{\mathbb{T}_{fg}}.
\end{gather*}

Recall from Lemma~\ref{calculate-cardinality}{\samepage
\[{\rm Gr}_{\e}(FM)=\coprod_{\e_1,\e_2}{\rm Gr}^{\epsilon'}_{\e}(FM,\e_1,\e_2)
=\coprod_{\e_1,\e_2}\psi^{-1}_{\epsilon',\e}\bigl({\rm Gr}^{\epsilon'}_{\e_1,\e_2}(F({\rm qt}\epsilon'),F({\rm st}\epsilon'),\e)\bigr)\]
with the dimension of the fiber being $l({\rm qt}\epsilon',{\rm st}\epsilon',M_1,M_2)$ and
\[{\rm Gr}_{\f}(FN)=\coprod_{\f_1,\f_2}{\rm Gr}^{\epsilon''}_{\f}(FN,\f_1,\f_2)
=\coprod_{\f_1,\f_2}\psi^{-1}_{{\rm qt}\epsilon'',{\rm st}\epsilon'',\epsilon'',\f}\bigl({\rm Gr}^{\epsilon''}_{\f_1,\f_2}(F({\rm qt}\epsilon''),F({\rm st}\epsilon''),\f)\bigr)\]
with the dimension of the fiber being $l({\rm qt}\epsilon'',{\rm st}\epsilon'',N_1,N_2)$.}

Then we have
\begin{gather*}
({f}_{\rm spec}*_{\mathbf{0}_{MN}} \mathbb{T}_{fg})*_{\mathbf{0}_{MN}} X_{M\oplus N}\\
\quad=\int_{\e}\int_{M_0\in{\rm Gr}_{\e}(FM)}X^{p(M,\e)}\int_{\f}\int_{N_0\in{\rm Gr}_{\f}(FN)}X^{p(M,\f)}\cdot q^{\mathbb{T}_{fg}}\\
\quad=\int_{\e}\int_{\e_1,\e_2}\int_{(M_1,M_2)\in{\rm Gr}^{\epsilon'}_{\e_1,\e_2}(F({\rm qt}\epsilon'),F({\rm st}\epsilon'),\e)}q^{l({\rm qt}\epsilon',{\rm st}\epsilon',M_1,M_2)}\cdot X^{p(M,\e)}\\
\phantom{\quad=}\cdot \!\int_{\f}\int_{\f_1,\f_2}\int_{(N_1,N_2)\in{\rm Gr}^{\epsilon''}_{\f_1,\f_2}(F({\rm qt}\epsilon''),F({\rm st}\epsilon''),\f)}q^{l({\rm qt}\epsilon'',{\rm st}\epsilon'',N_1,N_2)}\!\cdot\! X^{p(N,\f)}
\!\cdot\! q^{f(\psi_{\epsilon'}(M_0))+g(\psi_{\epsilon''}(N_0))}\\
\quad=\int_{\e}\int_{\e_1,\e_2}\int_{(M_1,M_2)\in{\rm Gr}^{\epsilon'}_{\e_1,\e_2}(F({\rm qt}\epsilon'),F({\rm st}\epsilon'),\e)}q^{l({\rm qt}\epsilon',{\rm st}\epsilon',M_1,M_2)}
\cdot q^{f(\psi_{\epsilon'}(M_0)) }\cdot X^{p(M,\e)}\\
\phantom{\quad=}\cdot \int_{\f}\int_{\f_1,\f_2}\int_{(N_1,N_2)\in{\rm Gr}^{\epsilon''}_{\f_1,\f_2}(F({\rm qt}\epsilon''),F({\rm st}\epsilon''),\f)}q^{l({\rm qt}\epsilon'',{\rm st}\epsilon'',N_1,N_2)}
\cdot q^{g(\psi_{\epsilon''}(N_0))}\cdot X^{p(N,\f)}\\
\quad=(f*_{\epsilon'}X_M)\cdot(g*_{\epsilon''}X_N). \tag*{\qed}
\end{gather*}\renewcommand{\qed}{}
\end{proof}

\subsection[The projectivization of Hom\_C(M, Sigma N)]{The projectivization of $\boldsymbol{\Hom_{\mathcal{C}}(M,\Sigma N)}$}
Since $\Hom_{\mathcal{C}}(M,\Sigma N)$ is a finite dimensional vector space, we can consider $\mathbb{P}\Hom_{\mathcal{C}}(M,\Sigma N)$. We denote the equivalence class of $\epsilon$ in $\mathbb{P}\Hom_{\mathcal{C}}(M,\Sigma N)$ by $\mathbb{P}\epsilon$.

In this subsection, we check that multiplication of weighted quantum cluster functions is still well defined if we replace $\epsilon$ by $\mathbb{P}\epsilon$.
We assume the parameter $\lambda$ is a non-zero element in $k$.

Recall the mapping
\[\psi_{MN}\colon\ \operatorname{EG}(M,N)\longrightarrow\coprod_{\e,\f}{\rm Gr}_\e(FM)\times{\rm Gr}_\f(FN)\]
with affine fibers.

By Lemma~\ref{vectorbundle3} and the linearity of $\alpha_{M_0,N_0}$, $\epsilon\!\in\! \Sigma p\KKer\alpha_{M_0,N_0}$ if and only if $\lambda\epsilon\!\in\! \Sigma p\KKer\alpha_{M_0,N_0}$. So in this case,
\[\big|\psi^{-1}_{\epsilon,\g}(M_0,N_0)\big|=q^{l(M,N,M_0,N_0)}=\big|\psi^{-1}_{\lambda\epsilon,\g}(M_0,N_0)\big|.\]

Otherwise, they are both zero. So we have the following.
\begin{Proposition}
Given $f*_{\epsilon}X_L$, set $f \in \mathbb{Z}_{\MG}(\lambda\epsilon)$ with
$f(M_0,N_0, \epsilon)=f(M_0,N_0, \lambda\epsilon)$.
Then $f*_{\epsilon}X_L=f*_{\lambda\epsilon}X_L$.
\end{Proposition}
\begin{proof}
\begin{align*}
f*_{\epsilon}X_L&=\int_{\g}\int_{\e,\f}\int_{(M_0,N_0)\in{\rm Gr}_\e(FM)\times{\rm Gr}_\f(FN)}\big|\psi_{\epsilon,\g}^{-1}(M_0,N_0)\big|\cdot q^{f(M_0,N_0, \epsilon)}\cdot X^{p(L,\g)}\\
&=\int_{\g}\int_{\e,\f}\int_{(M_0,N_0)\in{\rm Gr}_\e(FM)\times{\rm Gr}_\f(FN)}\big|\psi_{\lambda\epsilon_L,\g}^{-1}(M_0,N_0)\big|\cdot q^{f (M_0,N_0, \lambda\epsilon)}\cdot X^{p(L,\g)}\\
&= f*_{\lambda\epsilon}X_L. \tag*{\qed}
\end{align*} \renewcommand{\qed}{}
\end{proof}

Note that given $f*_{\epsilon'}X_M$ and $g*_{\epsilon''}X_N$, by Proposition~\ref {weight-mult-sum}, we have
\[(f*_{\epsilon'}X_M) \cdot(g*_{\epsilon''}X_N)=({f}_{\rm spec}*_{\mathbf{0}_{MN}} \mathbb{T}_{fg})*_{\mathbf{0}_{MN}} X_{M\oplus N},\]
where ${f}_{\rm spec} \in \mathbb{Z}_{\MG}$ and
$\mathbb{T}_{fg}(M_0,N_0, \epsilon)=f(\psi_{\epsilon'}(M_0)) +g(\psi_{\epsilon''}(N_0)).$
If we replace $\epsilon'$ and $\epsilon''$ by their scalar multiple on the left-hand side, since $\psi_{\epsilon'}(M_0)=\psi_{\lambda\epsilon'}(M_0)$ and $\psi_{\epsilon''}(N_0)=\psi_{\mu\epsilon''}(N_0)$ for any non-zero $\lambda$ and $\mu$ in $k$,
we have
\[(f*_{\epsilon'}X_M) \cdot(g*_{\epsilon''}X_N)=(f*_{\lambda\epsilon'}X_M) \cdot(g*_{\mu\epsilon''}X_N).\]
If we replace $\epsilon'$ and $\epsilon''$ by their scalar multiple on the right-hand side, since ${f}_{\rm spec} \in \mathbb{Z}_{\MG}$ is constant in any $\epsilon$ and $\lambda\mu\mathbf{0}_{MN}=\mathbf{0}_{MN}$,
the right-hand side also remains the same.

\subsection{Multiplication formula and balanced pairs}\label{sectionofmaintheorem2}
Firstly, we introduce several special weight functions and an important property.
\begin{Definition}\label{three-weight}
There are three weight functions in $\mathbb{Z}_{\MG}$ defined as
\begin{enumerate}\itemsep=0pt
\item[(1)] $g_{\rm skew}(N_0,M_0, \eta)=\lambda(p(M,\e),p(N,\f))$, where $\lambda(-,-)$ is the skew-symmetric bilinear form defined in Section \ref{sectionofskewpolynomial};
\item[(2)] $g^+_{\rm ext}(M_0,N_0, \epsilon)=\dim_k (\operatorname{Im}\alpha'_{N_0,M_0}\cap\operatorname{Hom}_{\mathcal{C}}(N,\Sigma M))$;
\item[(3)] $g^-_{\rm ext}(N_0,M_0, \eta)=\dim_k \Sigma p\KKer\alpha_{M_0,N_0}$
\end{enumerate}
for any $M$, $N$ in $\mathcal{C}$, $(M_0,N_0)\in\coprod_{\e,\f}{\rm Gr}_\e(FM)\times{\rm Gr}_\f(FN)$, $\epsilon \in \Hom_{\mathcal{C}} (M, \Sigma N)$ and $\eta \in \Hom_{\mathcal{C}} (N, \Sigma M)$.
\end{Definition}
Note that the three weight functions are all independent
of the extension $\eta$ or $\epsilon$.

\begin{Definition}\label{pointwise2}
Given a pair of weight functions $(g^+,g^-)$, set
 \[\sigma_1(g^+):=\int_{\mathbb{P}\epsilon\in\mathbb{P}\Sigma p \KKer\alpha_{M_0,N_0}}q^{g^+(M_0, N_0, \epsilon)}\]
and
\[\sigma_2(g^-):=\int_{\mathbb{P}\eta\in\mathbb{P}(\operatorname{Im}\alpha'_{N_0,M_0}\cap\operatorname{Hom}_{\mathcal{C}}(N,\Sigma M))}q^{g^-(N_0, M_0, \eta)}.\]
This pair is called pointwise balanced if
\[\frac{q^{\dim_k\operatorname{Hom}_{\mathcal{C}}(M,\Sigma N)}-1}{q-1}=\sigma_1 (g^+)+\sigma_2(g^-)\]
holds for any $M,N\in\mathcal{C}$ and $(M_0,N_0)\in\coprod_{\e,\f}{\rm Gr}_\e(FM)\times{\rm Gr}_\f(FN)$.
\end{Definition}

\begin{Proposition}\label{mainlemma21}
The following two pairs of weight functions
\begin{enumerate}\itemsep=0pt
\item[$(1)$] $\bigl(g^+_{\rm ext}, 0\bigr)$;
\item[$(2)$] $(0,g^-_{\rm ext})$
\end{enumerate}
are pointwise balanced.
\end{Proposition}
\begin{proof}
Recall from Lemma~\ref{orthogonal2} that for any $M,N\in\mathcal{C}$ and $(M_0,N_0)\in\coprod_{\e,\f}{\rm Gr}_\e(FM)\times{\rm Gr}_\f(FN)$,
 we have
\[\dim_k \Sigma p\KKer\alpha_{M_0,N_0}+\dim_k (\operatorname{Im}\alpha'_{N_0,M_0}\cap\operatorname{Hom}_{\mathcal{C}}(N,\Sigma M))=\dim_k\operatorname{Hom}_{\mathcal{C}}(M,\Sigma N).\]

After projectivization, we have
\begin{gather*}
\frac{q^{\dim_k\operatorname{Hom}_{\mathcal{C}}(M,\Sigma N)}-1}{q-1}\\
\quad=q^{\dim_k (\operatorname{Im}\alpha'_{N_0,M_0}\cap\operatorname{Hom}_{\mathcal{C}}(N,\Sigma M))}\cdot\frac{q^{\dim_k \Sigma p\KKer\alpha_{M_0,N_0}}-1}{q-1}\\
\phantom{\quad=}{}+\frac{q^{\dim_k (\operatorname{Im}\alpha'_{N_0,M_0}\cap\operatorname{Hom}_{\mathcal{C}}(N,\Sigma M))}-1}{q-1}\\
\quad=\frac{q^{\dim_k \Sigma p\KKer\alpha_{M_0,N_0}}-1}{q-1}
+q^{\dim_k \Sigma p\KKer\alpha_{M_0,N_0}}\cdot\frac{q^{\dim_k (\operatorname{Im}\alpha'_{N_0,M_0}\cap\operatorname{Hom}_{\mathcal{C}}(N,\Sigma M))}-1}{q-1}.
\end{gather*}

So by definition,
\begin{gather*}
\sigma_1(g^+_{\rm ext})+\sigma_2(0)\\
\quad=\int_{\mathbb{P}\epsilon\in\mathbb{P}\Sigma p \KKer\alpha_{M_0,N_0}}q^{g^+_{\rm ext}(M_0,N_0, \epsilon)}+\int_{\mathbb{P}\eta\in\mathbb{P}(\operatorname{Im}\alpha'_{N_0,M_0}\cap\operatorname{Hom}_{\mathcal{C}}(N,\Sigma M))}1\\
\quad=\int_{\mathbb{P}\epsilon\in\mathbb{P}\Sigma p \KKer\alpha_{M_0,N_0}}q^{\dim_k (\operatorname{Im}\alpha'_{N_0,M_0}\cap\operatorname{Hom}_{\mathcal{C}}(N,\Sigma M))}+\int_{\mathbb{P}\eta\in\mathbb{P}(\operatorname{Im}\alpha'_{N_0,M_0}\cap\operatorname{Hom}_{\mathcal{C}}(N,\Sigma M))}1\\
\quad=q^{\dim_k (\operatorname{Im}\alpha'_{N_0,M_0}\cap\operatorname{Hom}_{\mathcal{C}}(N,\Sigma M))}\cdot\frac{q^{\dim_k \Sigma p\KKer\alpha_{M_0,N_0}}\!-1}{q-1}
+\frac{q^{\dim_k (\operatorname{Im}\alpha'_{N_0,M_0}\cap\operatorname{Hom}_{\mathcal{C}}(N,\Sigma M))}\!-1}{q-1}\\
\quad=\frac{q^{\dim_k\operatorname{Hom}_{\mathcal{C}}(M,\Sigma N)}-1}{q-1}
\end{gather*}
and
\begin{gather*}
\sigma_1(0)+\sigma_2(g^-_{\rm ext})\\
\quad=\int_{\mathbb{P}\epsilon\in\mathbb{P}\Sigma p \KKer\alpha_{M_0,N_0}}1+\int_{\mathbb{P}\eta\in\mathbb{P}(\operatorname{Im}\alpha'_{N_0,M_0}\cap\operatorname{Hom}_{\mathcal{C}}(N,\Sigma M))}q^{g^-_{\rm ext}(N_0,M_0,\eta)}\\
\quad=\int_{\mathbb{P}\epsilon\in\mathbb{P}\Sigma p \KKer\alpha_{M_0,N_0}}1+\int_{\mathbb{P}\eta\in\mathbb{P}(\operatorname{Im}\alpha'_{N_0,M_0}\cap\operatorname{Hom}_{\mathcal{C}}(N,\Sigma M))}q^{\dim_k \Sigma p\KKer\alpha_{M_0,N_0}}\\
\quad=\frac{q^{\dim_k \Sigma p\KKer\alpha_{M_0,N_0}}-1}{q-1}+q^{\dim_k \Sigma p\KKer\alpha_{M_0,N_0}}\cdot\frac{q^{\dim_k (\operatorname{Im}\alpha'_{N_0,M_0}\cap\operatorname{Hom}_{\mathcal{C}}(N,\Sigma M))}-1}{q-1}\\
\quad=\frac{q^{\dim_k\operatorname{Hom}_{\mathcal{C}}(M,\Sigma N)}-1}{q-1}. \tag*{\qed}
\end{gather*} \renewcommand{\qed}{}
\end{proof}

\begin{Theorem}\label{mainlemma22}
If a pair of weight functions $\bigl(g^+,g^-\bigr)$ in $\mathbb{Z}_{\MG}$ is pointwise balanced, then for any weighted quantum cluster functions $f*_{\epsilon'}X_M$ and $g*_{\epsilon''}X_N$ such that ${\rm Hom}_{\mathcal{C}}(M,\Sigma N)\neq 0$, we have
\begin{gather*}
|\mathbb{P}\Hom_{\mathcal{C}}(M,\Sigma N)|(f*_{\epsilon'}X_M)\cdot(g*_{\epsilon''}X_N)\\
\quad=\int_{\mathbb{P}\epsilon\in\mathbb{P}\operatorname{Hom}_{\mathcal{C}}(M,\Sigma N)}(g^+*_\epsilon {f}_{\rm spec}*_\epsilon \mathbb{T}_{fg})*_\epsilon X_{{\rm mt}\epsilon}\\
\phantom{\quad=}{}+\int_{\mathbb{P}\eta\in\mathbb{P}\operatorname{Hom}_{\mathcal{C}}(N,\Sigma M)}(g^-*_\eta g_{\rm skew}*_\eta {f}_{\rm spec}*_\eta \mathbb{T}_{gf})*_\eta X_{{\rm mt}\eta}.
\end{gather*}
\end{Theorem}
\begin{proof}
We simplify the equality in the theorem to ${\rm l.h.s.}=\Sigma_1(g^+)+\Sigma_2(g^-)$. Just for simplicity, we omit some independent variables without causing ambiguity in following calculation.
For example, $\mathbb{T}_{fg}((M_0,N_0), \epsilon)$ is simplified as $\mathbb{T}_{fg}$, $l(M,N,M_0, N_0)$ as $l$ and ${f}_{\rm spec}((M_0,N_0), \epsilon)$ as~${f}_{\rm spec}$.

Direct calculation:
\begin{align*}
{\rm l.h.s.}&=|\mathbb{P}\operatorname{Hom}_{\mathcal{C}}(M,\Sigma N)|(f*_{\epsilon'}X_M)\cdot(g*_{\epsilon''}X_N)\\
&=|\mathbb{P}\operatorname{Hom}_{\mathcal{C}}(M,\Sigma N)|({f}_{\rm spec}*_{\mathbf{0}_{MN}} \mathbb{T}_{fg})*_{\mathbf{0}_{MN}}X_{M\oplus N}\\
&=|\mathbb{P}\operatorname{Hom}_{\mathcal{C}}(M,\Sigma N)|\int_{\g}\int_{\e,\f}\int_{(M_0,N_0)\in{\rm Gr}^{\mathbf{0}_{MN}}_{\e,\f}(FM,FN,\g)}q^{l+ \mathbb{T}_{fg} +{f}_{\rm spec}}\cdot X^{p(M\oplus N,\g)}\\
&=|\mathbb{P}\operatorname{Hom}_{\mathcal{C}}(M,\Sigma N)|\int_{\e,\f}\int_{(M_0,N_0)\in{\rm Gr}_\e(FM)\times{\rm Gr}_\f(FN)}q^{\mathbb{T}_{fg}}\cdot X^{p(M,\e)}\cdot X^{p(N,\f)}\\
&=\int_{\e,\f}\int_{(M_0,N_0)\in{\rm Gr}_\e(FM)\times{\rm Gr}_\f(FN)}q^{\mathbb{T}_{fg}}\cdot X^{p(M,\e)}\cdot X^{p(N,\f)}\cdot \frac{q^{\dim_k\operatorname{Hom}_{\mathcal{C}}(M,\Sigma N)}-1}{q-1}.
\end{align*}

Recall that the third equality in the above is based on the fact
\[\coprod_{\g}{\rm Gr}^{\mathbf{0}_{MN}}_{\e,\f}(FM,FN,\g)={\rm Gr}_\e(FM)\times{\rm Gr}_\f(FN).\]

For general $\epsilon\in\Hom_{\mathcal{C}}(M,\Sigma N)$ and $\eta\in\Hom_{\mathcal{C}}(N,\Sigma M)$, we only have
\begin{gather*}
\coprod_{\g}{\rm Gr}^{\epsilon}_{\e,\f}(FM,FN,\g)\subseteq{\rm Gr}_\e(FM)\times{\rm Gr}_\f(FN),\\
\coprod_{\g}{\rm Gr}^{\eta}_{\f,\e}(FN,FM,\g)\subseteq{\rm Gr}_\f(FN)\times{\rm Gr}_\e(FM),\end{gather*}
and whether a pair $(M_0,N_0)\in{\rm Gr}_\e(FM)\times{\rm Gr}_\f(FN)$ belongs to ${\rm Gr}^{\epsilon}_{\e,\f}(FM,FN,\g)$ or ${\rm Gr}^{\eta}_{\f,\e}(FN,\allowbreak FM,\g)$ is determined by Proposition \ref{vectorbundle20}.

So we can calculate the right-hand side as
\begin{align*}
\Sigma_1(g^+)&=\int_{\mathbb{P}\epsilon\in\mathbb{P}\operatorname{Hom}_{\mathcal{C}}(M,\Sigma N)}\bigl(g^+*_\epsilon {f}_{\rm spec}*_\epsilon \mathbb{T}_{fg}\bigr)*_\epsilon X_{{\rm mt}\epsilon}\\
&=\int_{\mathbb{P}\epsilon\in\mathbb{P}\operatorname{Hom}_{\mathcal{C}}(M,\Sigma N)}\int_{\g}\int_{\e,\f}\int_{(M_0,N_0)\in{\rm Gr}_{\e,\f}^{\epsilon}(FM,FN,\g)}q^{l+g^++\mathbb{T}_{fg}+{f}_{\rm spec}}\cdot X^{p(L,\g)}\\
&=\int_{\mathbb{P}\epsilon\in\mathbb{P}\operatorname{Hom}_{\mathcal{C}}(M,\Sigma N)}\int_{\e,\f}\int_{(M_0,N_0)\in\coprod_{\g}{\rm Gr}_{\e,\f}^{\epsilon}(FM,FN,\g)}q^{g^++\mathbb{T}_{fg}}\cdot X^{p(M,\e)}\cdot X^{p(N,\f)}\\
&=\int_{\e,\f}\int_{(M_0,N_0)\in{\rm Gr}_\e(FM)\times{\rm Gr}_\f(FN)}\int_{\mathbb{P}\epsilon\in\mathbb{P}\Sigma p \KKer\alpha_{M_0,N_0}}q^{g^++\mathbb{T}_{fg}}\cdot X^{p(M,\e)}\cdot X^{p(N,\f)}\\
&=\int_{\e,\f}\int_{(M_0,N_0)\in{\rm Gr}_\e(FM)\times{\rm Gr}_\f(FN)}q^{\mathbb{T}_{fg}}\cdot X^{p(M,\e)}\cdot X^{p(N,\f)}\cdot \sigma_1(g^+).
\end{align*}

Notice that in the following calculation of $\Sigma_2(g^-)$, we denote $l=l(N,M,N_0,M_0)$,
\begin{align*}
\Sigma_2(g^-)&=\int_{\mathbb{P}\eta\in\mathbb{P}\operatorname{Hom}_{\mathcal{C}}(N,\Sigma M)}(g^-*_\eta g_{\rm skew}*_\eta {f}_{\rm spec}*_\eta \mathbb{T}_{gf})*_\eta X_{{\rm mt}\eta}\\
&=\int_{\mathbb{P}\eta\in\mathbb{P}\operatorname{Hom}_{\mathcal{C}}(N,\Sigma M)}\int_{\g}\int_{\f,\e}\int_{(N_0,M_0)\in{\rm Gr}_{\f,\e}^{\eta}(FN,FM,\g)}q^{l+g^-+ \mathbb{T}_{gf}+{f}_{\rm spec}+g_{\rm skew}}\cdot X^{p(L,\g)}\\
&=\int_{\mathbb{P}\eta\in\mathbb{P}\operatorname{Hom}_{\mathcal{C}}(N,\Sigma M)}\int_{\f,\e}\int_{(N_0,M_0)\in\coprod_{g}{\rm Gr}_{\f,\e}^{\eta}(FN,FM,\g)}q^{g^-+ \mathbb{T}_{gf}}\cdot X^{p(M,\e)}\cdot X^{p(N,\f)}\\
&=\int_{\f,\e}\int_{(N_0,M_0)\in{\rm Gr}_\f(FN)\times{\rm Gr}_\e(FM)}\int_{\mathbb{P}\eta\in\mathbb{P}(\operatorname{Im}\alpha'_{N_0,M_0}\cap\operatorname{Hom}_{\mathcal{C}}(N,\Sigma M))}q^{g^-+ \mathbb{T}_{gf}}\\
&\phantom{=}\cdot X^{p(M,\e)}\cdot X^{p(N,\f)}\\
&=\int_{\f,\e}\int_{(N_0,M_0)\in{\rm Gr}_\f(FN)\times{\rm Gr}_\e(FM)}q^{ \mathbb{T}_{gf}}\cdot X^{p(M,\e)}\cdot X^{p(N,\f)}\cdot \sigma_2 (g^-).
\end{align*}
Then by Definitions~\ref{weight-func} and ~\ref{pointwise2}, we have ${\rm l.h.s.}=\Sigma_1(g^+)+\Sigma_2(g^-)$.
\end{proof}

From Proposition \ref{mainlemma21} and Theorem \ref{mainlemma22}, we have
\begin{Theorem}[multiplication formula of weighted quantum cluster functions]\label{maintheorem2}
Let $\mathcal{C}$ be a~${\rm Hom}$-finite, $2$-Calabi--Yau, Krull--Schmidt triangulated category over a finite field $k=\mathbb{F}_q$ with a~cluster tilting object $T$. For any weighted quantum cluster functions
 $f*_{\epsilon'}X_M$ and $g*_{\epsilon''}X_N$ such that ${\rm Hom}_{\mathcal{C}}(M,\Sigma N)\neq0$, we have
\begin{gather*}
|\mathbb{P}\operatorname{Hom}_{\mathcal{C}}(M,\Sigma N)|(f*_{\epsilon'}X_M)\cdot(g*_{\epsilon''}X_N)\\
\quad=\int_{\mathbb{P}\epsilon\in\mathbb{P}\operatorname{Hom}_{\mathcal{C}}(M,\Sigma N)}(g^+_{\rm ext}*_\epsilon {f}_{\rm spec}*_\epsilon \mathbb{T}_{fg})*_\epsilon X_{{\rm mt}\epsilon}\\
\phantom{\quad=}{}+\int_{\mathbb{P}\eta\in\mathbb{P}\operatorname{Hom}_{\mathcal{C}}(N,\Sigma M)}(g_{\rm skew}*_\eta {f}_{\rm spec}*_\eta \mathbb{T}_{gf})*_\eta X_{{\rm mt}\eta}\\
\quad=\int_{\mathbb{P}\epsilon\in\mathbb{P}\operatorname{Hom}_{\mathcal{C}}(M,\Sigma N)}({f}_{\rm spec}*_\epsilon \mathbb{T}_{fg})*_\epsilon X_{{\rm mt}\epsilon}\\
\phantom{\quad=}{}+\int_{\mathbb{P}\eta\in\mathbb{P}\operatorname{Hom}_{\mathcal{C}}(N,\Sigma M)}(g^-_{\rm ext}*_\eta g_{\rm skew}*_\eta {f}_{\rm spec}*_\eta \mathbb{T}_{gf})*_\eta X_{{\rm mt}\eta}.
\end{gather*}

\end{Theorem}

Moreover, we can define a balanced pair as
\begin{Definition}
A pair of weight functions $\bigl(g^+,g^-\bigr)$ in $\mathbb{Z}_{\MG}$ is called a balanced pair if for any weighted quantum cluster functions $f*_{\epsilon'}X_M$ and $g*_{\epsilon''}X_N$ such that
${\rm Hom}_{\mathcal{C}}(M,\Sigma N)\neq 0$,
\begin{gather*}
|\mathbb{P}\operatorname{Hom}_{\mathcal{C}}(M,\Sigma N)|(f*_{\epsilon'}X_M)\cdot(g*_{\epsilon''}X_N)\\
\quad=\int_{\mathbb{P}\epsilon\in\mathbb{P}\operatorname{Hom}_{\mathcal{C}}(M,\Sigma N)}\bigl(g^+*_\epsilon {f}_{\rm spec}*_\epsilon \mathbb{T}_{fg}\bigr)*_\epsilon X_{{\rm mt}\epsilon}\\
\phantom{\quad=}{}+\int_{\mathbb{P}\eta\in\mathbb{P}\operatorname{Hom}_{\mathcal{C}}(N,\Sigma M)}(g^-*_\eta g_{\rm skew}*_\eta {f}_{\rm spec}*_\eta \mathbb{T}_{gf})*_\eta X_{{\rm mt}\eta}
\end{gather*}
holds.
\end{Definition}
\subsection[The case dim\_k Hom\_C(M,Sigma N)=1]{The case $\boldsymbol{\dim_k\Hom_{\mathcal{C}}(M,\Sigma N)=1}$}\label{sectionofonedim2}

In this subsection, we assume
\[\dim_k\Hom_{\mathcal{C}}(M,\Sigma N)=\dim_k\Hom_{\mathcal{C}}(N,\Sigma M)=1\]
and the triangles
\[N\rightarrow L\rightarrow M\xrightarrow{\epsilon} \Sigma N\qquad \text{and}\qquad M\rightarrow L'\rightarrow N\xrightarrow{\eta} \Sigma M\]
are non-split. Thus $\mathbb{P}\Hom_{\mathcal{C}}(M,\Sigma N)$ and $\mathbb{P}\Hom_{\mathcal{C}}(N,\Sigma M)$ are both singleton sets represented by $\mathbb{P}\epsilon$ and $\mathbb{P}\eta$ respectively. Since $\Sigma p\KKer\alpha_{M_0,N_0}$ and $\im\alpha'_{N_0,M_0}\cap\operatorname{Hom}_{\mathcal{C}}(N,\Sigma M)$ are orthogonal, one of them is of dimension $1$ and the other is zero.
\begin{Lemma}
With assumptions above:
\begin{enumerate}\itemsep=0pt
\item[$(1)$] If $(M_0,N_0)\in\im \psi_{\epsilon}$, $g^+_{\rm ext}(M_0,N_0, \epsilon)=0$;
\item[$(2)$] If $(N_0,M_0)\in\im \psi_{\eta}$, $g^-_{\rm ext}(N_0,M_0, \eta)=0$.
\end{enumerate}
\end{Lemma}
\begin{proof}
If $(M_0,N_0)\in\im \psi_{\epsilon}$, then $\epsilon\in\Sigma p\KKer\alpha_{M_0,N_0}$. Since $\epsilon$ is non-zero, $\Sigma p\KKer\alpha_{M_0,N_0}$ is of dimension $1$. Thus
\[g^+_{\rm ext}(M_0,N_0,\epsilon)=\dim_k(\im\alpha'_{N_0,M_0}\cap\Hom_{\mathcal{C}}(N,\Sigma M))=0.\]
On the other hand, if $(N_0,M_0)\in\im \psi_{\eta}$, then $\eta\in\im\alpha'_{N_0,M_0}\cap\Hom_{\mathcal{C}}(N,\Sigma M)$. Since $\eta$ is non-zero, we have $g^-_{\rm ext}(N_0,M_0, \eta)=\dim_k\Sigma p\KKer\alpha_{M_0,N_0}=0$.\hfill $\qed$ \renewcommand{\qed}{}
\end{proof}

\begin{Corollary}
With assumptions above,
 \[(M_0,N_0)\in\im \psi_{\epsilon} \qquad \text{if and only if} \quad (N_0,M_0)\notin\im \psi_{\eta}.\]
\end{Corollary}
\begin{proof}
Both statements are equivalent to $\dim_k\Sigma p\KKer\alpha_{M_0,N_0}=1$. \qedhere
\end{proof}

By Lemma~\ref{orthogonal2}, we have
\[\dim_k \Sigma p\KKer\alpha_{M_0,N_0}+\dim_k (\operatorname{Im}\alpha'_{N_0,M_0}\cap\operatorname{Hom}_{\mathcal{C}}(N,\Sigma M))=\dim_k\operatorname{Hom}_{\mathcal{C}}(M,\Sigma N),\]
and
\begin{align*}
&\frac{q^{\dim_k\operatorname{Hom}_{\mathcal{C}}(M,\Sigma N)}-1}{q-1}\\
&\quad=q^{\dim_k (\operatorname{Im}\alpha'_{N_0,M_0}\cap\operatorname{Hom}_{\mathcal{C}}(N,\Sigma M))}\cdot\frac{q^{\dim_k \Sigma p\KKer\alpha_{M_0,N_0}}-1}{q-1}\\
&\phantom{\quad=}{}+\frac{q^{\dim_k (\operatorname{Im}\alpha'_{N_0,M_0}\cap\operatorname{Hom}_{\mathcal{C}}(N,\Sigma M))}-1}{q-1}\\
&\quad=\frac{q^{\dim_k \Sigma p\KKer\alpha_{M_0,N_0}}-1}{q-1}
+q^{\dim_k \Sigma p\KKer\alpha_{M_0,N_0}}\cdot\frac{q^{\dim_k (\operatorname{Im}\alpha'_{N_0,M_0}\cap\operatorname{Hom}_{\mathcal{C}}(N,\Sigma M))}-1}{q-1}.
\end{align*}
So in the case $\dim_k\Hom_{\mathcal{C}}(M,\Sigma N)=1$, one of $\dim_k \Sigma p\KKer\alpha_{M_0,N_0}$ and
\[\dim_k (\operatorname{Im}\alpha'_{N_0,M_0} \cap \operatorname{Hom}_{\mathcal{C}}(N,\Sigma M))\]
 is $1$ and the other is $0$.
 Then $g^+_{\rm ext}$ and $g^-_{\rm ext}$ are zero in the relevant domain of integration and the two balanced pairs given in Proposition~\ref{mainlemma21} both degenerate to $(0,g_{\rm skew})$.

\begin{Theorem}\label{onedim2}
Assume $\dim_k\Hom_{\mathcal{C}}(M,\Sigma N)=1$ and the triangles
\[N\rightarrow L\rightarrow M\xrightarrow{\epsilon} \Sigma N \qquad \text{and}\qquad M\rightarrow L'\rightarrow N\xrightarrow{\eta} \Sigma M\]
are non-split triangles. Then we have
\[(f*_{\epsilon'}X_M)\cdot(g*_{\epsilon''}X_N)
=({f}_{\rm spec}*_\epsilon \mathbb{T}_{fg})*_\epsilon X_L+(g_{\rm skew}*_\eta {f}_{\rm spec}*_\eta \mathbb{T}_{gf})*_\eta X_{L'}.\]
\end{Theorem}

\section{Connection with preprojective algebras}\label{sec5}
In this section, we provide a quantum analogue of the connections in
\cite{GLS2012} between Palu's multiplication formula for cluster characters \cite{Palu2012} and Geiss--Leclerc--Schr\"oer's multiplication formula for evaluation forms \cite{GLS2007}.

\subsection{Preprojective algebra and nilpotent modules}

Let $k$ be a finite field and $Q=(Q_0,Q_1,s,t)$ be a finite quiver where $Q_0=\{1,\dots,n\}$ is the vertex set, and for an arrow $\alpha\colon i\rightarrow j$ in $Q_1$, set $s(\alpha)=i$ and $t(\alpha)=j$. We can obtain a new quiver $\widetilde{Q}$ from $Q$ by adding a new arrow $\bar{\alpha}\colon j\rightarrow i$ for each arrow $\alpha\colon i\rightarrow j$ in $Q_1$. Define
\[c:=\sum\limits_{\alpha\in Q_1}\alpha\bar{\alpha}-\bar{\alpha}\alpha,\]
and let $\Lambda:=k\widetilde{Q}/\langle c\rangle$ be the preprojective algebra of $Q$. For each $1\leq i \leq n$, let $S_i$ be the simple $\Lambda$-module associated to the vertex $i$.
Denote the category of all nilpotent $\Lambda$-modules by~$\operatorname{nil}\Lambda$. Let $\hat{I_i}$ be the injective envelope of $S_i$ for $1\leq i \leq n$. Given an element $\omega$ in the Weyl group associated to $Q$,
Buan, Iyama, Reiten, and Scott~\cite{BIRS} have attached to $\omega$ a 2-Calabi--Yau Frobenius
subcategory $\mathcal{C}_{\omega}\subset \operatorname{nil}\Lambda$. One can refer to~\cite[Section~2.4]{GLS2011} for a detailed description of~$\mathcal{C}_{\omega}$. We fix the element $\omega$ and a reduced expression $\mathbf{i}$.

In $\operatorname{nil}\Lambda$, there is a classical definition of a flag.
\begin{Definition}
A flag $L_{\bullet}$ of $L$ in $\operatorname{nil}\Lambda$ is a series of submodules
\[0=L_m\subseteq L_{m-1}\subseteq\dots\subseteq L_1\subseteq L_0=L.\]

Moreover, a flag $L_{\bullet}$ is called of type $(\mathbf{i},\mathbf{a})$ if $L_{j-1}/L_j$ is isomorphic to $S_{i_j}^{\oplus a_j}$ for $1\leq j \leq m$.
Denote the set of all flags of $L$ of type $(\mathbf{i},\mathbf{a})$ by $\Phi_{\mathbf{i},\mathbf{a},L}$.
\end{Definition}

Note that any flag of type $(\mathbf{i},\mathbf{a})$ can be refined to a flag of type $(\mathbf{i}',\mathbf{a}')$ with $\mathbf{a}' \in \{0,1\}^m$.
As in Section~\ref{sec2}, set
\[\widetilde{\mathcal{F}}_{\mathbf{i},\mathbf{a}, L}^{\rm mono}:=\big\{ L_m \stackrel{\iota_{L,m}}{\longrightarrow} L_{m-1} {\longrightarrow} \cdots {\longrightarrow} L_1 \stackrel{\iota_{L,1}}{\longrightarrow} L_0=L \mid L_j \in \operatorname{nil}\Lambda, \iota_{L,j} \text{ is mono}, 1\leq j\leq m\big\}.\]

Consider the action of the group $\prod_{i=0}^ m \operatorname{Aut} L_i$ on $\widetilde{\mathcal{F}}_{\mathbf{i},\mathbf{a}, L}^{\rm mono}$ as follows.
For any
\[
\tilde{g}=(g_0, g_1, \dots, g_m) \in \displaystyle\prod_{i=0}^ m \operatorname{Aut} L_i\qquad \text{and} \qquad (\iota_{L,m}, \iota_{L, m-1}, \dots, \iota_{L,1}),
\]
define
\[\tilde{g}. (\iota_{L,m}, \iota_{L, m-1}, \dots, \iota_{L,1}):=\bigl(g^{-1}_{m-1}\iota_{L,m}g_m, g^{-1}_{m-2}\iota_{L, m-1}g_{m-1}, \dots, g^{-1}_0 \iota_{L,1}g_1\bigr),\]
which can be illustrated by the commutative diagram
\[
\begin{tikzcd}
L_m \arrow[r,"\iota_{L,m}"] &[2em] L_{m-1} \arrow[r,"\iota_{L,m-1}"]&[3.5em] \cdots\arrow[r]&[2em] L_1\arrow[r,"\iota_{L,1}"]&[2em]L_0 \,\, \\
L_m \arrow[r,"g^{-1}_{m-1}\iota_{L,m}g_m"] \arrow[u,"g_m"]& L_{m-1} \arrow[r,"g^{-1}_{m-2}\iota_{L, m-1}g_{m-1}"] \arrow[u,"g_{m-1}"] &\cdots\arrow[r]&L_1\arrow[r,"g^{-1}_0 \iota_{L,1}g_1"] \arrow[u,"g_1"] &L_0 \arrow[u,"g_0"] .
\end{tikzcd}
\]
Note that ${\mathcal{F}}_{\mathbf{i},\mathbf{a}, L}^{\rm mono}$ is the set of orbits of $\widetilde{\mathcal{F}}_{\mathbf{i},\mathbf{a}, L}^{\rm mono}$ under the action of $\prod_{i=0}^ m \operatorname{Aut} L_i$ and
$\Phi_{\mathbf{i},\mathbf{a}, L}$ is the set of orbits of $\widetilde{\mathcal{F}}_{\mathbf{i},\mathbf{a}, L}^{\rm mono}$ under the action of the group
\[
\biggl\{(g_0, g_1, \dots, g_m)\in \displaystyle\prod_{i=0}^ m \operatorname{Aut} L_i |g_0=\operatorname{id}_L\biggr\}.
\]
Hence ${\mathcal{F}}_{\mathbf{i},\mathbf{a}, L}^{\rm mono}$ can be viewed as the set of orbits of $\Phi_{\mathbf{i},\mathbf{a}, L}$ under the action of the group
\[
\biggl\{(g_0, \operatorname{id}_{L_1}, \dots, \operatorname{id}_{L_m})\in \displaystyle\prod_{i=0}^ m \operatorname{Aut} L_i\biggr\} \simeq \operatorname{Aut}L.
\]
Note that for $\mathbf{c}_L \in {\mathcal{F}}_{\mathbf{i},\mathbf{a}, L}^{\rm mono}$,
the cardinality of the stabilizer of $\mathbf{c}_L$ is $q^t$ for some $t\in \mathbb{N}$ (for example, see~\cite[Section 4.1]{Riedtmann}).

The above characterizes
the relationship between $\Phi_{\mathbf{i},\mathbf{a}, L}$ and $\mathcal{F}_{\mathbf{i},\mathbf{a}, L}^{\rm mono}$ when $\mathbf{a} \in \{0,1\}^m$. In the following, we will consider $\Phi_{\mathbf{i},\mathbf{a}, L}$
instead of $\mathcal{F}_{\mathbf{i},\mathbf{a}, L}^{\rm mono}$ as Geiss--Leclerc--Schr\"oer did in~\cite{GLS2012}.

For a short exact sequence $\epsilon \colon 0 \longrightarrow N \stackrel{i}{\longrightarrow} L \stackrel{p}{\longrightarrow} M \longrightarrow 0$, define a map
\[\bar{\phi}_\epsilon \colon\ \Phi_{\mathbf{i},\mathbf{a}, L} \longrightarrow \coprod_{\mathbf{a}'+\mathbf{a}''=\mathbf{a}} \Phi_{\mathbf{i},\mathbf{a}', M} \times \Phi_{\mathbf{i},\mathbf{a}'', N}\]
which maps a flag
\begin{gather*}
f_L:=(0=L_m\subseteq L_{m-1}\subseteq\dots\subseteq L_1\subseteq L_0 = L)\in \Phi_{\mathbf{i},\mathbf{a}, L} \quad\text{to}\\
(f_M:=(0=M_m\subseteq M_{m-1}\subseteq\dots\subseteq M_1\subseteq M_0 = M),\\ f_N:=(0=N_m\subseteq N_{m-1}\subseteq\dots\subseteq N_1\subseteq N_0 = N) \in \Phi_{\mathbf{i},\mathbf{a}', M} \times \Phi_{\mathbf{i},\mathbf{a}'', N}
\end{gather*}
 with $M_i =p(L_i)$ and $N_i=i^{-1}(L_i)$ for $0\leq i \leq m$. For any
$(\mathbf{c}', \mathbf{c}'') \in \coprod_{\mathbf{a}'+\mathbf{a}''=\mathbf{a}} \Phi_{\mathbf{i},\mathbf{a}', M} \times \Phi_{\mathbf{i},\mathbf{a}'', N}$, if $\bar{\phi}^{-1}_{\epsilon} (\mathbf{c}', \mathbf{c}'')\neq \varnothing$, then
$\big|\bar{\phi}^{-1}_{\epsilon} (\mathbf{c}', \mathbf{c}'')\big|=\big|\bar{\phi}^{-1}_{{\bf 0}_{MN}} (\mathbf{c}', \mathbf{c}'')\big|$ and $\bar{\phi}^{-1}_{{\bf 0}_{MN}} (\mathbf{c}', \mathbf{c}'')$ is a vector space (see ~\cite[Lemma 3.3.1]{GLS2007}).
We set \[\bar{k} (\mathbf{c}', \mathbf{c}''):=\operatorname{dim}_k \bar{\phi}^{-1}_{{\bf 0}_{MN}} (\mathbf{c}', \mathbf{c}'').\]

\subsection{Refined socle and top series}
In this subsection, we recall some notations and definitions from~\cite[Section~3.4]{GLS2012}.

Let $L$ be an $\Lambda$-module and $S$ be a simple $\Lambda$-module. Let ${\rm soc}_SL$ be the sum of all submodules of $L$ which are isomorphic to $S$. If there exists no such $U$, set ${\rm soc}_SL=0$.
Similarly, let ${\rm top}_SL:=L/V$ where $V$ is the intersection of all submodules $U$ of $L$ such that $L/U$ are isomorphic to $S$. If there exists no such submodule, then $V=L$ and ${\rm top}_SL=0$.
Define ${\rm rad}_SL:=V$.

For $\mathbf{i}=(i_1,\dots,i_m)$, there exists a unique chain of submodules
\[0=L_m\subseteq L_{m-1}\subseteq\dots\subseteq L_1\subseteq L_0 \subseteq L\]
such that $L_{j-1}/L_j={\rm soc}_{S_{i_j}}L/L_j$. We define ${\rm soc}_{\mathbf{i}}L:=L_0$ and $L_j^+:=L_j^{+,\mathbf{i}}:=L_j$.

Moreover, the chain is denoted by
\[L_{\bullet}^+:=(L_m^+\subseteq\dots\subseteq L_0^+).\]

In particular, if ${\rm soc}_{\mathbf{i}}L=L$, $L_{\bullet}^+$ is called the refined socle series of type $\mathbf{i}$ of $L$.

Similarly, there exists a unique chain of submodules
\[0\subseteq L_m\subseteq L_{m-1}\subseteq\dots\subseteq L_1\subseteq L_0= L\]
such that $L_{j-1}/L_j={\rm top}_{S_{i_j}}L_{j-1}$. We define ${\rm top}_{\mathbf{i}}L:=L/L_m$, ${\rm rad}_{\mathbf{i}}L:=L_m$ and $L_j^-:=L_j^{-,\mathbf{i}}:=L_j$.

Moreover, the chain is denoted by
\[L_{\bullet}^-:=(L_m^-\subseteq\dots\subseteq L_0^-).\]

In particular, if ${\rm rad}_{\mathbf{i}}L=0$, $L_{\bullet}^-$ is called the refined top series of type $\mathbf{i}$ of $L$.

\begin{Lemma}[{\cite[Lemma 3.5]{GLS2012}}]
For any $L\in\mathcal{C}_{\omega}$, we have ${\rm soc}_{\mathbf{i}}L=L$ and ${\rm rad}_{\mathbf{i}}L=0$.
\end{Lemma}

\begin{Lemma}[{\cite[Lemma 3.8]{GLS2012}}]
For any $L\in\mathcal{C}_{\omega}$ and any flag \[ L_{\bullet}=(0=L_m\subseteq\dots\subseteq L_0=L)\in\Phi_{\mathbf{i},\mathbf{a},L},\] we have $L_j^-\subseteq L_j\subseteq L_j^+$ for $0 \leq j \leq m$.
\end{Lemma}

The above lemma implies that the refined socle series and the refined top series are the maximal flag and the minimal flag, respectively.

\subsection{Construction of cluster tilting objects}
Assume for each $i\in\{1,\dots,n\}$, there exists a $j\in\{1,\dots,m\}$, such that $i=i_j$.
We define
\begin{gather*}
j^-:={\rm max}\{0,1\leq s\leq j-1\mid i_s=i_j\},\qquad
j^+:={\rm min}\{j+1\leq s\leq m,m+1\mid i_s=i_j\},\\
j_{{\rm max}}:={\rm max}\{1\leq s\leq r\mid i_s=i_j\},\qquad
j_{{\rm min}}:={\rm min}\{1\leq s\leq r\mid i_s=i_j\},\\
j_i:={\rm max}\{1\leq s\leq m\mid i_s=i\}.
\end{gather*}
We define $V_0:=0$ and $V_j:={\rm soc}_{(i_1,\dots,i_j)}\hat{I}_{i_j}$ for $1\leq j\leq m$. Moreover, let $V_{\mathbf{i}}:=\bigoplus^{m}_{l=1} V_l$. For $1\leq i\leq n$, let $I_{\mathbf{i},i}:=V_{j_i}$ and $I_{\mathbf{i}}:=\bigoplus^{n}_{l=1} I_{\mathbf{i},l}$.
\begin{Remark}[{\cite[Theorems 2.9 and 2.10]{GLS2011}}]
The module $V_{\mathbf{i}}$ is a $\mathcal{C}_{\omega}$-cluster-tilting object.
\end{Remark}

Since $\mathcal{C}_\omega$ is a Frobenius category, we can consider its stable category~\cite[Theorem 2.8]{GLS2011}. Let $\mathcal{I}(M,N)$ be the subspace of $\Hom_{\Lambda}(M,N)$ consisting of all morphisms factoring through objects in ${\rm add}I_{\mathbf{i}}$ and define
\[\overline{\Hom}_{\Lambda}(M,N)=\Hom_{\Lambda}(M,N)/\mathcal{I}(M,N).\]

Then we have
\begin{Proposition}[{\cite[Proposition 3.24]{GLS2012}}]
For any $L\in\mathcal{C}_{\omega}$ and $1\leq j\leq m$,
\[D\overline{\Hom}_{\Lambda}(L,V_j)\cong e_{i_j}\bigl(L^+_j/L^-_j\bigr),\]
where $e_i$ is the primitive idempotent in $\Lambda$ associated to the vertex $i$.
\end{Proposition}

\subsection{Quiver Grassmannians}
Define $\mathcal{E}_{\mathbf{i}}:=\End_{\Lambda}(V_{\mathbf{i}})^{{\rm op}}$ and $\underline{\mathcal{E}}_{\mathbf{i}}:=\underline{\End}_{\mathcal{C}_{\omega}}(V_{\mathbf{i}})^{\rm op}$. If $k$ is algebraically closed, then $\mathcal{E}_{\mathbf{i}}$ is a finite-dimensional basic algebra and the corresponding quiver can be constructed explicitly. In this section, we make the assumption that over our
fixed finite field $k$, the algebra $\mathcal{E}_{\mathbf{i}}$ is presented by the same quiver with relations as in the case that $k$ is algebraically closed.

Define the quiver $Q_{\mathbf{i}}$ as follows: the vertex set is $\{1,\dots,m\}$; for each pair of subscripts $1\leq k$, $j\leq m$ satisfying $k^+\geq j^+\geq k> j$ and each arrow $\alpha\colon i_j\rightarrow i_k$, there is an arrow $\gamma^{k,j}_{\alpha}\colon j\rightarrow k$ called the ordinary arrow; for each $1\leq j\leq m$, there is an arrow $\gamma_j\colon j\rightarrow j^-$ if $j^-\neq 0$ called the horizontal arrow.

\begin{Proposition}[{{\cite[Proposition 3.25]{GLS2012}}}]
There is an isomorphism of quivers $Q_{\mathbf{i}}\rightarrow Q_{\mathcal{E}_{\mathbf{i}}}$ which maps $j$ to $V_j$.
\end{Proposition}

Moreover, for $L\in\mathcal{C}_{\omega}$, the $\mathcal{E}_{\mathbf{i}}$-module $D\overline{\Hom}_A(L,V_{\mathbf{i}})$ can be realized as follows:
the vector space at the vertex $j$ is \smash{$D\overline{\Hom}_A(L,V_j)=e_{i_j}\bigl(L_j^+/L_j^-\bigr)$}; for the ordinary arrow \smash{$\gamma^{k,j}_{\alpha}\colon j\rightarrow k$}, the linear map is given by
\smash{$e_{i_j}\bigl(L_j^+/L_j^-\bigr)\xrightarrow{\alpha\cdot}e_{i_k}\bigl(L_k^+/L_k^-\bigr)$}; for the horizontal arrow $\gamma_j\colon j\rightarrow j^-$, the linear map is given by
\smash{$e_{i_j}\bigl(L_j^+/L_j^-\bigr)\xrightarrow{e_{i_j}\cdot}e_{i_j}\bigl(L_{j^-}^+/L_{j^-}^-\bigr)$}.

\subsection{Bijection}



Given $L\in\mathcal{C}_{\omega}$, define the map
\begin{align*}
d_{\mathbf{i},L}\colon\ \{\mathbf{a}\in\mathbb{N}^n\mid \Phi_{\mathbf{i},\mathbf{a},L}\neq\varnothing\}&\longrightarrow\big\{g\in\mathbb{N}^n\mid {\rm Gr}^{\mathcal{E}_{\mathbf{i}}}_g(FL)\neq\varnothing\big\},\\
 (a_1,\dots,a_m) & \longmapsto (g_1,\dots g_m),
\end{align*}
where $g_j=(a_j^--a_j)+(a^-_{j^-}-a_{j_-})+\dots+(a^-_{j_{{\rm min}}}-a_{j_{{\rm min}}})$ and the map
\[{\mathrm{FG}}_{\mathbf{i},\mathbf{a},L}\colon\ \Phi_{\mathbf{i},\mathbf{a},L}\longrightarrow {\rm Gr}^{\underline{\mathcal{E}}_{\mathbf{i}}}_{d_{\mathbf{i},L}(\mathbf{a})}(FL),\]
where $FL=D\overline{\Hom}_A(L,V_{\mathbf{i}})\cong\EExt^1_{\Lambda}(W_{\mathbf{i}},L)$.
For a given $L_{\bullet}=(0=L_m\subseteq\dots\subseteq L_0=L)\in \Phi_{\mathbf{i},\mathbf{a},L}$, the image ${\mathrm{FG}}_{\mathbf{i},\mathbf{a},L}(L_{\bullet})$ is a submodule of $FL$ which can be realized as follows: the vector space at the vertex $j$ is \smash{$e_{i_j}(L_j/L^-_j)$}; for the ordinary arrow \smash{$\gamma^{k,j}_{\alpha}\colon j\rightarrow k$}, the linear map is given by
\smash{$
e_{i_j}(L_j/L_j^-)\xrightarrow{\alpha\cdot}e_{i_k}(L_k/L_k^-)
$};
for the horizontal arrow $\gamma_j\colon j\rightarrow j^-$, the linear map is given by
\smash{$
e_{i_j}(L_j/L_j^-)\xrightarrow{e_{i_j}\cdot}e_{i_j}(L_{j^-}/L_{j^-}^-)
$}.

\begin{Theorem}[{\cite[Theorem 3.27]{GLS2012}}]\label{for-comu}
For any $L\in\mathcal{C}_{\omega}$ and $\mathbf{a}\in\mathbb{N}^n$ such that $\Phi_{\mathbf{i},\mathbf{a},L}\neq\varnothing$,
the maps $d_{\mathbf{i},L}$ and ${\mathrm{FG}}_{\mathbf{i},\mathbf{a},L}$ are bijective.
\end{Theorem}

For the proof of this theorem, one can refer to~\cite[Theorem 1]{GLS2012} where Geiss, Leclerc and Schr\"oer proved an isomorphism of algebraic varieties from $\Phi_{\mathbf{i},\mathbf{a},L}$ to a certain quiver Grassmannian over the complex field $\mathbb{C}$. This result degenerates to a bijection between finite sets on a finite field $k$, therefore here we state the result without providing a detailed proof. 




The above theorem provides the relation between flags and Grassmannians. Furthermore, in order to discover relations between multiplication formulas, we need the following commutative diagram.
\begin{Theorem}\label{commutativediagram}
Let $M,N,L\in\mathcal{C}_{\omega}$ where $L={\rm mt}\epsilon$ for some $[\epsilon]\in\EExt^1_{\Lambda}(M,N)$. The short exact sequence $\epsilon$ provides a triangle $N\rightarrow L\rightarrow M\xrightarrow{\epsilon}\Sigma N$ after stabilization which we still denote by~$\epsilon$. Then there is a commutative diagram
\[
\begin{tikzcd}
{\Phi_{\mathbf{i},\mathbf{a},L}} \arrow[dd, "{\bar{\phi}_{\epsilon}}"'] \arrow[rrr, "{{\mathrm{FG}}_{\mathbf{i},\mathbf{a},L}}"] & & & {\rm Gr}^{\underline{\mathcal{E}}_{\mathbf{i}}}_{d_{\mathbf{i},L}(\mathbf{a})}(FL) \arrow[dd, "{\psi_{\epsilon}}"] \\
                                                                                                                                                                                                              & & &  \\
{\coprod\limits_{\mathbf{a}'+\mathbf{a}''=\mathbf{a}} \Phi_{\mathbf{i},\mathbf{a}',M}\times\Phi_{\mathbf{i},\mathbf{a}'',N}} \arrow[rrr, "{\coprod {\mathrm{FG}}_{\mathbf{i},\mathbf{a}',M}\times {\mathrm{FG}}_{\mathbf{i},\mathbf{a}'',N}\,\,}"] & & &{\coprod\limits_{\mathbf{a}'+\mathbf{a}''=\mathbf{a}} {\rm Gr}^{\mathcal{E}_{\mathbf{i}}}_{d_{\mathbf{i},M}(\mathbf{a}')}(FM)\times {\rm Gr}^{\mathcal{E}_{\mathbf{i}}}_{d_{\mathbf{i},N}(\mathbf{a}'')}(FN)}.
\end{tikzcd}
\]
\end{Theorem}

\begin{proof}
Consider the following diagram:
\[
\begin{tikzcd}
{\Phi_{\mathbf{i},\mathbf{a},L}} \arrow[dd, "{\bar{\phi}_{\epsilon}}"'] \arrow[rrr, "{{\mathrm{FG}}_{\mathbf{i},\mathbf{a},L}}"] & & & {{\rm Gr}^{\mathcal{E}_{\mathbf{i}}}_{d_{\mathbf{i},L}(\mathbf{a})}(FL)} \arrow[dd, "{\psi_{\epsilon}}"] \\
                                                                                                                                                                                                                  & & &  \\
{\coprod\limits_{\mathbf{a}'+\mathbf{a}''=\mathbf{a}} \Phi_{\mathbf{i},\mathbf{a}',M}\times\Phi_{\mathbf{i},\mathbf{a}'',N}} \arrow[rrr, "{\coprod {\mathrm{FG}}_{\mathbf{i},\mathbf{a}',M}\times {\mathrm{FG}}_{\mathbf{i},\mathbf{a}'',N}\,\,}"] & & & {\coprod\limits_{\mathbf{a}', \mathbf{a}''} {\rm Gr}^{\mathcal{E}_{\mathbf{i}}}_{d_{\mathbf{i},M}(\mathbf{a}')}(FM)\times {\rm Gr}^{\mathcal{E}_{\mathbf{i}}}_{d_{\mathbf{i},N}(\mathbf{a}'')}(FN)}.
\end{tikzcd}
\]
Denote ${\mathrm{FG}}_{\mathbf{i},\mathbf{a},L}$, ${\mathrm{FG}}_{\mathbf{i},\mathbf{a}',M}$ and ${\mathrm{FG}}_{\mathbf{i},\mathbf{a}'',N}$ by $\bar{F}_L$, $\bar{F}_M$ and $\bar{F}_N$ respectively. Given
$L_{\bullet}\in \Phi_{\mathbf{i},\mathbf{a},L}$, we claim that
\[
\bigl(\bar{F}_M\times \bar{F}_N\bigr) \bigl(\bar{\phi}_{\epsilon} (L_{\bullet})\bigr)=\psi_{\epsilon} \bigl(\bar{F}_L\bigr).
\]
For $\epsilon\colon N \stackrel{p}{\longrightarrow} L \stackrel{q}{\longrightarrow} M {\longrightarrow} \Sigma N$, we have
\begin{gather*}
\bar{\phi}_{\epsilon} (L_{\bullet})=\bigl(q(L_{\bullet}), p^{-1} (L_{\bullet})\bigr),\qquad
F(\epsilon)\colon\ FN \stackrel{Fp}{\longrightarrow} FL \stackrel{Fq}{\longrightarrow} FM \longrightarrow,\\
e_{i_k}(F(\epsilon))\colon\ e_{i_k}\bigl(N_k^+/N_k^-\bigr) \stackrel{p}{\longrightarrow} e_{i_k}\bigl(L_k^+/L_k^-\bigr) \stackrel{q}{\longrightarrow} e_{i_k}\bigl(M_k^+/M_k^-\bigr) {\longrightarrow}.
\end{gather*}
 Denoted by $\bar{p}$ the induced map of $p$ and $\bar{p}_k:=e_{i_k} \bar{p}$. Consider the following commutative diagram
\[
\begin{tikzcd}[sep=1.4cm]
& e_{i_k}(N_k^+/N_k^-) \arrow[r, "\bar{p}_k\,\, "] & e_{i_k}(L_k^+/L_k^-) \arrow[r, "\bar{q}_k"] & e_{i_k}(M_k^+/M_k^-) \arrow[r] & \, \, \\
& e_{i_k}(p^{-1}(L_k)/N_k^-)\arrow[r, "\bar{p}|_{p^{-1}({L_k})}"] \arrow[u, hook] &e_{i_k}(L_k/L_k^-) \arrow[r, "\bar{q}|_{L_k}"] \arrow[u, hook] & e_{i_k}(q(L_k)/M_k^-) \arrow[r] \arrow[u, hook] &\, .
\end{tikzcd}
\]
To prove $\bigl(\bar{F}_M\times \bar{F}_N\bigr) \bigl(\bar{\phi}_{\epsilon} (L_{\bullet})\bigr)=\psi_{\epsilon} \bigl(\bar{F}_L\bigr)$, it suffices to show that
\begin{enumerate}\itemsep=0pt
\item[(a)] $\bar{q}|_{L_k}$ is surjective which follows from the surjectivity of $q|_{L_k}\colon L_k \longrightarrow q(L_k)$;
\item[(b)] since $p\colon N \!\longrightarrow L$ is injective, $\Ker \bar{p}_k = p^{-1}(L_k^-) \cap N_k^+/ N_k^- =p^{-1}(L_k^-)/ N_k^- =\Ker (\bar{p}|_{p^{-1}({L_k})})$;
\item[(c)] $\Ker (\bar{q}|_{L_k})=\operatorname{Im} (\bar{p}|_{p^{-1}(L_k)})$ which follows by chasing the above diagram.
\end{enumerate}
Therefore, the above diagram is commutative.
We can obtain the commutative diagram described in the theorem.
\end{proof}

By Theorem \ref{commutativediagram}, we have
\begin{Corollary}\label{compare-two}
For any $(M_{\bullet},N_{\bullet})\in\coprod_{\mathbf{a}'+\mathbf{a}''=\mathbf{a}} \Phi_{\mathbf{i},\mathbf{a}',M}\times\Phi_{\mathbf{i},\mathbf{a}'',N}$, under the same assumptions as in Theorem~{\rm \ref{commutativediagram}}, the map
\[{\mathrm{FG}}_{\mathbf{i},\mathbf{a},L}\colon \ \bar{\phi}_{\epsilon}^{-1}(M_{\bullet},N_{\bullet})\longrightarrow\psi_{\epsilon}^{-1}({\mathrm{FG}}_{\mathbf{i},\mathbf{a}',M}(M_{\bullet}),{\mathrm{FG}}_{\mathbf{i},\mathbf{a}'', N}(N_{\bullet}))
\]
is a bijection. In particular,
\[\big|\bar{\phi}_{\epsilon}^{-1}(M_{\bullet},N_{\bullet})\big|=\big|\psi_{\epsilon}^{-1}({\mathrm{FG}}_{\mathbf{i},\mathbf{a}',M}(M_{\bullet}),{\mathrm{FG}}_{\mathbf{i},\mathbf{a}'', N}(N_{\bullet}))\big|.\]
\end{Corollary}

In order to keep compatibility with Sections \ref{sec2} and \ref{sec3}, throughout the rest of this section, we only handle $\Phi_{\mathbf{i},\mathbf{a}, L}$ with $\mathbf{a} \in \{0,1\}^m$.

\subsection[The skew-polynomial corresponding to Delta\_L]{The skew-polynomial corresponding to $\boldsymbol{\Delta_L}$}
In this subsection, we introduce a variant of the quantum cluster function defined in Section~\ref{sec3}.

Notice that $\mathcal{C}_{\omega}$ is a 2-Calabi--Yau Frobenius subcategory of $\operatorname{nil}\Lambda$ with $\Ext$-symmetry,
so that after stabilization, $\underline{\mathcal{C}_{\omega}}$ is a 2-Calabi--Yau triangulated category.

We define
\[\operatorname{FL}:=\biggl\{(\mathbf{c}',\mathbf{c}'')\mid  (\mathbf{c}', \mathbf{c}'') \in \coprod_{\mathbf{a}',\, \mathbf{a}''}\Phi_{\mathbf{i},\,\mathbf{a}',M} \times \Phi_{\mathbf{i},\mathbf{a}'',N},\,  M, N\in \mathcal{C}_{\omega}, \, \mathbf{a}', \mathbf{a}'',\mathbf{a}'+ \mathbf{a}''\in \{0,1\}^m\biggr\}\]
and set
\[\mathbb{Z}_{\operatorname{FL}}:=\{ f\colon\operatorname{FL} \times \operatorname{Exact}_{\mathcal{C}_{\omega}} \rightarrow Z\mid  f(\mathbf{c}', \mathbf{c}'', \epsilon)=0 \text{ unless }\mathbf{c}'_0={\rm qt}\epsilon, \mathbf{c}''_0={\rm st}\epsilon \},\]
where $\operatorname{Exact}_{\mathcal{C}_{\omega}}=\{\epsilon\colon 0\rightarrow N \rightarrow L \rightarrow M \rightarrow 0\mid  \epsilon \text{ is a short exact sequence, } N, L, M\in \mathcal{C}_{\omega}\}$. The functions in $\mathbb{Z}_{\operatorname{FL}}$ are called weight functions. Given $\epsilon \in \operatorname{Exact}_{\mathcal{C}_{\omega}}$, we define
\[\mathbb{Z}_{\operatorname{FL}}[\epsilon]:= \{f\in \mathbb{Z}_{\operatorname{FL}}\mid  f(\mathbf{c}', \mathbf{c}'', \rho)=0 \text{ if } \rho\neq \epsilon\}.\]
For $f\in \mathbb{Z}_{\operatorname{FL}}[\epsilon]$, we write $f(\mathbf{c}', \mathbf{c}'', \epsilon)$ instead as $f(\mathbf{c}', \mathbf{c}'')$.


\begin{Definition}\label{skew-poly}
Given a weight function $f \in \mathbb{Z}_{\operatorname{FL}}[\epsilon]$, $f*_{\epsilon}\Delta_{\mathbf{i},L}$ is the skew-polynomial in $A_{m,\lambda}$ defined by
\[f*_{\epsilon}\Delta_{\mathbf{i},L}:=\int_{\mathbf{a}}\int_{\mathbf{a}'+\mathbf{a}''=\mathbf{a}}\int_{(\mathbf{c}',\mathbf{c}'')\in \bar{\phi}_{\epsilon} (\Phi_{\mathbf{i},\mathbf{a},L})}q^{\bar{k}(\mathbf{c}',\mathbf{c}'')}\cdot q^{f(\mathbf{c}',\mathbf{c}'')}
\cdot X^{p(L,d_{\mathbf{i},L}(\mathbf{a}))},\]
where $p(L,d_{\mathbf{i},L}(\mathbf{a}))$ is defined as in Definition~\ref{pld} with $L$ being treated as an object in $\underline{\mathcal{C}_{\omega}}$.
\end{Definition}

For the convenience of calculation in the following, we use the notation \[f*_\epsilon \delta_L(\Phi_{\mathbf{i},\mathbf{a}})\colon=\int_{\mathbf{a}'+\mathbf{a}''=\mathbf{a}}\int_{(\mathbf{c}',\mathbf{c}'')\in \bar{\phi}_{\epsilon} (\Phi_{\mathbf{i},\mathbf{a},L})}q^{\bar{k}(\mathbf{c}',\mathbf{c}'')}\cdot q^{f(\mathbf{c}',\mathbf{c}'')}
\cdot X^{p(L,d_{\mathbf{i},L}(\mathbf{a}))} .\]

Notice that $A_{m,\lambda}$ is defined in Section \ref{sectionofskewpolynomial} and we take the same skew polynomial algebra for~$f*_{\epsilon}\Delta_{\mathbf{i},L}$ as the one where $f*_{\epsilon}X_L$ is
located. The skew-polynomial $f*_{\epsilon}\Delta_{\mathbf{i},L}$ is also called the weighted quantum cluster function of $L$ here.

To state the corresponding multiplication formula of $f*_{\epsilon_L}\Delta_{\mathbf{i},L}$, we need two specific weight functions.
\begin{Definition}
For any objects $M$, $N$ in $\mathcal{C}_{\omega}$, define the weight function ${f}_{\rm spec}$ by
\[{f}_{\rm spec}(M_{\bullet},N_{\bullet}, \epsilon)=\frac{1}{2}\lambda(p(M,d_{\mathbf{i},M}(\mathbf{a}')),p(N,d_{\mathbf{i},N}(\mathbf{a}'')))-\bar{k}(M_{\bullet},N_{\bullet})\]
if $M_{\bullet}\in\coprod_{\mathbf{a}'}\Phi_{\mathbf{i}, \mathbf{a}',M}$, $N_{\bullet}\in\coprod_{\mathbf{a}''}\Phi_{\mathbf{i},\mathbf{a}'',N}$, and $0$ otherwise.
\end{Definition}
\begin{Definition}
For any objects $M$, $N$ in $\mathcal{C}_{\omega}$, define the weight function $f_{\rm skew}$ as
\[f_{\rm skew}(N_{\bullet},M_{\bullet}, \epsilon)=\lambda(p(M,d_{\mathbf{i},M}(\mathbf{a}')),p(N,d_{\mathbf{i},N}(\mathbf{a}'')))\]
if $M_{\bullet}\in\coprod_{\mathbf{a}'}\Phi_{\mathbf{i},\mathbf{a}',M}$, $N_{\bullet}\in\coprod_{\mathbf{a}''}\Phi_{\mathbf{i},\mathbf{a}'',N}$, and $0$ otherwise.
\end{Definition}

Then we have
\begin{Proposition}
For any objects $M,N\in\mathcal{C}_{\omega}$ and any weighted quantum cluster functions $f*_{\epsilon'}\Delta_{\mathbf{i},M}$ and $g*_{\epsilon''}\Delta_{\mathbf{i},N}$, in $A_{m,\lambda}$
we have
\[(f*_{\epsilon'}\Delta_{\mathbf{i},M} )\cdot (g*_{\epsilon''}\Delta_{\mathbf{i},N})=({f}_{\rm spec}*_{\mathbf{0}_{MN}}\mathbb{S}_{fg}) *_{\mathbf{0}_{MN}} \Delta_{\mathbf{i},M\oplus N}.\]
\end{Proposition}
\begin{proof}
For simplicity, without causing ambiguity, we omit some variables of weight functions in following calculation. For example, ${f}_{\rm spec}((M_{\bullet},N_{\bullet}), \epsilon)$ is simplified as ${f}_{\rm spec}$.
Direct calculation shows
\begin{gather*}
(f*_{\epsilon'}\Delta_{\mathbf{i},M} )\cdot (g*_{\epsilon''}\Delta_{\mathbf{i},N})\\
\quad=\int_{\mathbf{a}',\mathbf{a}''}(f*_{\epsilon'}\delta_M(\Phi_{\mathbf{i},\mathbf{a}'}))\cdot (g*_{\epsilon''}\delta_N(\Phi_{\mathbf{i},\mathbf{a}''}))\cdot X^{p(M,d_{\mathbf{i},M}(\mathbf{a}'))}\cdot X^{p(N,d_{\mathbf{i},N}(\mathbf{a}''))}\\
\quad=\int_{\mathbf{a}}\int_{\mathbf{a}'+\mathbf{a}''=\mathbf{a}}(f*_{\epsilon'}\delta_M(\Phi_{\mathbf{i},\mathbf{a}'}))\cdot (g*_{\epsilon''}\delta_N(\Phi_{\mathbf{i},\mathbf{a}''}))\cdot q^{\frac{1}{2}\lambda(p(M,d_{\mathbf{i},M}(\mathbf{a}')),p(N,d_{\mathbf{i},N}(\mathbf{a}'')))}\\
\phantom{\quad=}{}\cdot X^{p(M,d_{\mathbf{i},M}(\mathbf{a}'))+p(N,d_{\mathbf{i},N}(\mathbf{a}''))}\\
\quad=\int_{\mathbf{a}}\int_{\mathbf{a}'+\mathbf{a}''=\mathbf{a}} (f*_{\epsilon'}\delta_M(\Phi_{\mathbf{i},\mathbf{a}'}))\cdot (g*_{\epsilon''}\delta_N(\Phi_{\mathbf{i},\mathbf{a}''}))\cdot q^{\frac{1}{2}\lambda(p(M,d_{\mathbf{i},M}(\mathbf{a}')),p(N,d_{\mathbf{i},N}(\mathbf{a}'')))}\\
\phantom{\quad=}{}\cdot X^{p(M\oplus N,d_{\mathbf{i},M\oplus N}(\mathbf{a}))}.
\end{gather*}

On the other hand,
\begin{gather*}
({f}_{\rm spec}*_{\mathbf{0}_{MN}} \mathbb{S}_{fg}) *_{\mathbf{0}_{MN}} \Delta_{\mathbf{i},M\oplus N}\\
\quad=\int_{\mathbf{a}} ({f}_{\rm spec}*_{\mathbf{0}_{MN}} \mathbb{S}_{fg})*_{\mathbf{0}_{MN}}\delta_{M\oplus N}(\Phi_{\mathbf{i},\mathbf{a}})\cdot X^{p(M\oplus N,d_{\mathbf{i},M\oplus N}(\mathbf{a}))}\\
\quad=\int_{\mathbf{a}}\int_{\mathbf{a}'+\mathbf{a}''=\mathbf{a}}\int_{(M_{\bullet},N_{\bullet})\in\Phi_{\mathbf{i},\mathbf{a}',M}\times\Phi_{\mathbf{i},\mathbf{a}'',N}}q^{\bar{k}(M_{\bullet},N_{\bullet})}\cdot q^{f+g+{f}_{\rm spec}}\cdot X^{p(M\oplus N,d_{\mathbf{i},M\oplus N}(\mathbf{a}))}\\
\quad=\int_{\mathbf{a}}\int_{\mathbf{a}'+\mathbf{a}''=\mathbf{a}}\int_{(M_{\bullet},N_{\bullet})\in\Phi_{\mathbf{i},\mathbf{a}',M}\times\Phi_{\mathbf{i},\mathbf{a}'',N}} q^{f+g}q^{\frac{1}{2}\lambda(p(M,d_{\mathbf{i},M}(\mathbf{a}')),p(N,d_{\mathbf{i},N}(\mathbf{a}'')))}\\
\phantom{\quad=}{}\cdot X^{p(M\oplus N,d_{\mathbf{i},M\oplus N}(\mathbf{a}))}.
\end{gather*}
Then by definitions of $f*_{\epsilon'}\delta_M$ and $g*_{\epsilon''}\delta_N$, we can obtain the equality.
\end{proof}

\begin{Theorem}\label{maintheorem3}
For any weighted quantum cluster functions $f*_{\epsilon'}\Delta_{\mathbf{i},M}$ and $g*_{\epsilon''}\Delta_{\mathbf{i},N}$ such that $\EExt^1_{\mathcal{C}_{\omega}}(M,N)\neq 0$, in $A_{m,\lambda}$ we have
\begin{gather*}
\big|\mathbb{P}\EExt_{\mathcal{C}_{\omega}}^1(M,N)\big| (f*_{\epsilon'}\Delta_{\mathbf{i},M} )\cdot (g*_{\epsilon''}\Delta_{\mathbf{i},N}) \\
\quad=\int_{\mathbb{P}\epsilon\in\mathbb{P}\EExt_{\mathcal{C}_{\omega}}^1(M,N)}(f^+_{\rm ext}*_\epsilon {f}_{\rm spec}*_\epsilon \mathbb{S}_{fg})*_\epsilon \Delta_{\mathbf{i},{\rm mt}\epsilon}\\
\phantom{\quad=}{}+\int_{\mathbb{P}\eta\in\mathbb{P}\EExt_{\mathcal{C}_{\omega}}^1(N,M)}(f_{\rm skew}*_\eta {f}_{\rm spec}*_\eta \mathbb{S}_{gf})*_\eta \Delta_{\mathbf{i},{\rm mt}\eta}\\
\quad=\int_{\mathbb{P}\epsilon\in\mathbb{P}\EExt_{\mathcal{C}_{\omega}}^1(M,N)}({f}_{\rm spec}*_\epsilon \mathbb{S}_{fg})*_\epsilon \Delta_{\mathbf{i},{\rm mt}\epsilon}\\
\phantom{\quad=}{}+\int_{\mathbb{P}\eta\in\mathbb{P}\EExt_{\mathcal{C}_{\omega}}^1(N,M)}(f^-_{\rm ext}*_\eta f_{\rm skew}*_\eta {f}_{\rm spec}*_\eta \mathbb{S}_{gf})*_\eta \Delta_{\mathbf{i},{\rm mt}\eta},
\end{gather*}
where $f^+_{\rm ext}$, $f^-_{\rm ext}$ and $f_{\rm skew}$ are defined as in Definition~{\rm \ref{three-weight}}.
\end{Theorem}
\begin{proof}
We only prove the first equality. The calculation for the second one is similar.

Recall that
\begin{equation*}
\begin{aligned}
{f}_{\rm spec}(M_{\bullet},N_{\bullet}, \epsilon)&=\frac{1}{2}\lambda(p(M,d_{\mathbf{i},M}(\mathbf{a}')),p(N,d_{\mathbf{i},N}(\mathbf{a}'')))-\bar{k}(M_{\bullet},N_{\bullet}),\\
{f}_{\rm spec}(N_{\bullet},M_{\bullet}, \eta)&=\frac{1}{2}\lambda(p(N,d_{\mathbf{i},N}(\mathbf{a}'')),p(M,d_{\mathbf{i},M}(\mathbf{a}')))-\bar{k}(N_{\bullet},M_{\bullet}),\\
f_{\rm skew}(N_{\bullet},M_{\bullet}, \eta )&=\lambda(p(M,d_{\mathbf{i},M}(\mathbf{a}')),p(N,d_{\mathbf{i},N}(\mathbf{a}''))),\\
f_{\rm hom}(N_{\bullet},M_{\bullet}, \eta)&=\bar{k}(M_{\bullet},N_{\bullet})-\bar{k}(N_{\bullet},M_{\bullet}).
\end{aligned}
\end{equation*}

So we have \[f_{\rm skew}(N_{\bullet},M_{\bullet}, \eta )+ {f}_{\rm spec}(N_{\bullet},M_{\bullet}, \eta)=f_{\rm hom}(N_{\bullet},M_{\bullet}, \eta)+{f}_{\rm spec}(M_{\bullet},N_{\bullet}, \epsilon).\]

Then we have
\begin{gather*}
\big|\mathbb{P}\EExt_{\mathcal{C}_{\omega}}^1(M,N)\big|
(f*_{\epsilon'}\Delta_{\mathbf{i},M} )\cdot (g*_{\epsilon''}\Delta_{\mathbf{i},N})\\
\quad=\big|\mathbb{P}\EExt_{\mathcal{C}_{\omega}}^1(M,N)\big| ({f}_{\rm spec}*_{\mathbf{0}_{MN}} \mathbb{S}_{fg})*_{\mathbf{0}_{MN}} \Delta_{\mathbf{i},M\oplus N}\\
\quad=\int_{\mathbf{a}}\big|\mathbb{P}\EExt_{\mathcal{C}_{\omega}}^1(M,N)\big| ({f}_{\rm spec}*_{\mathbf{0}_{MN}} \mathbb{S}_{fg})*_{\mathbf{0}_{MN}} \delta_{M\oplus N}(\Phi_{\mathbf{i},\mathbf{a}})\cdot X^{p(M\oplus N,d_{\mathbf{i},M\oplus N}(\mathbf{a}))},\\
\int_{\mathbb{P}\epsilon\in\mathbb{P}\EExt_{\mathcal{C}_{\omega}}^1(M,N)}\bigl(f^+_{\rm ext}*_\epsilon {f}_{\rm spec}*_\epsilon \mathbb{S}_{fg}\bigr)*_\epsilon \Delta_{\mathbf{i},{\rm mt}\epsilon}\\
\quad=\int_{\mathbf{a}}\int_{\mathbb{P}\epsilon\in\mathbb{P}\EExt_{\mathcal{C}_{\omega}}^1(M,N)}\bigl(f^+_{\rm ext}*_\epsilon {f}_{\rm spec}*_\epsilon \mathbb{S}_{fg}\bigr)*_\epsilon \delta_{{\rm mt}\epsilon}(\Phi_{\mathbf{i},\mathbf{a}})\cdot X^{p({\rm mt}\epsilon,d_{\mathbf{i},{\rm mt}\epsilon}(\mathbf{a}))},
\end{gather*}
and
\begin{gather*}
\int_{\mathbb{P}\eta\in\mathbb{P}\EExt_{\mathcal{C}_{\omega}}^1(N,M)}(f_{\rm skew}*_\eta {f}_{\rm spec}*_\eta \mathbb{S}_{gf})*_\eta \Delta_{\mathbf{i},{\rm mt}\eta}\\
\quad=\int_{\mathbf{a}}\int_{\mathbb{P}\eta\in\mathbb{P}\EExt_{\mathcal{C}_{\omega}}^1(N,M)}(f_{\rm skew}*_\eta {f}_{\rm spec}*_\eta \mathbb{S}_{gf})*_\eta \delta_{{\rm mt}\eta}(\Phi_{\mathbf{i},\mathbf{a}})\cdot X^{p({\rm mt}\eta,d_{\mathbf{i},{\rm mt}\eta}(\mathbf{a}))}.
\end{gather*}
We can rewrite
\begin{equation*}
\int_{\mathbf{a}}\int_{\mathbb{P}\eta\in\mathbb{P}\EExt_{\mathcal{C}_{\omega}}^1(N,M)}(f_{\rm skew}*_\eta {f}_{\rm spec}*_\eta \mathbb{S}_{gf})*_\eta \delta_{{\rm mt}\eta}(\Phi_{\mathbf{i},\mathbf{a}})\cdot X^{p({\rm mt}\eta,d_{\mathbf{i},{\rm mt}\eta}(\mathbf{a}))}
\end{equation*}
as
\begin{equation*}
\int_{\mathbf{a}}\int_{\mathbb{P}\eta\in\mathbb{P}\EExt_{\mathcal{C}_{\omega}}^1(N,M)}(f_{\rm hom}*_\eta {f}_{\rm spec}*_\eta \mathbb{S}_{gf})*_\eta \delta_{{\rm mt}\eta}(\Phi_{\mathbf{i},\mathbf{a}})\cdot X^{p({\rm mt}\eta,d_{\mathbf{i},{\rm mt}\eta}(\mathbf{a}))}.
\end{equation*}
But in the calculation of $q$-powers, we can replace $q^{(f_{\rm skew}*(\mathbb{S}_{gf}, \eta)*{f}_{\rm spec})}$ by $q^{(f_{\rm hom}*(\mathbb{S}_{gf},\eta)*{f}_{\rm spec})}.$ Then the desired equality is a direct consequence of Theorem \ref{maintheorem2}.
\end{proof}

\begin{Definition}
We denote the $\mathbb{Q}$-algebra generated by all weighted quantum cluster functions $f*_{\epsilon}\Delta_{\mathbf{i},L}$, where $L$ runs over $\mathcal{C}_{\omega}$ and $f \in \mathbb{Z}_{\operatorname{FL}}[{\epsilon}]$ with ${\rm mt}\epsilon=L$, by $A_q^p(\mathcal{C}_{\omega})$.
\end{Definition}

\begin{Definition}
We denote the $\mathbb{Q}$-algebra generated by $f*_{\epsilon}X_L$, where $L$ runs over $\underline{\mathcal{C}_{\omega}}$ and $f \in \mathbb{Z}_{\MG}(\epsilon)$, by $A_q(\underline{\mathcal{C}_{\omega}})$.
\end{Definition}

\subsection{Connection between two multiplication formulas}

Now we have introduced multiplication rules and proved multiplication formulas for both $f*_{\epsilon}\Delta_{\mathbf{i},L}$ and $g*_{\epsilon}X_L$. To end this section, we give the relationship between weight functions.






Given a weight function $f \in \mathbb{Z}_{\operatorname{FL}}[\epsilon]$, we define a weight function $\mathrm{FG}(f) \in \mathbb{Z}_{\operatorname{MG}}(\epsilon)$ by
\[\mathrm{FG}(f) (M_0, N_0, \epsilon) =f\bigl({\mathrm{FG}}^{-1}_{\mathbf{i},\mathbf{a}',M}(M_0) , {\mathrm{FG}}^{-1}_{\mathbf{i},\mathbf{a}'',N}(N_0)\bigr)\]
for any
\[(M_0,N_0)\in\coprod_{\e,\f}{\rm Gr}_\e(FM)\times{\rm Gr}_\f(FN),\]
where
$\mathbf{a}'=d_{{\mathbf{i}}, M}^{-1} (\e)$ and $\mathbf{a}''=d_{{\mathbf{i}}, N}^{-1} (\f)$.

Similarly, given a weight function $g \in \mathbb{Z}_{\operatorname{MG}}(\epsilon)$, we define a weight function
$\mathrm{FG}^* (g) \in \mathbb{Z}_{\operatorname{FL}}[\epsilon]$ by
\[\mathrm{FG}^* (g)(M_{\bullet},N_{\bullet}, \epsilon)=g( {\mathrm{FG}}_{\mathbf{i},\mathbf{a}',M}(M_{\bullet}) , {\mathrm{FG}}_{\mathbf{i},\mathbf{a}'',N}(N_{\bullet}), \epsilon) \]
for any $(M_{\bullet},N_{\bullet})\in \coprod_{\mathbf{a}', \mathbf{a}''}\Phi_{\mathbf{i},\mathbf{a}',M} \times \Phi_{\mathbf{i},\mathbf{a}'',N}.$

Moreover, we can extend the definitions of $\mathrm{FG}$ and $\mathrm{FG}^*$ to $A_q^p(\mathcal{C}_{\omega})$ and
$A_q(\underline{\mathcal{C}_{\omega}})$.

\begin{Definition}\label{FGandGF}
Define $\mathrm{FG}$ as a map from $A_q^p(\mathcal{C}_{\omega})$ to $A_q(\underline{\mathcal{C}}_{\omega})$ by
\[
\mathrm{FG}(f*_{\epsilon}\Delta_{\mathbf{i},L})=\mathrm{FG}(f)*_{\epsilon} X_L
\]
and $\mathrm{FG}^*$ as a map from $A_q(\underline{\mathcal{C}}_{\omega})$ to $A_q^p(\mathcal{C}_{\omega})$ by
\[
\mathrm{FG}^*(g*_{\epsilon}X_L)=\mathrm{FG}^*(g)*_{\epsilon}\Delta_{\mathbf{i},L}.
\]
\end{Definition}

\begin{Theorem}
The map $\mathrm{FG}$ is an isomorphism of algebras with the inverse $\mathrm{FG}^*$.
\end{Theorem}
\begin{proof}
In fact, if we consider $A_q^p(\mathcal{C}_{\omega})$ and $A_q(\underline{\mathcal{C}}_{\omega})$ as subalgebras of a fixed $A_{m,\lambda}$, $\mathrm{FG}$ and $\mathrm{FG}^*$ are identities. Set $\mathrm{FG}(f*_{\epsilon}\Delta_{\mathbf{i},L})=g*_{\epsilon}X_L$.

By Definition~\ref{skew-poly},
\begin{align*}
f*_{\epsilon}\Delta_{\mathbf{i},L}&=\int_{\mathbf{a}}f*_{\epsilon}\delta_L(\Phi_{\mathbf{i},\mathbf{a}})\cdot X^{p(L,d_{\mathbf{i},L}(\mathbf{a}))}\\
&=\int_{\mathbf{a}}\int_{\mathbf{a}'+\mathbf{a}''=\mathbf{a}}\int_{(M_{\bullet},N_{\bullet})\in \bar{\phi}_{\epsilon} (\Phi_{\mathbf{i},\mathbf{a},L})}q^{\bar{k}(M_{\bullet},N_{\bullet})}\cdot q^{f(M_{\bullet},N_{\bullet})}\cdot X^{p(L,d_{\mathbf{i},L}(\mathbf{a}))}
\end{align*}
and
\[g*_{\epsilon}X_L
=\int_\g\int_{\e,\f}\int_{(M_0,N_0)\in{\rm Gr}^{\epsilon}_{\e,\f}(FM,FN,\g)}q^{l(M,N,M_0,N_0)}\cdot q^{g(M_0,N_0)}\cdot X^{p(L,\g)}.\]

By Corollary~\ref{compare-two}, $\bar{k}(M_{\bullet}, N_{\bullet})=l(M,N,M_0, N_0)$ and by Definition~\ref{FGandGF},
\[q^{f(M_{\bullet}, N_{\bullet})}=q^{\mathrm{FG}(f)(M_0, N_0)}\] for
$M_0={\mathrm{FG}}_{\mathbf{i},\mathbf{a}',M}(M_{\bullet})$ and $N_0={\mathrm{FG}}_{\mathbf{i},\mathbf{a}'',N}(N_{\bullet})$.
Then we have $f*_{\epsilon}\Delta_{\mathbf{i},L}=g*_{\epsilon}X_L.$ \qedhere
\end{proof}

\begin{Remark}
By Definition~\ref{skew-poly}, when $L$ is a $\mathcal{C}_{\omega}$-projective-injective object, $FG(L)\in \mathbb{Q}$.

\end{Remark}

\section{Special version in hereditary case}\label{sec6}

In this section, we consider the cluster category from a hereditary algebra and the corresponding multiplication formula.

\subsection{Cluster category from a hereditary algebra}

Let $m\geq n$ be two positive integers and $\widetilde{Q}$ an acyclic quiver with the vertex set $\{1, \dots, m\}$. Let $Q$ be the full subquiver of $\widetilde{Q}$
with the vertex set $\{1, \dots, n\}$. Given a finite field $k$, set $\widetilde{A}=k\widetilde{Q}$ and~$A=kQ$. For any vertex $i$ of $\widetilde{Q}$ (respectively $Q$), denote by
$S_i$ the simple $\widetilde{A}$-module (respectively $A$-module) at $i$ and by $P_i$ the indecomposable projective $\widetilde{A}$-module corresponding to $i$.

Let $\widetilde{\mathcal{A}}$ be the category of finite-dimensional left $\widetilde{{A}}$-modules and the cluster category of
$\widetilde{\mathcal{A}}$ introduced by Buan--Marsh--Reineke--Reiten--Todorov \cite{BMRRT}, is defined as \smash{$\mathcal{C}\colon=\mathcal{C}_{\widetilde{\mathcal{A}}}=D^b\bigl(\widetilde{\mathcal{A}}\bigr)/\tau^{-1}\Sigma$} where \smash{$D^b\bigl(\widetilde{\mathcal{A}}\bigr)$} is the bounded derived category of $\widetilde{\mathcal{A}}$, $\tau$ is the Auslander--Reiten translation and~$\Sigma$ is the shift functor. Respectively, one can define $\mathcal{A}$ and $\mathcal{C}_{\mathcal{A}}=D^b(\mathcal{A})/\tau^{-1}\Sigma$.

Let $\widetilde{B}=(b_{ij})$ be an $m\times n$-matrix where
\[
b_{ij}=\dim_k\Ext^1_{\widetilde{\mathcal{A}}}(S_i,S_j)-\dim_k\Ext^1_{\widetilde{\mathcal{A}}}(S_j,S_i)
\]
for $1\leq i \leq m, 1\leq j \leq n.$
Assume there exists a skew-symmetric $m\times m$-matrix $\Lambda$ such that
\[\Lambda(-\widetilde{B})={\begin{bmatrix}
 I_n \\
 0
 \end{bmatrix}}_{m\times n},
\]
where $I_n$ is the $n\times n$ identity matrix.

In the following, the bilinear form $\lambda$ is always given by
$\lambda(\e,\f)=\e^T\Lambda \f$.

The cluster category $\mathcal{C}$ is a 2-Calabi--Yau triangulated category~\cite{Keller} with a cluster tilting object $T=\Sigma \widetilde{A}$.
There is a natural functor $F\colon =\operatorname{Hom}_{\mathcal{C}} (\widetilde{A}, -): \mathcal{C} \longrightarrow \widetilde{\mathcal{A}}$ which induces an
equivalence of categories $\mathcal{C}/\operatorname{add} \Sigma \widetilde{A} \stackrel{\simeq}{\longrightarrow}\widetilde{\mathcal{A}}$.

Moreover, all iso-classes of indecomposable objects in $\mathcal{C}$ can be classified by
\[
\operatorname{ind}\mathcal{C}=\operatorname{ind}\widetilde{\mathcal{A}}\cup\{\Sigma P_1,\dots,\Sigma P_m\},
\]
where $\operatorname{ind}\widetilde{\mathcal{A}}$ is the iso-classes of all indecomposable objects in $\widetilde{\mathcal{A}}$. We say $L \in \mathcal{C}$ is located in the fundamental domain if $L\in \widetilde{\mathcal{A}}$.
An object $M\in \mathcal{C}$ is called coefficient-free if $M$ does not contain a direct summand $P_i[1]$, $i>n$.

In particular, if $M, N \in \operatorname{ind}\mathcal{A}$ and $\dim\Hom_\mathcal{C}(M,\Sigma N)=1$, there exist two non-split triangles
\[
N\rightarrow L\rightarrow M\rightarrow \Sigma N \qquad\text{and}\qquad
M\rightarrow L'\rightarrow N\rightarrow \Sigma M,
\]
where one of $L$ and $L'$ is located in the fundamental domain and the other is not~\cite{Hubery}. Without loss of generality, in the following, we always assume $L'$ is located in the fundamental domain.

Now we can introduce weighted quantum cluster functions and corresponding multiplication formulas given in Section~\ref{sec4} to the cluster category $\mathcal{C}$.

\subsection{A special weighted quantum cluster function}
Let $L$ be an object with $\Bell:=\underline{\dim}_k FL$ in $\mathcal{C}$ and consider the trivial triangle
\[
L\xrightarrow{\, =\,} L\rightarrow 0 \xrightarrow{\sigma_L}\Sigma L.
\]

\begin{Definition}
Given $L$ as above, define a special weight function ${f}_{L}$ as
\[
{f}_{L}((0,L_0), \sigma_L):=-\frac{1}{2}\langle\underline{\dim}_kL_0,\Bell-\underline{\dim}_kL_0\rangle
\]
for any $L_0\in\coprod_\g{\rm Gr}_\g FL$ and $0$ otherwise, where $\langle -, -\rangle$ is the Euler form of $\widetilde{\mathcal{A}}$.
\end{Definition}

\begin{Definition} With respect to ${f}_L$, define the weighted quantum cluster function 
\[
\tilde{X}_L:={f}_L *_{\sigma_L}X_L.
\]
\end{Definition}

Notice that, for a given dimension vector $\g$ of a submodule of $FL$, the appropriate dimension vectors of ${\rm qt}\sigma_L$ and ${\rm st}\sigma_L$ are unique and obviously ${\mathbf{0}}$ and $\g$. So we can simplify the calculation of $\tilde{X}_L$ to
\[
\tilde{X}_L
=\int_\g\int_{(\mathbf{0},\g)}\int_{(0,L_0)}q^{{f}_L ((0, L_0), {\sigma_L})}\cdot X^{p(L,\g)}
=\int_\g|{\rm Gr}_\g FL|\cdot q^{-\frac{1}{2}\langle \g, \l-\g\rangle}\cdot X^{p(L,\g)}.\]

Note that Rupel \cite{Rup2011} firstly gave the above definition of quantum cluster characters over cluster categories of hereditary algebras over finite fields.
Qin provided an alternative definition of quantum cluster characters via Serre polynomials \cite{fanqin}.

\begin{Remark}
Note that $f*_{\epsilon}X_L=(f *_\epsilon {f}_L)*_{\epsilon}X_L$, with support area still decided by $\epsilon$ instead of $\sigma_L$.
\end{Remark}

\begin{Lemma}
For any $M,N\in\mathcal{C}$, we have
\[
\tilde{X}_M\cdot\tilde{X}_N=\tilde{{f}}_{MN}*_{\mathbf{0}_{MN}}\tilde{X}_{M\oplus N}
\]
 in $A_{n,\lambda}$ where $\tilde{{f}}_{MN} \in\mathbb{Z}_{\MG}(\mathbf{0}_{MN})$ is defined as
\begin{align*}
\tilde{{f}}_{MN}(M_0,N_0, \mathbf{0}_{MN})
={}&\frac{1}{2}\lambda(p(M,\e),p(N,\f))-l(M,N,M_0,N_0)+{f}_{M}((0,M_0),\sigma_M)\\
& + {f}_{N}((0,N_0),\sigma_N)
-{f}_{M\oplus N} ((0,M_0\oplus N_0), \sigma_{M\oplus N})
\end{align*}
for any $M_0\in\coprod_\e{\rm Gr}_\e(FM)$, $N_0\in\coprod_\f{\rm Gr}_\f(FN)$ and $0$ otherwise.
\end{Lemma}
\begin{proof}
We can calculate both sides as
\begin{align*}
\tilde{X}_M\cdot \tilde{X}_N
={}&\int_{\e}|{\rm Gr}_\e(FM)|\cdot q^{{f}_M((0, M_0), {\sigma_M})} \cdot X^{p(M,\e)}\cdot\int_{\f}|{\rm Gr}_\f(FN)|\cdot q^{ {f}_N((0, N_0), {\sigma_N})}\cdot X^{p(N,\f)}\\
={}&\int_{\e,\f}|{\rm Gr}_\e(FM)||{\rm Gr}_\f(FN)|\\
&\cdot q^{\frac{1}{2}\lambda(p(M,\e),p(N,\f))+{f}_M((0, M_0), {\sigma_M})+{f}_N((0, N_0), {\sigma_N}) }\cdot X^{p(M,\e)+p(N,\f)}\\
={}&\int_{\e,\f}\int_{(M_0,N_0)\in{\rm Gr}_\e(FM)\times{\rm Gr}_\f(FN)}q^{\frac{1}{2}\lambda(p(M,\e),p(N,\f))+{f}_M((0, M_0), {\sigma_M})+{f}_N((0, N_0), {\sigma_N})}\\
&\cdot X^{p(M,\e)+p(N,\f)}
\end{align*}
and
\begin{align*}
\tilde{X}_{M\oplus N}
={}&\int_{\g}\int_{\e,\f}\int_{(M_0,N_0)\in{\rm Gr}_\e(FM)\times{\rm Gr}_\f(FN)}\big|\psi_{\mathbf{0}_{MN},\g}^{-1}(M_0,N_0)\big|\cdot q^{
{f}_0((0, M_0\oplus N_0), \sigma_{M\oplus N})}\\
={}&\int_{\g}\int_{\e,\f}\int_{(M_0,N_0)\in{\rm Gr}_{\e,\f}^{\mathbf{0}_{MN}}(FM,FN,g)}q^{l(M,N,M_0,N_0)+{f}_{M\oplus N}((0, M_0\oplus N_0), \sigma_{M\oplus N})}\\
&\cdot X^{p(M\oplus N,\g)}.
\end{align*}

Then by definition of weighted quantum cluster functions, one can easily check such defined~$\tilde{{f}}_{MN}$ is the appropriate weight function to satisfy the lemma.
\end{proof}

Moreover, if we set $\underline{\dim}_kFM=\m$, $\underline{\dim}_kFN=\n$, then
\begin{Lemma}[{\cite{fanqin}}]\label{qinformula}
With the notation above, we have
\[
\frac{1}{2}\lambda(p(M,\e),p(N,\f))
=\frac{1}{2}\lambda(\operatorname{ind}M,\operatorname{ind}N)+\frac{1}{2}\langle \f,\m \rangle-\frac{1}{2}\langle \e,\n\rangle+\frac{1}{2}\langle \e,\f\rangle-\frac{1}{2}\langle \f,\e\rangle.
\]
\end{Lemma}

By this lemma, we can simplify the calculation of $\tilde{{f}}_{MN}$ to
\begin{align*}
\tilde{{f}}_{MN}(M_0,N_0, \mathbf{0}_{MN})
={}&\frac{1}{2}\lambda(p(M,\e),p(N,\f))-l(M,N,M_0,N_0)+{{f}_{M}((0,M_0), \sigma_M)}\\
&+
{{f}_{N}((0,N_0), \sigma_N)}-{{f}_{M\oplus N}((0,M_0\oplus N_0),\sigma_{M\oplus N})}\\
={}&\frac{1}{2}\lambda(\operatorname{ind}M,\operatorname{ind}N)+\frac{1}{2}\langle \f,\m \rangle-\frac{1}{2}\langle \e,\n\rangle+\frac{1}{2}\langle \e,\f\rangle-\frac{1}{2}\langle \f,\e\rangle\\
&-l(M,N,M_0,N_0)-\frac{1}{2}\langle \e,\m-\e\rangle-\frac{1}{2}\langle \f,\n-f\rangle\\
&+\frac{1}{2}\langle \e+\f,\m+\n-\e-\f\rangle\\
={}&\frac{1}{2}\lambda(\operatorname{ind}M,\operatorname{ind}N)+\langle \f,\m-\e\rangle-l(M,N,M_0,N_0).
\end{align*}

\subsection{Special version of multiplication formula}

Recall that given any balanced pair of weight functions $\bigl(g^+,g^-\bigr)$ and weighted quantum cluster functions
 $f*_{\epsilon'}X_M$ and $g*_{\epsilon''}X_N$ such that ${\rm Hom}_{\mathcal{C}}(M,\Sigma N)\neq0$, we have
\begin{gather*}
|\mathbb{P}\operatorname{Hom}_{\mathcal{C}}(M,\Sigma N)|(f*_{\epsilon'}X_M)\cdot(g*_{\epsilon''}X_N)\\
\quad=\int_{\mathbb{P}\epsilon\in\mathbb{P}\operatorname{Hom}_{\mathcal{C}}(M,\Sigma N)}\bigl(g^+*_\epsilon {f}_{\rm spec}*_\epsilon \mathbb{T}_{fg}\bigr)*_\epsilon X_{{\rm mt}\epsilon}\\
\phantom{\quad=}{}+\int_{\mathbb{P}\eta\in\mathbb{P}\operatorname{Hom}_{\mathcal{C}}(N,\Sigma M)}(g^-*_\eta g_{\rm skew}*_\eta {f}_{\rm spec}*_\eta \mathbb{T}_{gf})*_\eta X_{{\rm mt}\eta}.
\end{gather*}

In particular, we can take $\tilde{X}_M={f}_M *_{\sigma_M} X_M$ and $\tilde{X}_N={f}_N *_{\sigma_N} X_N$, then we have
\begin{align*}
|\mathbb{P}\Hom_{\mathcal{C}}(M,\Sigma N)|\tilde{X}_M\cdot\tilde{X}_N
={}&\int_{\mathbb{P}\epsilon\in\mathbb{P}\Hom_{\mathcal{C}}(M,\Sigma N)}\bigl(g^+*_\epsilon {{f}_{\rm spec}}*_\epsilon {\mathbb{T}_{{f}_M{f}_N}}\bigr)*_\epsilon X_{{\rm mt}\epsilon}\\
&+\int_{\mathbb{P}\eta\in\mathbb{P}\Hom_{\mathcal{C}}(N,\Sigma M)}(g^-*_\eta g_{\rm skew}*_\eta {{f}_{\rm spec}}*_\eta {\mathbb{T}_{{f}_N{f}_M}})*_\eta X_{{\rm mt}\eta}.
\end{align*}

To express the right-hand side in terms of $\tilde{X}_{{\rm mt}\epsilon}$ and $\tilde{X}_{{\rm mt}\eta}$, we need to extend $\tilde{{f}}_{MN}$ to the non-split case.
\begin{Definition}
Given $M, N\in\mathcal{C}$ and $\epsilon\in\Hom_{\mathcal{C}}(M,\Sigma N)$ with ${\rm mt}\epsilon=L$, define $\tilde{{f}}_{\epsilon} \in\mathbb{Z}_{\MG}(\epsilon)$ as
\begin{align*}
 \tilde{{f}}_{\epsilon}(M_0,N_0, \epsilon)
={}&\frac{1}{2}\lambda(p(M,\e),p(N,\f))-l(M,N,M_0,N_0)\\
&+{f}_M ((0,M_0),\sigma_M)+{f}_N ((0,N_0),\sigma_N)-{f}_L ((0,L_0),\sigma_L)
\end{align*}
for $(M_0,N_0)\in\operatorname{Im}\psi_{\epsilon}$ where $L_0$ satisfies $\psi_{\epsilon} (L_0)=(M_0,N_0)$, and $0$ otherwise.
\end{Definition}
\begin{Remark}\label{threeprop}
Note that
\begin{enumerate}\itemsep=0pt
\item[(1)] although the submodule $L_0$ which satisfies $\psi_{\epsilon}(L_0)=(M_0,N_0)$ is not unique, according to Corollary \ref{commondimension}, the dimension vector of $L_0$ is independent of choice, so $\tilde{{f}}_{\epsilon}$ is well defined;
\item[(2)] by definition, $\tilde{{f}}_{\mathbf{0}_{MN}}=\tilde{{f}}_{MN}$;
\item[(3)] if the triangle induced by $\epsilon$ is mapped by $F$ to a short exact sequence
\[
0\rightarrow FN\rightarrow FL\rightarrow FM\rightarrow 0,
\]
then for any $(M_0,N_0)\in\operatorname{Im}\psi_{\epsilon}$,
\[
\tilde{{f}}_{\epsilon}(M_0,N_0,\epsilon)=\tilde{{f}}_{MN}(M_0,N_0, \mathbf{0}_{MN}).
\]
\end{enumerate}
\end{Remark}

\begin{Lemma}
For any balanced pair $\bigl(g^+,g^-\bigr)$ and morphisms $\epsilon\in\Hom_{\mathcal{C}}(M,\Sigma N)$, $\eta\in\Hom_{\mathcal{C}}(N,\allowbreak\Sigma M)$ with ${\rm mt}\epsilon=L,{\rm mt}\eta=L'$, we have
\begin{gather*}
\bigl(g^+*_\epsilon {f}_{\rm spec} *_\epsilon \mathbb{T}_{{f}_M{f}_N} \bigr) *_\epsilon X_L=\bigl(g^+*_\epsilon \tilde{{f}}_{\epsilon}\bigr)*_\epsilon \tilde{X}_L,\\
 (g^-*_\eta g_{\rm skew}*_\eta {f}_{\rm spec} *_\eta \mathbb{T}_{{f}_N{f}_M}) *_\eta X_{L'}=\bigl(g^+*_\eta \tilde{{f}}_{\eta}\bigr)*_\eta \tilde{X}_{L'}.
\end{gather*}
\end{Lemma}
\begin{proof}
We only prove the first equality. The proof of the second equality is similar. By definition, the right-hand side is equal to $\bigl(g^+*_\epsilon \tilde{{f}}_{\epsilon}*_\epsilon {f}_L\bigr)*_\epsilon {X}_L $. We just need to compare the values that
both weight functions take at any $((M_0,N_0), \epsilon)$ satisfying $(M_0,N_0)\in\operatorname{Im}\psi_{\epsilon}$. We have
\begin{align*}
g^+*_\epsilon {f}_{\rm spec} *_\epsilon \mathbb{T}_{{f}_M{f}_N} (M_0,N_0,\epsilon)
={}&g^+(M_0,N_0, \epsilon)+{f}_M ((0,M_0),\sigma_M) +{f}_N ((0,N_0),\sigma_N)\\
&+\frac{1}{2}(p(M,\e),p(N,\f))-l(M,N,M_0,N_0).
\end{align*}
On the other hand,
\begin{align*}
g^+*_\epsilon \tilde{{f}}_{\epsilon}*_\epsilon {f}_L (M_0,N_0, \epsilon)
={}&g^+(M_0,N_0,\epsilon)+\frac{1}{2}\lambda(p(M,\e),p(N,\f))-l(M,N,M_0,N_0)\\
&+{f}_M ((0,M_0),\sigma_M) +{f}_N ((0,N_0),\sigma_N) \\
&- {f}_L ((0,L_0),\sigma_L) +{f}_L ((0,L_0), \sigma_L) \\
={}&g^+(M_0,N_0,\epsilon)+\frac{1}{2}\lambda(p(M,\e),p(N,\f))-l(M,N,M_0,N_0)\\
&+{f}_M ((0,M_0),\sigma_M) +{f}_N ((0,N_0),\sigma_N).
\end{align*}

Recall the definition of weighted quantum cluster function
\[
f*_{\epsilon}X_L
=\int_{\g}\int_{\e,\f}\int_{(M_0,N_0)\in{\rm Gr}_{\e,\f}^{\epsilon}(FM,FN,\g)}q^{l(M,N,M_0,N_0)}\cdot q^{f(M_0,N_0)}\cdot X^{p(L,\g)}.
\]
Since $\operatorname{Im}\psi_{\epsilon}=\coprod_\g\coprod_{\e,\f}{\rm Gr}^\epsilon_{\e,\f}(FM,FN,\g)$, the weight functions
\[g^+*_\epsilon {f}_{\rm spec} *_\epsilon \mathbb{T}_{{f}_M{f}_N}\qquad \text{and}\qquad g^+*_\epsilon \tilde{{f}}_{\epsilon}*_\epsilon {f}_L\]
have the same values on the whole domain of integration. Hence the weighted quantum cluster functions are the same.
\end{proof}

Then we can express the right-hand side of the multiplication formula in terms of $\tilde{X}_{{\rm mt}\epsilon}$ and~$\tilde{X}_{{\rm mt}\eta}$.
\begin{Theorem}
Given a balanced pair of weight functions $\bigl(g^+,g^-\bigr)$ and two objects $M,N\in\mathcal{C}$ such that ${\rm Hom}_{\mathcal{C}}(M,\Sigma N)\neq 0$, we have
\begin{align*}
&|\mathbb{P}\Hom_{\mathcal{C}}(M,\Sigma N)|\tilde{X}_M\cdot\tilde{X}_N\\
&\quad=\int_{\mathbb{P}\epsilon\in\mathbb{P}\Hom_{\mathcal{C}}(M,\Sigma N)} \bigl(g^+*_\epsilon \tilde{{f}}_{\epsilon}\bigr)*_\epsilon \tilde{X}_{{\rm mt}\epsilon}+
\int_{\mathbb{P}\eta\in\mathbb{P}\Hom_{\mathcal{C}}(N,\Sigma M)}\bigl(g^-*_\eta g_{\rm skew}*_\eta \tilde{{f}}_{\eta}\bigr)*_\eta \tilde{X}_{{\rm mt}\eta}.
\end{align*}
\end{Theorem}
Recall that Proposition \ref{mainlemma21} provides two pointwise balanced pairs of weight functions.
\begin{Corollary}
Given two objects $M$, $N$ in $\mathcal{C}$ such that ${\rm Hom}_{\mathcal{C}}(M,\Sigma N)\neq 0$, we have
\begin{gather*}
|\mathbb{P}\Hom_{\mathcal{C}}(M,\Sigma N)|\tilde{X}_M\cdot\tilde{X}_N\\
\quad=\int_{\mathbb{P}\epsilon\in\mathbb{P}\Hom_{\mathcal{C}}(M,\Sigma N)} \bigl(g^+_{\rm ext}*_\epsilon \tilde{{f}}_{\epsilon}\bigr)*_\epsilon \tilde{X}_{{\rm mt}\epsilon}+
\int_{\mathbb{P}\eta\in\mathbb{P}\Hom_{\mathcal{C}}(N,\Sigma M)}(g_{\rm skew}*_\eta \tilde{{f}}_{\eta})*_\eta \tilde{X}_{{\rm mt}\eta}\\
\quad=\int_{\mathbb{P}\epsilon\in\mathbb{P}\Hom_{\mathcal{C}}(M,\Sigma N)} \tilde{{f}}_{\epsilon}*_\epsilon \tilde{X}_{{\rm mt}\epsilon}+
\int_{\mathbb{P}\eta\in\mathbb{P}\Hom_{\mathcal{C}}(N,\Sigma M)}\bigl(g^-_{\rm ext} *g_{\rm skew}*\tilde{{f}}_{\eta}\bigr)*_\eta \tilde{X}_{{\rm mt}\eta}.
\end{gather*}
\end{Corollary}

\begin{Remark}
As this work was being completed, we became aware of a
similar result by Chen, Ding and Zhang \cite{CDZ} for hereditary categories. The exact
relationship between these two results will be investigated in the near future.
\end{Remark}

In particular, if $\dim_k\Hom_{\mathcal{C}}(M,\Sigma N)=1$, both balanced pairs degenerate to $(0,g_{\rm skew})$. Thus we have
\begin{Theorem}\label{onedimtilde}
Given two objects $M,N\in\mathcal{C}$ with $\dim_k\Hom_{\mathcal{C}}(M,\Sigma N)=1$ and two non-split triangles 
\[
N\rightarrow L\rightarrow M\xrightarrow{\epsilon}\Sigma N\qquad \text{and}\qquad M\rightarrow L'\rightarrow N\xrightarrow{\eta}\Sigma M,
\]
we have
\[\tilde{X}_M\cdot\tilde{X}_N=\tilde{{f}}_{\epsilon}*_\epsilon \tilde{X}_L+\bigl(g_{\rm skew}*\tilde{{f}}_{\eta}\bigr)*_\eta \tilde{X}_{L'}.
\]
\end{Theorem}

\subsection{Recalculation and simplification}
Theorem \ref{onedimtilde} is a special case of Theorem~\ref{onedim2} which was proved generally in Section \ref{sectionofonedim2}. But in this subsection, we calculate the right-hand side again to obtain a simple expression.
We always assume $M$, $N$ are indecomposable coefficient-free rigid objects in $\widetilde{\mathcal{A}}$ with ${\dim_k\Hom_{\mathcal{C}}(M,\Sigma N)=1}$ and
the triangles
\[N\rightarrow L\rightarrow M\xrightarrow{\epsilon} \Sigma N \qquad \text{and} \qquad M\rightarrow L'\rightarrow N\xrightarrow{\eta} \Sigma M\]
are non-split, where $L'$ is located in the fundamental domain.

By definition of weighted quantum cluster functions and Lemma \ref{qinformula},
\begin{gather*}
\tilde{{f}}_{\epsilon}*_{\epsilon} \tilde{X}_L\\
\quad=\int_{\g}\int_{\e,\f}\int_{(M_0,N_0)\in{\rm Gr}_{\e,\f}^{\epsilon}(FM,FN,\g)}q^{\frac{1}{2}\lambda(p(M,\e),p(N,\f))+ {f}_M ((0, M_0), \sigma_M) +{f}_N ((0, N_0), \sigma_N) }\cdot X^{p(L,\g)}\\
\quad=\int_{\g}\int_{\e,\f}\int_{(M_0,N_0)\in{\rm Gr}_{\e,\f}^{\epsilon}(FM,FN,\g)}\!q^{\frac{1}{2}\lambda(\operatorname{ind}M,\operatorname{ind}N)+
\langle \f,\m-\e\rangle+{f}_{M\oplus N} ((0, M_0\oplus N_0), \sigma_{M \oplus N} )} \cdot X^{p(L,\g)}.
\end{gather*}
On the other hand,
\begin{gather*}
\big(g_{\rm skew}*_\eta \tilde{{f}}_{\eta}\big)*_\eta \tilde{X}_{L'}\\
\quad=\int_{\g}\int_{\f,\e}\int_{(N_0,M_0)\in{\rm Gr}_{\f,\e}^{\eta}(FN,FM,\g)}q^{\frac{1}{2}\lambda(p(M,\e),p(N,\f))+{f}_N ((0, N_0), \sigma_N) + {f}_M ((0, M_0), \sigma_M) }\cdot X^{p(L',\g)}\\
\quad=\int_{\g}\int_{\f,\e}\int_{(N_0,M_0)\in{\rm Gr}_{\f,\e}^{\eta}(FN,FM,\g)}\!q^{\frac{1}{2}\lambda(\operatorname{ind}M,\operatorname{ind}N)+
\langle \f,\m-\e\rangle+ {f}_{N\oplus M} ((0, N_0\oplus M_0), \sigma_{N \oplus M} ) }\!\cdot X^{p(L',\g)}.
\end{gather*}

\begin{Remark}
One can check the left-hand side in Theorem \ref{onedimtilde} is
\begin{gather*}
\tilde{X}_M\cdot\tilde{X}_N\\
\quad=\int_{\e,\f}\int_{(M_0,N_0)\in{\rm Gr}_\e FM\times{\rm Gr}_\f FN}q^{\frac{1}{2}\lambda(p(M,\e),p(N,\f))+{f}_M ((0, M_0), \sigma_M) +{f}_N ((0, N_0), \sigma_N) }\\
\phantom{\quad=}{}\cdot X^{p(M,\e)+p(N,\f)}\\
\quad=\int_{\e,\f}\int_{(M_0,N_0)\in{\rm Gr}_\e FM\times{\rm Gr}_\f FN}q^{\frac{1}{2}\lambda(\operatorname{ind}M,\operatorname{ind}N)+
\langle \f,\m-\e\rangle+ {f}_{M\oplus N} ((0, M_0\oplus N_0), \sigma_{M \oplus N} )}\\
\phantom{\quad=}{}\cdot X^{p(M,\e)+p(N,\f)}.
\end{gather*}

Notice that in the case when $\dim_k\Hom_{\mathcal{C}}(M,\Sigma N)=1$, ${\rm Gr}_{\e,\f}^{\epsilon}(FM,FN,\g)$ is complementary to ${\rm Gr}_{\f,\e}^{\eta}(FN,FM,\g)$. That is to say,
\[
{\rm Gr}_\e FM\times{\rm Gr}_\f FN=\coprod_\g({\rm Gr}_{\e,\f}^{\epsilon}(FM,FN,\g)\cup{\rm Gr}_{\f,\e}^{\eta}(FN,FM,\g)).
\]
Thus the above calculation also provides a direct proof of Theorem \ref{onedimtilde}.
\end{Remark}

Now we focus on
\begin{gather*}
\bigl(g_{\rm skew}*_\eta \tilde{{f}}_{\eta}\bigr)*_\eta \tilde{X}_{L'}\\
\quad=\int_{\g}\int_{\f,\e}\int_{(N_0,M_0)\in{\rm Gr}_{\f,\e}^{\eta}(FN,FM,\g)}\!\!q^{\frac{1}{2}\lambda(\operatorname{ind}M,\operatorname{ind}N)+\langle \f,\m-\e\rangle+{f}_{N\oplus M} ((0, N_0\oplus M_0), \sigma_{N \oplus M} ) }\cdot X^{p(L',\g)}.
\end{gather*}
Recall the assumption that $L'$ is located in the fundamental domain, so
\[0\rightarrow FM\rightarrow FL'\rightarrow FN\rightarrow 0
\]
is a short exact sequence.

The following result is implicitly implied by an argument used by Qin in~\cite[Proposition~5.4.1]{fanqin}.

\begin{Lemma}
With the assumptions above, if $\psi_{\eta}(L'_0)=(N_0,M_0)$ as in the following diagram
\[
\begin{tikzcd}
& & FM/M_0 & & FN/N_0 \\
&\epsilon\colon & FM \arrow[u, "p_M", two heads] \arrow[r] & FL' \arrow[r] & FN \arrow[u, "p_N", two heads] \\
&\epsilon_0\colon & M_0 \arrow[u, "i_M", hook] \arrow[r] & L'_0 \arrow[r] \arrow[u, hook] & N_0 \arrow[u, "i_N", hook],
\end{tikzcd}
\]
then $\dim_k\Ext_{\widetilde{\mathcal{A}}}^1(N_0,FM/M_0)=0$.
\end{Lemma}
\begin{proof}
Notice that $\widetilde{A}$ is hereditary, thus $\Ext_{\widetilde{\mathcal{A}}}^2(-,-)=0$. Applying $\Hom_{\widetilde{\mathcal{A}}}(N_0,-)$ on
\[
0\rightarrow M_0\rightarrow FM\rightarrow FM/M_0\rightarrow 0,
\]
we get
\[
\Ext_{\widetilde{\mathcal{A}}}^1(N_0,M_0)\xrightarrow{\Ext_{\widetilde{\mathcal{A}}}^{1}(N_0,i_M)} \Ext_{\widetilde{\mathcal{A}}}^1(N_0,FM)\twoheadrightarrow \Ext_{\widetilde{\mathcal{A}}}^1(N_0,FM/M_0)\rightarrow0.
\]
Similarly, applying $\Hom_{\widetilde{\mathcal{A}}}(-,FM)$ on
\[
0\rightarrow N_0\rightarrow FN\rightarrow FN/N_0\rightarrow 0,
\]
we get
\[
\Ext_{\widetilde{\mathcal{A}}}^1(FN/N_0,FM)\rightarrow\Ext_{\widetilde{\mathcal{A}}}^1(FN,FM)\xrightarrow{\Ext_{\widetilde{\mathcal{A}}}^1(i_N,FM)}\Ext_{\widetilde{\mathcal{A}}}^1(N_0,FM)\rightarrow 0.
\]
Since $\dim_k \Ext_{\widetilde{\mathcal{A}}}^1(FN,FM)=1$, we have that $\dim_k \Ext_{\widetilde{\mathcal{A}}}^1(N_0,FM/M_0)$ is at most one. If $\dim_k \Ext_{\widetilde{\mathcal{A}}}^1(N_0,FM/M_0)=1$, then $\dim_k \Ext_{\widetilde{\mathcal{A}}}^1(N_0,FM)=1$, $\Ext_{\widetilde{\mathcal{A}}}^{1}(N_0,i_M)$ is zero and $\Ext_{\widetilde{\mathcal{A}}}^{1} (i_N,\allowbreak FM)$ is a linear isomorphism between one dimensional vector spaces. However, considering the construction of $\alpha'_{N_0,M_0}$, such a commutative diagram implies that $\Ext_{\widetilde{\mathcal{A}}}^1(N_0,i_M)(\epsilon_0)$ and $\Ext_{\widetilde{\mathcal{A}}}^1(i_N,FM)(\epsilon)$ coincide when $\epsilon$ is non-zero, as we are assuming. This is a contradiction.
\end{proof}

Recalling Remark \ref{threeprop}, for a short exact sequence
\[
0\rightarrow FM\rightarrow FL'\rightarrow FN\rightarrow 0
\]
and $(N_0,M_0)\in\operatorname{Im}\psi_{\eta}$, we have
\[
{f}_{N\oplus M} ((N_0,M_0), \sigma_{N\oplus M})={f}_{L'} ((N_0,M_0), \sigma_{L'}).
\]

Finally we can simplify the calculation of $\bigl(g_{\rm skew}*_\eta \tilde{{f}}_{\eta}\bigr)*_\eta \tilde{X}_{L'}$.
\begin{Lemma}\label{rhs1}
We have
\[
\bigl(g_{\rm skew}*_\eta \tilde{{f}}_{\eta}\bigr)*_\eta \tilde{X}_{L'}=q^{\frac{1}{2}\lambda(\operatorname{ind}M,\operatorname{ind}N)}\cdot \tilde{X}_{L'}.
\]
\end{Lemma}
\begin{proof}
A direct calculation shows
\begin{gather*}
 \bigl(g_{\rm skew}*_\eta \tilde{{f}}_{\eta}\bigr)*_\eta \tilde{X}_{L'} \\
\quad=\int_{\g}\int_{\f,\e}\int_{(N_0,M_0)\in{\rm Gr}_{\f,\e}^{\eta}(FN,FM,\g)}\!q^{\frac{1}{2}\lambda(\operatorname{ind}M,\operatorname{ind}N)+\langle \f,\m-\e\rangle+
{f}_{N\oplus M} ((0, N_0\oplus M_0), \sigma_{N \oplus M} ) }\cdot X^{p(L',\g)}\\
\quad=\int_{\g}\int_{\f,\e}\int_{(N_0,M_0)\in{\rm Gr}_{\f,\e}^{\eta}(FN,FM,\g)}q^{\frac{1}{2}\lambda(\operatorname{ind}M,\operatorname{ind}N)+l(N,M,N_0,M_0)+{f}_{L'} ((0, L'_0), \sigma_{L'} )}\cdot X^{p(L',\g)}\\
\quad=q^{\frac{1}{2}\lambda(\operatorname{ind}M,\operatorname{ind}N)}\cdot\int_{\g}\int_{\f,\e}\int_{(N_0,M_0)\in{\rm Gr}_{\f,\e}^{\eta}(FN,FM,\g)}q^{l(N,M,N_0,M_0)+{f}_{L'} ((0, L'_0), \sigma_{L'} )}\cdot X^{p(L',\g)}\\
\quad=q^{\frac{1}{2}\lambda(\operatorname{ind}M,\operatorname{ind}N)}\cdot {f}_{L'}*_{\sigma_{L'}}X_{L'}=q^{\frac{1}{2}\lambda(\operatorname{ind}M,\operatorname{ind}N)}\cdot \tilde{X}_{L'}. \tag*{\qed}
\end{gather*} \renewcommand{\qed}{}
\end{proof}

Notice that in the last expression above, the $q$-power is independent of the specific choice of submodules, and depends only on $M$ and $N$.

Now we analyze
\begin{gather*}
\tilde{{f}}_{\epsilon}*_{\epsilon} \tilde{X}_L\\
\quad=\int_{\g}\int_{\e,\f}\int_{(M_0,N_0)\in{\rm Gr}_{\e,\f}^{\epsilon}(FM,FN,\g)}\! q^{\frac{1}{2}\lambda(\operatorname{ind}M,\operatorname{ind}N)+\langle \f,\m-\e\rangle+
{f}_{M\oplus N} ((0, M_0\oplus N_0), \sigma_{M \oplus N} ) }\cdot X^{p(L,\g)},
\end{gather*}
where the middle term $L$ in the triangle $N\rightarrow L\rightarrow M\xrightarrow{\epsilon} \Sigma N$ is not located in the fundamental domain. In this case, there is an explicit construction of $FL$ given as follows.
Let $U=\operatorname{Im}F\epsilon$, then there exists a short exact sequence
\[
0\rightarrow V\rightarrow FM\xrightarrow{F\epsilon}U\rightarrow 0.
\] Then we have an exact sequence
\[
0\rightarrow U\rightarrow F\Sigma N\rightarrow \tau W'\oplus I\rightarrow 0.
\]
 Let $W=W'\oplus P$, there exists a sequence
 \[
 FN\xrightarrow{Fi} W\oplus V\xrightarrow{Fp} FM\xrightarrow{F\epsilon} F\Sigma N.
 \]
Finally, $FL=W\oplus V$.




Applying $F$ to the triangle $N\rightarrow L\rightarrow M\xrightarrow{\epsilon} \Sigma N$ leads to the commutative diagram
\[
\begin{tikzcd}
 & 0 & 0 & 0 & \\
 & FN/N_0 \arrow[u] \arrow[r] & FL/L_0 \arrow[u] \arrow[r] & FM/M_0 \arrow[u] \arrow[r] & \operatorname{Coker}Fp \\
 & FN \arrow[r, "Fi"] \arrow[u] & FL \arrow[r, "Fp"] \arrow[u] & FM \arrow[r, "F\epsilon"] \arrow[u] & F\Sigma N . \\
\operatorname{Ker} Fi \arrow[r] & N_0 \arrow[u] \arrow[r] & L_0 \arrow[u] \arrow[r] & M_0 \arrow[u] & \\
 & 0 \arrow[u] & 0 \arrow[u] & \,\, 0 \arrow[u] &
\end{tikzcd}
\]
Since $FL=W\oplus V$, we have a commutative diagram
\[
\begin{tikzcd}
 & 0 & 0 & 0 & \\[1mm]
0 \arrow[r] & V/J \arrow[r] \arrow[u] & FL/L_0\arrow[r] \arrow[u] & W/K \arrow[r] \arrow[u] & 0 \,\\[1mm]
0 \arrow[r] & V \arrow[r] \arrow[u] & W\oplus V \arrow[r] \arrow[u] & W\arrow[r] \arrow[u] & 0 \, \\[1mm]
0 \arrow[r] & J \arrow[r] \arrow[u] & L_0 \arrow[r] \arrow[u] & K\arrow[r] \arrow[u] & 0 .\\[1mm]
 & 0 \arrow[u] & 0 \arrow[u] & 0 \arrow[u] &
\end{tikzcd}
\]

Considering $J$ as a submodule of $V$, we have the commutative diagram
\[
\begin{tikzcd}
 & 0 \arrow[d] & 0 \arrow[d] & \\[1mm]
 & J \arrow[r, equal] \arrow[d] & J \arrow[d] & \\[1mm]
0 \arrow[r] & V \arrow[r] \arrow[d] & FM \arrow[d] \arrow[r, dashed] & U \arrow[d, equal] \\[1mm]
0 \arrow[r] & V/J \arrow[r] \arrow[d] & FM/J \arrow[d] \arrow[r, dashed] & U . \\[1mm]
 & 0 & \,\, 0 &
\end{tikzcd}
\]
On the other hand, given submodule $K$ of $W$, we can obtain an injection
\[
\tau K\rightarrow \tau W'\oplus I
\] by applying $\tau$ and using the commutative diagram
\[
\begin{tikzcd}
 & & \tau (W/K) & \\[1mm]
U \arrow[r, dashed] & F\Sigma N \arrow[r] & \tau W'\oplus I \arrow[u] \arrow[r] & 0 \\[1mm]
U \arrow[r, dashed] \arrow[u, equal] & \tau N_0 \arrow[r] \arrow[u] & \tau K \arrow[u] \arrow[r] & \, 0, \\[1mm]
 & 0 \arrow[u] \arrow[r, equal] & 0 \arrow[u] &
\end{tikzcd}
\]
where $\tau N_0$ is given by the pullback which is unique up to isomorphism.

\begin{Lemma}
There is a commutative diagram
\[
\begin{tikzcd}
 & J \arrow[r, equal] \arrow[d] & J \arrow[d] & & \tau(W/K) & \\
0 \arrow[r] & V \arrow[r] \arrow[d] & FM \arrow[d] \arrow[r, "F\epsilon"] & F\Sigma N \arrow[r] & \tau W'\oplus I \arrow[r] \arrow[u] & 0 \\
0 \arrow[r] & V/J \arrow[r] & FM/J \arrow[r, "\overline{F\epsilon}"] & \tau N_0 \arrow[u] \arrow[r] & \tau K \arrow[r] \arrow[u] &\, 0, \\
 & & & 0 \arrow[u] \arrow[r, equal] & \,\,0 \arrow[u] &
\end{tikzcd}
\]
where $\overline{F\epsilon}$ is induced by ${F\epsilon}$.
\end{Lemma}
\begin{proof}
By definition, $U=\operatorname{Im}F\epsilon$, therefore the composition of morphisms
\[
\begin{tikzcd}
FM \arrow[r, dashed] & U \arrow[r, dashed] & F\Sigma N
\end{tikzcd}
\]
is exactly $F\epsilon$.
Thus we denote the composition of morphisms
\[
\begin{tikzcd}
FM/J \arrow[r, dashed] & U \arrow[r, dashed] & \tau N_0
\end{tikzcd}
\]
by $\overline{F\epsilon}$ where $\overline{F\epsilon}(\bar{m})=F\epsilon(m)$.

Since $J\subseteq V=\operatorname{Ker}F\epsilon$, $F\epsilon(J)=0$ and $\overline{F\epsilon}(FM/J)=F\epsilon(FM)=U\subseteq \tau N_0$. Therefore, $\overline{F\epsilon}$ is well defined and the diagram
\[
\begin{tikzcd}
FM \arrow[d] \arrow[r, "F\epsilon"] & F\Sigma N \\
FM/J \arrow[r, "\overline{F\epsilon}"] & \tau N_0 \arrow[u]
\end{tikzcd}
\]
commutes. Moreover, $\Ker \overline{F\epsilon}=\Ker F\epsilon/J=V/J$ and $\im \overline{F\epsilon}=\im F\epsilon=U$.

Thus the third row in the above four-row diagram is a long exact sequence.
\end{proof}

Let $\g=\underline{\dim}\, L_0$, $\n=\underline{\dim}\,N$, $\w=\underline{\dim}\,W$, $\k=\underline{\dim}\,K $ and $\j=\underline{\dim}\,J$.
We can construct a~correspondence
\begin{align*}
\phi\colon\ \coprod_\g\coprod_{\k,\j}{\rm Gr}^{\mathbf{0}_{WV}}_{\k,\j}(W,V,\g)&\longrightarrow\coprod_\g\coprod_{\e,\f}{\rm Gr}^{\epsilon}_{\e,\f}(FM,FN,\g),\\
(K,J)&\longmapsto(J,N_0),
\end{align*}
where $J=M_0$ and $N_0$ is determined by $\tau N_0$ in the second row of the above two-row diagram.

\begin{Lemma}
The map $\phi$ is bijective.
\end{Lemma}
\begin{proof}
First we prove $\phi$ is surjective. Given a long exact sequence in $\widetilde{A}$
\[0\longrightarrow V \stackrel{i}{\longrightarrow} FM \stackrel{F\epsilon}{\longrightarrow} F\Sigma N \stackrel{\left(\begin{smallmatrix} \pi_1 \\ \pi_2 \end{smallmatrix}\right)}{\longrightarrow} \tau W' \oplus I \longrightarrow 0.\]
It gives a triangle in $\mathcal{C}$
\[N \longrightarrow W\oplus I[-1] \oplus V \longrightarrow M \longrightarrow \tau N.\]
Applying $F$ to the above triangle, we obtain the exact sequence
\[
\begin{tikzcd}
& FN \arrow[r," {\left(\begin{smallmatrix}\pi \\ 0 \end{smallmatrix}\right)}"] & W\oplus V\arrow[r, "{\left(\begin{smallmatrix} 0, i \end{smallmatrix}\right)}"] & FM.
\end{tikzcd}
\]
Given any $(J, N_0) \in \coprod_\g\coprod_{\e,\f}{\rm Gr}^{\epsilon}_{\e,\f}(FM,FN,\g)$, by definition (see Section \ref{sec4.4}), there exists
$(W_0, V_0)\in \coprod_\g\coprod_{\k,\j}{\rm Gr}^{\mathbf{0}_{WV}}_{\k,\j}(W,V,\g)$ such that
the following diagram is commutative:
\[
\begin{tikzcd}
& FN \arrow[r," {\left(\begin{smallmatrix} \pi \\ 0 \end{smallmatrix}\right)}"] & W\oplus V\arrow[r, "{\left(\begin{smallmatrix} 0, i \end{smallmatrix}\right)}"] & FM \, \\
& N_0 \arrow[u, "i_M", hook] \arrow[r] & W_0 \oplus V_0 \arrow[r] \arrow[u, hook] & J \arrow[u, "i_N", hook] .
\end{tikzcd}
\]
Hence, $W_0=\pi_1(N_0)$, $V=J$ and then $\phi(\pi_1(N_0), J)=(J, N_0)$.

We now prove $\phi$ is injective. By the definition of $\phi$, $N_0=\pi_1^{-1} (K)$ and $K=\pi_1(N_0)$. Therefore, there is a unique $(K, J)$ mapped to the given
$(J, N_0)$. \hfill $\qed$ \renewcommand{\qed}{}
\end{proof}

\begin{Lemma}
We have that $\underline{\dim}\, \phi (K, J):=(\underline{\dim}\,J, \underline{\dim}\, N_0)=(\j,\n-(\w-\k))$.
\end{Lemma}
\begin{proof}
Given $K \subseteq W$, consider the commutative diagram
\[
\begin{tikzcd}
0 \arrow[r] & U \arrow[r] & F\Sigma N \arrow[r] & \tau W'\oplus I \arrow[r] & 0 \\
0 \arrow[r] & U \arrow[r] \arrow[u, equal] & \tau N_0 \arrow[r] \arrow[u] & \tau K \arrow[u] \arrow[r] & \,\,0.
\end{tikzcd}
\]
Since both rows are short exact sequences, we can compute dimension vectors after applying~$\tau^{-1}$~as
\[
\underline{\dim}_kN_0=\underline{\dim}_kFN+\underline{\dim}_kK-\underline{\dim}_kW=\n-(\w-\k). \tag*{\qed}
\]\renewcommand{\qed}{}
\end{proof}

To complete the final calculation, we need an identity given in \cite{fanqin}.
\begin{Lemma}[{\cite{fanqin}}]
With the notation above, we have
\begin{gather*}
\langle \n-(\w-\k),\m-\j\rangle-\frac{1}{2}\langle \e+\f,\m+\n-\e-\f\rangle\\
\quad=\langle \k,\v-\j\rangle-\frac{1}{2}\langle \j+\k,\w+\v-\j-\k \rangle-\frac{1}{2}.
\end{gather*}
\end{Lemma}

\begin{Remark}
Recalling the definitions of ${f}_{M\oplus N}$ and ${f}_{L}$, we can rewrite the above identity as
\begin{equation*}
\begin{aligned}
&\langle \n-(\w-\k),\m-\j\rangle+{f}_{M\oplus N} ((0, M_0\oplus N_0), \sigma_{M\oplus N})
=\langle \k,\v-\j\rangle+{f}_{L} ((0, L_0), \sigma_{L})-\frac{1}{2}.
\end{aligned}
\end{equation*}
\end{Remark}

Finally, we can simplify the calculation of $\tilde{{f}}_{\epsilon}*_{\epsilon} \tilde{X}_L$.
\begin{Lemma}\label{rhs2}
We have $ \tilde{{f}}_{\epsilon}*_{\epsilon} \tilde{X}_L =q^{\frac{1}{2}\lambda(\operatorname{ind}M,\operatorname{ind}N)-\frac{1}{2}}\cdot\tilde{X}_L.$
\end{Lemma}

\begin{proof}
A direct calculation shows
\begin{gather*}
 \tilde{{f}}_{\epsilon}*_{\epsilon} \tilde{X}_L \\
\quad=\int_{\g}\int_{\e,\f}\int_{(M_0,N_0)\in{\rm Gr}_{\e,\f}^{\epsilon}(FM,FN,\g)}q^{\frac{1}{2}\lambda(\operatorname{ind}M,\operatorname{ind}N)+\langle \f,\m-\e\rangle+
{f}_{M\oplus N} ((0, M_0\oplus N_0), \sigma_{M\oplus N})}\cdot X^{p(L,\g)}\\
\quad=q^{\frac{1}{2}\lambda(\operatorname{ind}M,\operatorname{ind}N)}\\
\phantom{\quad=}{}\cdot\int_{\g}\int_{\e,\f}\int_{(M_0,N_0)\in{\rm Gr}_{\e,\f}^{\epsilon}(FM,FN,\g)}q^{\langle \f,\m-\e\rangle
+ {f}_{M\oplus N} ((0, M_0\oplus N_0), \sigma_{M\oplus N}) }\cdot X^{p(L,\g)}\\
\quad=q^{\frac{1}{2}\lambda(\operatorname{ind}M,\operatorname{ind}N)}\\
\phantom{\quad=}{}\cdot\int_\g\int_{\k,\j}\int_{(K,J)\in{{\rm Gr}^{\mathbf{0}_{WV}}_{\k,\j}(K,J,\g)}}q^{\langle \n-(\w-\k),\m-\j\rangle+{f}_{M\oplus N} ((0, M_0\oplus N_0), \sigma_{M\oplus N})}\cdot X^{p(L,g)}\\
\quad=q^{\frac{1}{2}\lambda(\operatorname{ind}M,\operatorname{ind}N)-\frac{1}{2}}\cdot\int_\g\int_{\k,\j}\int_{(K,J)\in{{\rm Gr}^{\mathbf{0}_{WV}}_{\k,\j}(K,J,\g)}}q^{\langle \k,\v-\j\rangle+
{f}_{L} ((0, L_0), \sigma_{L})}\cdot X^{p(L,\g)}.
\end{gather*}

Now we focus on the exponent $\langle \k,\v-\j\rangle$ of $q$ in the above integration. Since $M$ and $N$ are rigid, so is $W\oplus V\oplus \Sigma^{-1}I$. Notice that
$
0\rightarrow V\rightarrow W\oplus V\rightarrow W\rightarrow 0$
is a split short exact sequence. Therefore in this case, we also have
$
\dim_k \Ext_{\mathcal{A}}^1(K,V/J)=0
$
for any $(K,J)\in{\rm Gr}_\k W\times{\rm Gr}_\j V$ and hence
$
\langle \k,\v-\j\rangle=\dim_k \Hom_{\mathcal{A}}(K,V/J).
$
On the other hand, we can consider the mapping
\[
\coprod_\g{\rm Gr}_\g FL\longrightarrow\coprod_{\k,\j}{\rm Gr}_\k W\times{\rm Gr}_\j V
\]
with affine fibers, induced by this split exact sequence, and observe that $\dim_k \Hom_{\mathcal{A}}(K,V/J)$ is exactly the dimension of the fiber at $(K,J)$.
 Thus we have
\begin{gather*}
 \tilde{{f}}_{\epsilon}*_{\epsilon} \tilde{X}_L \\
\quad=q^{\frac{1}{2}\lambda(\operatorname{ind}M,\operatorname{ind}N)-\frac{1}{2}}\cdot\int_\g\int_{\k,\j}\int_{(K,J)\in{{\rm Gr}^{\mathbf{0}_{WV}}_{\k,\j}(K,J,\g)}}q^{\langle \k,\v-\j\rangle+
{f}_{L} ((0, L_0), \sigma_{L})}\cdot X^{p(L,\g)}\\
\quad=q^{\frac{1}{2}\lambda(\operatorname{ind}M,\operatorname{ind}N)-\frac{1}{2}}\cdot\int_\g\int_{\k,\j}\int_{(K,J)\in{{\rm Gr}^{\mathbf{0}_{WV}}_{\k,\j}(K,J,\g)}}q^{\dim_k\Hom_{\mathcal{A}}(K,V/J)+
{f}_{L} ((0, L_0), \sigma_{L})}\cdot X^{p(L,\g)}\\
\quad=q^{\frac{1}{2}\lambda(\operatorname{ind}M,\operatorname{ind}N)-\frac{1}{2}}\cdot\int_\g\int_{L_0\in{\rm Gr}_\g FL}q^{
{f}_{L} ((0, L_0), \sigma_{L})}\cdot X^{p(L,\g)}\\
\quad=q^{\frac{1}{2}\lambda(\operatorname{ind}M,\operatorname{ind}N)-\frac{1}{2}}\cdot\int_\g |{\rm Gr}_\g FL|\cdot q^{-\frac{1}{2}\langle \g,\l-\g\rangle}\cdot X^{p(L,\g)}\\
\quad=q^{\frac{1}{2}\lambda(\operatorname{ind}M,\operatorname{ind}N)-\frac{1}{2}}\cdot\tilde{X}_L. \tag*{\qed}
\end{gather*} \renewcommand{\qed}{}
\end{proof}

By Theorem \ref{onedimtilde}, Lemmas~\ref{rhs1} and~\ref{rhs2}, we have
\begin{Theorem}[{\cite[Proposition 5.4.1]{fanqin}}]\label{maintheorem4}
In the cluster category $\mathcal{C}=D^b(\widetilde{\mathcal{A}})/{\tau^{-1}\Sigma}$ of a hereditary algebra $\widetilde{A}$, given two indecomposable coefficient-free rigid
objects $M, N \in \widetilde{\mathcal{A}}$ with \[\dim_k\Hom_{\mathcal{C}}(M,\Sigma N)=1\] and two non-split triangles
\[
N\rightarrow L\rightarrow M\xrightarrow{\epsilon}\Sigma N \qquad\text{and}\qquad
M\rightarrow L'\rightarrow N\xrightarrow{\eta}\Sigma M,
\]
where $L'$ is located in the fundamental domain, then we have
\[
\tilde{X}_M\cdot\tilde{X}_N=q^{\frac{1}{2}\lambda(\operatorname{ind}M,\operatorname{ind}N)-\frac{1}{2}}\cdot\tilde{X}_L+q^{\frac{1}{2}\lambda(\operatorname{ind}M,\operatorname{ind}N)}\cdot\tilde{X}_{L'}.
\]
\end{Theorem}

\subsection*{Acknowledgements}
The research was supported by the National Natural Science Foundation of China (no.~11771217 and no.~12031007). We greatly appreciate the referees' extraordinarily useful and detailed comments and suggestions which helped us to improve our manuscript. We are grateful to Xueqing Chen and Ming Ding for indicating many mistakes in the preliminary version of this paper and many valuable suggestions.

\pdfbookmark[1]{References}{ref}
\LastPageEnding

 \end{document}